\documentclass{article}
\usepackage[english]{babel}
\usepackage[latin1]{inputenc}
\usepackage[T1]{fontenc}
\usepackage{amsmath, amsfonts, amssymb, amsthm}
\usepackage[round, authoryear]{natbib}
\usepackage{multicol}

\newtheorem{defi}{Definition}
\newtheorem{lemm}{Lemma}[section]
\newtheorem{rema}{Remark}[section]
\newtheorem{prop}{Proposition}[section]
\newtheorem{theo}{Theorem}[section]
\newtheorem{nota}{Notation}
\newtheorem{exam}{Example}
\newcommand \N {\mathbb{N}}
\newcommand \Pn {\mathrm{P}_n}
\newcommand \Pk {\mathrm{P}_k}
\newcommand \und {\underline}

\newcommand \mcP {\mathcal{P}}
\newcommand \mbP {\mathbf{P}}
\newcommand \mbE {\mathbf{E}}
\newcommand \mbbP {\mathbb{P}}
\newcommand \mbbE {\mathbb{E}}
\newcommand \ind {\mathbb{I}}

\newcommand \vt {\vartheta}

\begin{document}

\begin{center}
\textbf{COALESCENT PROCESSES IN SUBDIVIDED POPULATIONS SUBJECT TO RECURRENT MASS 
EXTINCTIONS}
\end{center}
\bigskip

\begin{multicols}{2}
\begin{center}
\textbf{Jesse E. Taylor}\footnote{Supported by EPSRC grant EP/E010989/1.}\\
\medskip
Department of Statistics, University of Oxford\\
1 South Parks Road\\
Oxford OX1 3TG, United Kingdom\\
\verb"jtaylor@stats.ox.ac.uk"
\end{center}

\begin{center}
\textbf{Amandine V\'eber}\footnote{Corresponding author.}\\
\medskip
D\'epartement de Math\'ematiques,\\
Universit\'e Paris-Sud \\
91405 Orsay Cedex, France\\
\verb"amandine.veber@math.u-psud.fr"
\end{center}
\end{multicols}
\bigskip
\begin{abstract}
We investigate the infinitely many demes limit of the genealogy of 
a sample of individuals from a subdivided population that experiences
sporadic mass extinction events. By exploiting a separation of
time scales that occurs within a class of structured population
models generalizing Wright's island model, we show that as the 
number of demes tends to infinity, the limiting form of the genealogy
can be described in terms of the alternation of instantaneous \emph{scattering} 
phases that depend mainly on local demographic processes, and 
extended \emph{collecting} phases that are dominated by global processes.  
When extinction and recolonization events are local, the genealogy
is described by Kingman's coalescent, and the scattering phase influences
only the overall rate of the process.  In contrast, if the demes left 
vacant by a mass extinction event are recolonized by individuals 
emerging from a small number of demes, then the limiting genealogy 
is a coalescent process with simultaneous multiple mergers (a
$\Xi$-coalescent).  In this case, the details of the within-deme
population dynamics influence not only the overall rate of the
coalescent process, but also the statistics of the complex mergers
that can occur within sample genealogies.  These results suggest 
that the combined effects of geography and disturbance could
play an important role in producing the unusual patterns of genetic 
variation documented in some marine organisms with high
fecundity.


\end{abstract}

\section{Introduction}\label{section intro}
In this article, we investigate a class of population genetics models
that describe a population of individuals subdivided into $D$ demes 
which are subject to sporadic mass extinction events.  In general, we 
will think of these demes as corresponding to geographically distinct
subpopulations such as occur in Wright's island model \citep{WRI1931}, 
but this structure could also arise in other ways, such as through the 
association of homologous chromosomes within different individuals 
of a diploid species.  Whatever the source of the structure, many species 
are subject to recurrent disturbances which, if severe enough, can 
result in the extinction of a large proportion of the population \citep{SOU1984}.  
Important sources of widespread disturbance include fire, severe storms,
drought, volcanic eruptions, earthquakes, insect outbreaks, and
disease epidemics.  Our goal in this paper is to characterize the
effects that such events have on the genealogy of a sample of 
individuals or genes collected from the entire population.  Specifically, 
we will identify a set of conditions which will guarantee that in the limit
of infinitely many demes, the genealogy of the sample converges to
a process which alternates between two phases: an extended phase 
during which ancestral lineages occupy distinct demes, and an effectively 
instantaneous phase that begins each time two or more lineages are 
gathered into the same deme and ends when these are again scattered 
into different subpopulations through a combination of mergers and 
migrations.  The existence of this limit is a consequence of the separation 
of time scales between demographic events occurring within individual 
demes and those affecting the global dynamics of the population.

This study was partly motivated by recent investigations of the population 
genetics of several marine organisms whose genealogies appear to 
depart significantly from Kingman's coalescent (see Section \ref{intro-coalescent}). 
Based on their analysis of sequence polymorphism in a population of the 
Pacific oyster, \citet{EW2006} and more recently \citet{SW2008} suggest 
that the genealogies of some marine organisms with high fecundity and 
sweepstakes recruitment may be better described by a class of coalescent
processes that generalize Kingman's coalescent by allowing for 
simultaneous multiple mergers.  Indeed, in such organisms, the capacity 
of individuals to spawn millions of offspring makes it possible, in theory 
at least, for a substantial fraction of the population to be descended 
from a single parent.  However, depending on the life history and 
ecology of the species in question, this could happen in several 
different ways.  One possibility is that on rare occasions, individuals 
give birth to such a large number of offspring that even with random, 
independent survival of young, these cohorts constitute a sizable 
proportion of the next generation.  Such a scenario has been studied 
by \citet{SCH2003}, who showed that coalescents with multiple 
and simultaneous mergers arise naturally when the offspring 
distribution has a polynomial tail.  Another possibility is that on rare 
occasions a small number of individuals could contribute disproportionately
many of the surviving offspring not because they are exceptionally fecund, 
but because of mass reproductive failure or death in other parts of the
population.  What distinguishes these two scenarios is whether individuals 
win the recruitment sweepstakes by producing an exceptionally large 
number of offspring relative to the long-term average, or by simply giving 
birth to an average (or even below-average) number of offspring at a time 
when most other individuals experience an exceptional failure of 
reproduction.

The multiple mergers that occur in the models investigated in this
paper arise through a combination of both of these factors: mass
extinctions create large swathes of unoccupied territory which is
then instantaneously repopulated by individuals emerging from
finitely many demes. Of course, one weakness of this study is that
we do not identify the biological mechanism responsible for
restricting recolonization in this way, and in fact it seems difficult to
formulate such a mechanism that is both realistic and consistent with 
the metapopulation models considered in this paper. However, there 
are several scenarios under which similar dynamics could arise in 
a spatially extended population in which disturbance events tend to 
affect contiguous demes. For example, if dispersal distances are 
short, then recolonization of vacant habitat in a one-dimensional
population such as along a shoreline or a riparian corridor could be
dominated by individuals recently descended from the small number of
demes bordering the affected area.  Similar reasoning might also apply
to organisms with fractal-like distributions, such as aquatic or littoral
species in estuarine environments or possibly even HIV-1 populations
in the lymphatic system of an infected host. Alternatively, if regrowth 
from the margins is slow or even impossible (e.g., because
surviving demes are separated from vacant demes by inhospitable
habitat), then a few long-distance migrants could be responsible for
repopulating empty demes even in species with two-dimensional
distributions. Furthermore, in this case, we might also predict that
the number of demes contributing recolonizers would be negatively
correlated with the fecundity of the organism, since less time would
be available for additional migrants to enter the affected region
before the first migrant propagule had completely
repopulated the region. Although the mathematical analysis is much
more challenging than that given here, spatially-explicit models
incorporating these features are currently under development
(Alison Etheridge, pers.\ comm.).

\subsection{Wright's island model with mass extinctions}\label{first example}
To motivate both the class of models studied in this paper as well as the 
separation of time scales phenomenon that leads to the infinitely-many 
demes limit, let us begin by considering a version of Wright's island model 
with mass extinctions.  Suppose that a population of haploid organisms is 
subdivided into $D$ demes, each of which contains $N$ individuals.  We 
will assume that individuals reproduce continuously, i.e., generations are 
overlapping, and that at rate $1$ each individual gives birth to a single 
offspring which settles in that same deme with probability $1-m$ and
otherwise migrates to one of the other $D-1$ demes, chosen uniformly
at random.  In either case, we will assume that the deme size is constant 
and that a newborn individual immediately replaces one of the existing
$N$ members of the deme in which it settles.  Notice that if $m = 0$, then 
this model reduces to a collection of $D$ independent Moran models in 
populations of constant size $N$, whereas if $N = m = 1$, it describes the 
usual Moran model in a single population of size $D$.  However, in the 
following discussion we will assume that $m > 0$ and that $D$ is very much 
larger than $N$.

Before we account for mass extinctions, let us consider the genealogy of a sample 
of $n$ individuals chosen uniformly at random from the entire population.  We
first observe that, looking backwards in time, each lineage migrates out of its 
current deme at rate $(D-1)Nm/ND \approx m$.  Furthermore, if two lineages 
occupy different demes, then for these to coalesce, one of the two must migrate 
into the deme where the other lineage currently resides, an event that occurs 
approximately at rate $m/D$; here we have neglected terms of order $D^{-2}$
and will continue to do so without further comment.  When two lineages are
collected in the same deme, then they can either coalesce immediately, which 
happens with probability $1/N$, or they can cohabit within that deme for some random 
period of time until either they coalesce or they migrate into different demes.  
Since two lineages occupying the same deme coalesce at rate $2(1-m)/N$, 
and each lineage, independently of the other, migrates out of the deme at rate 
$m$, the probability that the two coalesce rather than migrate is $\chi = 
(1-m)/(1 - m + Nm)$.  Putting these observations together, it follows that every 
time two lineages are collected within the same deme by migration, the total probability 
that they coalesce rather than migrate into different demes is $1/N + (1-1/N)\chi = 
1/(1 - m + Nm)$, and the time that elapses between entry into the same deme
and either coalescence or escape is a mixture of a point mass at $0$ (in case
they coalesce at the entry time) and an exponential random variable with 
mean $N/(2mN + 2(1-m))$.  In particular, notice that typically much less time 
is required for two lineages occupying the same deme to either coalesce 
or escape (of order $N$) than is needed for two lineages occupying different 
demes to be collected into the same deme (of order $D$).  It is this disparity
between the rate of events happening within individual demes and the rate
at which lineages are gathered together that gives rise to a separation of 
time scales in the island model.  If we rescale time by a factor of $D$ 
and let the number of demes tend to infinity, then the time required for two 
lineages sampled from different demes to coalesce is exponentially distributed 
with mean $(1 - m + Nm)/2m$.  

To complete our description of the coalescent process in this model, we need  
to consider the possibility of more complex coalescent events.  We first observe
that if $n$ individuals are sampled from $D$ demes, then the probability that
all of these individuals reside in different demes will be close to one if $D$ is 
much greater than $n$.  Furthermore, because lineages occupying different 
demes coalesce and migrate independently of one another, it is straightforward 
to show that the probability that three or more lineages ancestral to our sample 
are collected into the same deme is of order $D^{-2}$ or smaller.  Likewise, it 
can be shown that the probability of having multiple pairs of lineages collected 
into several demes at the same time is similarly negligible.  From these 
observations, it follows that only pairwise coalescence events matter in the 
infinitely-many demes limit, and that if there are $n$ ancestral lineages, then at 
rate ${n \choose 2} 2m/(1 - m + Nm)$, two of these, chosen uniformly at random, 
coalesce, leaving $n-1$ ancestral lineages.  In other words, the genealogy for 
this model can be approximated by a scalar time change of Kingman's coalescent, 
with a rate that depends on both the migration rate and the deme size.  This
result is essentially due to \cite{WA2001}, who considered a similar model 
with non-overlapping generations and Wright-Fisher sampling.

Now let us introduce mass extinction events into this model.  Fix $e > 0$ 
and $y \in [0,1)$, and suppose that at rate $e/D$, the metapopulation suffers 
a disturbance which causes each deme to go extinct, independently of all 
others, with probability $y$.   For example, we could consider a model in which 
the demes represent small islands or keys in the Caribbean and the disturbances 
are hurricanes that completely inundate those islands lying in their path.  Here 
we are reverting to the original time units, i.e., time has not yet been rescaled
by a factor of $D$, and we have chosen the disturbance rate so that mass 
extinctions occur at rates commensurate with coalescence in the pure island
model.  In keeping with the assumption that deme size is constant, we will 
assume that all of the islands that are left vacant by a mass extinction are 
immediately recolonized by offspring dispersing out of a single source deme 
that is chosen uniformly at random from among the demes unaffected by the 
disturbance.  In addition, we will assume that the parent of each colonizing 
individual is chosen uniformly at random from the $N$ members of the source 
deme.  Of course, the entire metapopulation could be extirpated by a mass 
extinction if $y > 0$, but the probability of this outcome is exponentially small 
in $D$ and can be disregarded as $D$ tends to infinity. 

Suppose that a mass extinction occurs at a time when there are $n$ ancestral 
lineages occupying distinct demes.  Bearing in mind that we are now looking 
backwards in time, all of the lineages belonging to demes that are affected 
by the disturbance will move into the source deme, where those sharing the 
same parent will immediately coalesce.  Thus, one reason that multiple mergers 
can occur in this model is because of the very highly skewed distribution of 
recolonizing offspring contributed both by individuals and demes following a 
mass extinction.  Suppose that there are $n_{1}$ distinct lineages remaining in 
the source deme once we account for this initial set of coalescences.  These 
lineages will undergo a random sequence of migration and coalescence 
events until there is only one lineage remaining within the source deme.  For 
example, if $n_{1} = 4$, then one possible outcome would see one lineage
migrate out of the deme followed by a pair of binary mergers, leaving 
only one lineage within the source deme.  Whatever the sequence, the 
amount of time required to scatter the lineages into different demes will be 
of order $O(1)$, whereas the time until either the next mass extinction 
event or the next binary merger involving lineages outside of the source
deme will be of order $O(D)$.  Thus, if we again rescale time by a factor of 
$D$, then any sequence of coalescence and migration events involving a 
source deme will effectively be instantaneous when we let $D$ tend to infinity.  
This is the second way in which multiple merger events can arise in this model.  
Furthermore, varying the migration rate and deme size changes not only the 
overall rate of coalescence, but also the relative rates of the different kinds of 
multiple merger events that can occur.  For example, if $Nm$ is very small, 
then the coalescent process will be close to a $\Lambda$-coalescent (which 
has multiple mergers, but not simultaneous multiple mergers) because most 
lineages that are collected into a source deme by a mass extinction event will 
coalesce before any escape by migration.  However, as $N$ increases, so will 
the probability that multiple lineages enter into and then escape from the source 
deme without coalescing.  This suggests that at moderate values of $Nm$,
mass extinctions may be likely to result in simultaneous mergers (i.e., the coalescent 
is a $\Xi$-coalescent), while for very large values of $Nm$, multiple mergers 
of all types will be unlikely and the coalescent process will tend towards Kingman's 
coalescent.

\subsection{Neutral genealogies and coalescents}\label{intro-coalescent}
In the last twenty years, coalescent processes have taken on increasingly
important role in both theoretical and applied population genetics, where 
their relationship to genealogical trees has made them powerful tools 
to study the evolution of genetic diversity within a population.  Under the 
assumption of neutrality, allelic types do not influence the reproduction of 
individuals and it is therefore possible to separate `type' and `descent'. 
This allows us to study the genealogy of a sample of individuals on its 
own and then superimpose a mechanism describing how types are transmitted 
from parent to offspring, justifying the interest in investigating genealogical
processes corresponding to particular reproduction mechanisms without
explicit mention of types. We refer to \citet{NOR2001} for a review
of coalescent theory in population genetics.

Beginning with the coalescent process introduced by \citet{KIN1982}
to model the genealogy of a sample of individuals from a large
population, three increasingly general classes of coalescent
processes have been described. A key feature shared by all three classes
is the following consistency property: the process induced on the
set of all partitions of $\{1,\ldots,n\}$ by the coalescent acting
on the partitions of $\{1,\ldots,n+k\}$ (obtained by considering
only the blocks containing elements of $\{1,\ldots,n\}$) has the
same law as the coalescent acting on the partitions of
$\{1,\ldots,n\}$. In terms of genealogies, this property means that
the genealogy of $n$ individuals does not depend on the size of the
sample that contains them. To describe these continuous-time Markov
processes, it will be convenient to introduce some notation. For all 
$n\in\N$, we denote the set of all partitions of $[n]\equiv \{1,\ldots,n\}$ 
by $\Pn$. In the following, the index $n$ of the set of partitions in which 
we are working will be referred to as the \textbf{sample size}, an
element of $\{1,\ldots,n\}$ will be called an \textbf{individual},
and `block' or `\textbf{lineage}' will be equivalent terminology to
refer to an equivalence class. If $\zeta\in \bigcup_n\Pn$, then
$|\zeta|=k$ means that the partition $\zeta$ has $k$ blocks. Also,
for $\zeta,\eta \in \Pn$ and $k_1,\ldots,k_r \geq 2$, we will write $\eta
\subset_{k_1,\ldots,k_r} \zeta$ if $\eta$ is obtained from $\zeta$
by merging exactly $k_1$ blocks of $\zeta$ into one block, $k_2$
into another block, and so on. Kingman's coalescent is defined on
$\Pn$ for all $n\geq 1$, as a Markov process with the following
$Q$-matrix: if $\zeta,\eta \in \Pn$,
$$q_K(\zeta\rightarrow \eta) = \left\{\begin{array}{cl}1 & \mathrm{if\ }\eta \subset_2
\zeta, \\
-\binom{|\zeta|}{2} & \mathrm{if\ }\eta=\zeta,\\
0 & \mathrm{otherwise}.\end{array}\right.$$ A more general class of
exchangeable coalescents, allowing mergers of more than two blocks
at a time, was studied by \citet{PIT1999} and \citet{SAG1999}. These
coalescents with multiple mergers (or $\Lambda$-coalescents) are in
one-to-one correspondence with the finite measures on $[0,1]$ in the
following manner: for a given coalescent, there exists a unique
finite measure $\Lambda$ on $[0,1]$ such that the entries
$q_{\Lambda}(\zeta \rightarrow \eta)$ of the $Q$-matrix of the
coalescent, for $\zeta,\eta \in \Pn$, are given by
$$q_{\Lambda}(\zeta\rightarrow \eta) = \left\{\begin{array}{ll}\int_0^1
\Lambda(dx)x^{k-2}(1-x)^{b-k} & \mathrm{if\ }\eta \subset_k
\zeta \mathrm{\ and\ }|\zeta|=b, \\
- \int_0^1 \Lambda(dx)x^{-2} \big(1-(1-x)^{b-1}(1-x+bx)\big)&
\mathrm{if\ }\eta=\zeta \mathrm{\ and\ } |\zeta|=b,\\
0 & \mathrm{otherwise}.\end{array}\right.$$ Kingman's coalescent is
recovered by taking $\Lambda=\delta_0$, the point mass at $0$.
Lastly, a third and wider class of coalescents was introduced by
\citet{MS2001} and \citet{SCH2000}, for which mergers involving more
than one ancestor are allowed. These coalescents with simultaneous
multiple mergers (or $\Xi$-coalescents) are characterized in
\citet{SCH2000} by a finite Borel measure on the infinite ordered simplex
$$\Delta=\Big\{(x_1,x_2,\ldots): x_1\geq x_2\geq \ldots \geq 0,\sum_{i=1}^{\infty}x_i\leq 1\Big\}.$$
Indeed, to each coalescent corresponds a unique finite measure $\Xi$
on $\Delta$ of the form $\Xi=\Xi_0+a\delta_{\mathbf{0}}$, where
$\Xi_0$ has no atom at zero and $a\in [0,\infty)$, such that the
transition rates of the coalescent acting on $\Pn$ are given by
$$q_{\Xi}(\zeta\rightarrow \eta) = \int_{\Delta} \frac{\Xi_0(d\mathbf{x})}
{\sum_{j=1}^{\infty}x_j^2}\bigg(\sum_{l=0}^s \sum_{i_1\neq
\ldots\neq i_{r+l}} \binom{s}{l}x_{i_1}^{k_1}\ldots
x_{i_r}^{k_r}x_{i_{r+1}} \ldots x_{i_{r+l}} \big(1-\sum_{j=1}^\infty
x_j \big)^{s-l}\bigg) + a\ \ind_{\{r=1,k_1=2\}}$$ if $\eta
\subset_{k_1,\ldots,k_r} \zeta$ and $s \equiv |\zeta|-\sum_{i=1}^r
k_i$. The other rates (for $\eta \neq \zeta$) are equal to zero. The
$\Lambda$-coalescents are particular cases of $\Xi$-coalescents, for
which $\Xi(x_2>0)=0$.

As mentioned above, coalescent processes can be used to describe the
genealogy of large populations. Indeed, a large body of literature
has been devoted to describing conditions on the demography
of a population of finite size $N$ that guarantee that the genealogical 
process of a sample of individuals converges to a coalescent as $N$ 
tends to infinity. Such limiting results for populations with discrete 
non-overlapping generations are reviewed in \citet{MOH2000}, and 
some examples can be found for instance in \citet{SCH2003}, 
\citet{EW2006} and \cite{SW2008}.  In these examples, the shape
of the limiting coalescent is related to the propensity of individuals 
to produce a non-negligible fraction of the population in the
next generation.

However, the representation of the genealogy as a coalescent
requires in particular that any pair of lineages has the same chance
to coalesce.  This condition breaks down when the population is
structured into subpopulations, since then coalescence will occur
disproportionately often between lineages belonging to the same deme.  
To model these kinds of scenarios, structured analogues of coalescent 
processes were introduced \citep[see e.g.][]{NOT1990, WIL1998},
which allow lineages both to move between demes as well as coalesce
within demes.  Various state spaces have been used to describe a structured
coalescent, such as vectors in which the $i$'th component gives the
lineages (or their number) present in deme $i$, or vectors of pairs
`block $\times$ deme label'.  All these representations of a
structured genealogy take into account the fact that the
reproductive or dispersal dynamics may differ between demes, hence
the need to keep track of the location of the lineages. In contrast, 
several papers investigate models where the
structure of the genealogy collapses on an appropriate time scale,
i.e., the limiting genealogy no longer sees the geographical division
of the population. In \cite{COX1989}, demes are located at the sites
of the torus $\mathbb{T}(D)\subset \mathbb{Z}^d$ of size $D$ and
each site can contain at most one lineage. Lineages move between
sites according to a simple random walk, and when one of them lands
on a site already occupied, it merges instantaneously with the
inhabitant of this `deme'. These coalescing random walks, dual to the
voter model on the torus, are proved to converge to Kingman's
coalescent as $D\rightarrow \infty$.  More precisely, Cox shows that
if $n<\infty$ lineages start from $n$ sites independently and
uniformly distributed over $\mathbb{T}(D)$, then the process counting 
the number of distinct lineages converges to the pure death process 
that describes the number of lineages in Kingman's coalescent. This 
analysis is generalized in \cite{CD2002} and \cite{ZCD2005}, where each 
site of the torus now contains $N\in \N$ individuals and a Moran-type
reproduction dynamics occurs within each deme. Again, the limiting
genealogy of a finite number of particles sampled at distant sites
is given by Kingman's coalescent, and convergence is in the same
sense as for Cox' result.  Other studies of systems of particles
moving between discrete subpopulations and coalescing do not 
require that the initial locations of the lineages be thinned out. In
\cite{GLW2007}, demes are distributed over the grid $\mathbb{Z}^2$
and the process starts with a Poisson-distributed number of lineages
on each site of a large box of size $D^{\alpha/2}$, for some $\alpha\in (0,1]$. 
The authors show that the total number of lineages alive at times of the 
form $D^t$ converges in distribution as a process (indexed by $t\geq \alpha$)
to a time-change of the block counting process of Kingman's coalescent, 
as $D\rightarrow \infty$. See \cite{GLW2007} for many other references 
related to these ideas.

Our emphasis in this paper will be on the separation of time scales 
phenomenon and the way in which local and global demographic
processes jointly determine the statistics of the limiting coalescent
process.  Consequently, we shall always assume that the demes comprising 
our population are exchangeable, i.e., the same demographic processes
operate within each deme, and migrants are equally likely to come 
from any one of the $D$ demes.  In this simplified setting, we only need 
to know how lineages are grouped into demes, but not the labels of
these demes.

\subsection{Separation of time scales}
A separation of time scales can be said to occur whenever 
different components of a stochastic process evolve at rates which
greatly differ in their magnitudes.  This concept is usually invoked
when there is a sequence of stochastic processes $(X^D_t,t\geq 0)$ 
on a space $E$ as well as a function $\eta:E\rightarrow E'$
and an increasing sequence $r_D \rightarrow \infty$ such that the
processes $(X^D_{r_D^{-1}t},t\geq 0)$ have a non-trivial limit
$(X^{\mathrm{fast}}_t,t\geq 0)$ determined by the \textbf{fast}
time scale, while the processes $(\eta(X^D_t),t\geq 0)$ (which are
only weakly influenced by the fast evolution) have another
non-trivial limit $(X^{\mathrm{slow}}_t,t\geq 0)$ determined by the
\textbf{slow} time scale.  It is often the case that the processes 
$(X^{D}_{t}, t \geq 0)$ have the Markov property but do not 
converge to a limit, while the slow processes $(\eta(X^D_t),t\geq 0)$
do converge, but are not Markovian.

Separation of time scales techniques were first introduced into
population genetics by \citet{EN1980}, and since then have been used
to study the genealogical processes of structured populations in
several different settings. \citet{NK2002} consider a population of
total size $N$, evolving according to a Wright-Fisher model
\citep[see][]{FIS1930, WRI1931} and distributed over $D<\infty$
demes. These demes are in turn structured into groups of demes,
within which individuals migrate faster (at a rate of order $N^{-\alpha}$ 
for an $\alpha \in [0,1]$) than from one group to another (which occurs 
at a rate $O(N^{-1})$). When all demes are connected by fast migration, 
they show that structured genealogy collapses to an unstructured 
Kingman's coalescent as $N$ tends to infinity, due to the fact that
migration is so fast compared to the coalescence rate (of order
$N^{-1}$) that the population becomes well-mixed before the first
coalescence event occurs. When several groups of demes are 
connected by slow migration, the genealogical process converges to 
a structured coalescent, in which groups of demes act as panmictic
populations and coalescence of lineages within a group is faster
than between two groups. These results are made possible by the fact
that the blocks of the partition induced by the genealogy are not
affected by a migration event. Since only migration occurs on the
fast time scale  and coalescence is on the slow time scale, forgetting
about the location of the lineages gives a sequence of (non-Markov)
processes which converge on the slow time scale to a Markov process.

Another kind of separation of time scales was studied by Wakeley and
co-authors in a series of papers \citep[see in
particular][]{WAK1998,WAK1999,WAK2004,WA2001}. In these models, a
population evolving in discrete non-overlapping generations occupies
$D$ demes, labeled $1,\ldots,D$. Deme $i$ contains a population of
$N_i$ adults and receives $M_i$ migrants each generation. Then, a
Wright-Fisher resampling within each deme brings the population
sizes back to their initial values. Other mechanisms can also be
taken into account, such as extinction of a group of demes followed
by instantaneous recolonization. Allowing $D$ to tend to infinity
greatly simplifies the analysis of the genealogical processes, and
in particular gives rise to a decomposition of the genealogy of a
sample of individuals into two different phases, occurring on two
time scales. Following the terminology introduced in \citet{WAK1999},
the first phase to occur is the \textbf{scattering phase}, in which
lineages occupying the same deme coalesce or move to \textbf{empty}
demes (`empty' meaning that none of the sampled lineages are in this
deme). In the limit, this phase occurs on the fast time scale  and is
therefore viewed as instantaneous. At the end of the scattering
phase, all remaining lineages lie in different demes. The
\textbf{collecting phase} is the following period of time during which
lineages are gathered together into the same demes by migration or
extinction/recolonization, where they may merge. The limiting
genealogical process is a coalescent on the slow time scale, which
ends when the number of lineages reaches one.

As we have already mentioned, apart from an initial instantaneous burst 
of mergers (which only occurs if multiple individuals are sampled from 
the same deme) all of the genealogical processes obtained in this setting 
are scalar time changes of Kingman's coalescent. Indeed, in the forwards 
in time evolution, migrants and colonizers are assumed to come from the 
whole population or from a non-vanishing fraction of the demes and so, 
with probability one, only two of the finitely many lineages of the sample 
are brought into the same deme at a time in the limit. Subsequently, the 
two lineages either coalesce or are scattered again, but in any case the 
outcome is at most a binary merger. In this paper, we shall study 
coalescent processes that arise in population models which
include mass extinctions and general recolonization mechanisms, and
describe the conditions in which it corresponds to an unstructured
$\Xi$-coalescent on the slow time scale. To this end, we will speak
of `scattering' and `collecting' phases in a more general sense.  We
prefer to call the `collecting phase' the period of time during which
lineages wander among empty demes until a migration or extinction
event brings several lineages into the same deme. We shall show
that, once such a `geographical collision' has occurred, an
instantaneous scattering phase follows at the end of which all
lineages have merged or moved to empty demes.  Another collecting
phase then starts and so on until the most recent common ancestor
of the sample has been reached and there is only one lineage 
remaining.

\subsection{Framework and main results}\label{section resultats}
Fix $n\in \N$ and consider the genealogy of a sample of $n$
individuals from a population of $D>n$ demes (the following
framework also allows $D= \infty$). In the following, we shall
suppose that demes are exchangeable in the sense given in Section
\ref{intro-coalescent}. We shall work in the space $\Pn^s$ defined as
follows:
\begin{defi}\label{structured partitions}
Let $\tilde{\mathrm{P}}_n^s$ be the set \setlength \arraycolsep{1pt}
\begin{eqnarray*}
\tilde{\mathrm{P}}_n^s &\equiv &
\big\{\big(\{B_1,\ldots,B_{i_1}\},\ldots,
\{B_{i_{n-1}+1},\ldots,B_{i_n}\}\big):\ 0\leq i_1 \leq \ldots\leq
i_n \leq n ,\\ & & \qquad \qquad \emptyset\neq B_j\subset [n]\
\forall j\in \{1,\ldots,i_n\},\ \{B_1,\ldots, B_{i_n}\}\in \Pn\big\}
\end{eqnarray*}
 of
$n$-tuples of sets (we allow some of the components of the $n$-tuple
to be empty), and let us define the equivalence relation $\sim$ on
$\tilde{\mathrm{P}}_n^s$ by $\xi \sim \xi'$ if and only if there
exists a permutation $\sigma$ of $[n]$ such that, if
$\mathcal{B}_1\equiv\{B_1,\ldots,B_{i_1}\}, \ldots,\
\mathcal{B}_n\equiv \{B_{i_{n-1}+1},\ldots,B_{i_n}\}$ are the
components of the vector $\xi$, then $\xi'=
\big(\mathcal{B}_{\sigma(1)},\ldots,\mathcal{B}_{\sigma(n)}\big)$.
The quotient of $\tilde{\mathrm{P}}_n^s$ by $\sim$ is denoted by
$\Pn^s$.

We call any
$\big(\{B_1,\ldots,B_{i_1}\},$ $\ldots,\{B_{i_{n-1}+1},\ldots,
B_{i_n}\}\big) \in \Pn^s$ an \textbf{unordered structured partition of $[n]$}.
\end{defi}
In view of the application we have in mind, each component
$\mathcal{B}_j$ represents a particular deme containing some of the
lineages ancestral to the sample, and the blocks $B_k$ (for $k\in
\{1,\ldots,i_n\}$) specify the partition of the sample 
determined by the ancestors alive at a particular time.  Empty 
components are used to guarantee a constant vector size, $n$, 
independent of the index $D$ used later.  In the following, we 
omit the term `unordered' when referring to the structured partitions 
of Definition~\ref{structured partitions}.

The finite set $\Pn^s$ is endowed with the discrete topology, which
is equivalent to the quotient by $\sim$ of the discrete topology on
$\tilde{\mathrm{P}}_n^s$.
\begin{defi}
A Markov process $\mcP$ on $\Pn^s$ for which blocks can only merge
and change component is called a \textbf{structured genealogical
process}.
\end{defi}
To illustrate the possible transitions, let us take $n=5$ and consider the 
following sequence of events:
\begin{eqnarray*}\big(\{\{1\}\},\{\{2\}\},\{\{3\}\},\{\{4\}\},
\{\{5\}\}\big)& \stackrel{(i)}{\rightarrow} &
\big(\{\{1\},\{2\}\},\{\{3,4\}\},
\{\{5\}\},\emptyset, \emptyset\big)\\
&\stackrel{(ii)}{\rightarrow} &
\big(\{\{1\}\},\{\{2\}\},\{\{3,4,5\}\},
\emptyset, \emptyset\big)\\
&\stackrel{(iii)}{\rightarrow} & \big(\{\{1,2,3,4,5\}\},
\emptyset,\emptyset,\emptyset, \emptyset\big).
\end{eqnarray*}
In this example, we start from the configuration in $\mathrm{P}_5^s$
where each lineage is alone in its deme. During transition $(i)$,
either $\{1\}$ or $\{2\}$ changes component and both blocks end up in
the same deme (which creates an empty component in our
representation), but remain distinct. In contrast, either $\{3\}$ 
or $\{4\}$ also moves (emptying another component), but then
the two blocks merge into a single block $\{3,4\}$ which is not
allowed to split during later transitions. Block $\{5\}$ remains
alone in its component. During transition $(ii)$, lineages $\{1\}$
and $\{2\}$ are scattered again into two different demes by the
movement of one of them, while one of the lineages $\{3,4\}$ or
$\{5\}$ changes component and the two blocks merge. Eventually, all
the remaining blocks are gathered into the same deme and merge into
a single block. Since elements of $\mathrm{P}_5^s$ are defined up to
a permutation of their components and since a block is not allowed
to split, no other change is possible from the state reached after
transition $(iii)$.

\begin{rema}\label{remark consistency}
Movements and mergers of blocks do not alter the sample
size. However, this does not guarantee that the structured
genealogies are consistent in the sense given in Section
\ref{intro-coalescent} as we would expect from a reasonable
genealogical process. In fact, several conditions will be imposed on
the models we consider so that this property holds: see Lemma
\ref{consistence} for the consistency of the fast genealogical
process, and the set of conditions (\ref{consist lambda g}) imposed
on the geographical gatherings in Proposition \ref{xi-coal}.
Proposition \ref{collapse} states in particular that the latter
conditions are necessary and sufficient for the genealogies to be
consistent on both time scales and that when they are fulfilled, the
unstructured genealogical process on the slow time scale is 
a $\Xi$-coalescent.
\end{rema}

Let us order the components of a given structured partition by the
smallest element belonging to a block contained in the component (if
it is non-empty). Empty components come last. For each $k \leq n$
and $\zeta\in \Pn^s$, let us write $|\zeta|_a=k$ if the $a$'th
component (in the order just defined) of the structured partition
$\zeta$ contains $k$ blocks, and define a subset $\Pi_n$ of $\Pn^s$
by
\begin{equation}\label{definition
Pi}\Pi_n\equiv \big\{\zeta \in \Pn^s: |\zeta|_a \leq 1 \quad
\forall\ a \in \{1,\ldots,n\}\big\}.\end{equation}
$\Pi_n$ is the
set of all structured partitions of $[n]$ in which each deme
contains at most one lineage. These sets will appear naturally in
the description of the limiting processes.

Recall from the example given in Section \ref{first example} that 
the rate at which lineages are collected together in the same deme
is much smaller than the rate at which lineages already occupying
the same deme either coalesce or are scattered into different demes.
Furthermore, as in that example, we will continue to assume that catastrophic 
extinction-recolonization events occur rarely, in fact, at rates that
are of the same order of magnitude as the rate at which lineages occupying
different demes are brought together by ordinary migration.  With 
these points in mind, let us consider a sequence $(\mcP^D_s,s\geq 0)$ 
of structured genealogical processes for a finite sample from the 
whole population, which consists of the following kinds of events:
\begin{itemize}
\item within-deme coalescence and movement of lineages to empty demes
at rates of order $O(1)$;
\item movement of groups of lineages initially occupying different demes 
into the same deme, possibly followed by mergers of some of these lineages, 
at  rates of order $O(r_D^{-1})$. 
\end{itemize}
Let us rescale time by a factor of $r_{D}$ so that the coalescence rate of 
two individuals in different demes is of order $O(1)$ as $D$ tends to infinity. 
Of course, within-deme coalescence and migration now occur at increasing 
rates of order $O(r_D)$. This implies that, for a given
sample size $n$, the generator $G^D$ of the genealogical process
acting on $\Pn^s$ has the form
$$G^D=r_D \Psi +\Gamma + R_D,$$
where $\Psi,\Gamma$ and $R_D$ are bounded linear operators,
$\langle R_D\rangle \rightarrow 0$ as $D\rightarrow \infty$, and we do not 
record the dependence of the operators on the sample size $n$. Here, 
if $\|\cdot\|$ stands for the supremum norm on the space of functions $f :
\Pn^s\rightarrow \mathbb{R}$, then $\langle R\rangle$ is defined by
\begin{equation}\label{norm}
\langle R\rangle=\sup_{f\neq 0}\frac{\|Rf\|}{\|f\|}.
\end{equation}

Because $r_D\rightarrow\infty$, the sequence $(G^D)_{D\geq 1}$ is
unbounded, even when applied to functions of the unstructured
partition induced by $\mcP^D$, and so we do not expect the structured 
coalescent processes corresponding to these generators to converge 
pathwise. Nevertheless, our heuristic description of the fast dynamics 
suggests that elements of $\Pi_n$ will be unaffected by the `fast' events
corresponding to $\Psi$, which will indeed be the case under the assumptions  
made in Section \ref{section construction}. 
Furthermore, we will show (cf. Lemma \ref{description xi}) that 
the process generated by $\Psi$ on $\Pn^s$ and starting at $\zeta \in \Pn^s$ 
a.s.\ reaches a random final state $\und{\zeta}$ in $\Pi_{n}$ in a finite 
number of steps.  Since the rates of the events generated by $\Psi$ 
grow to infinity, increasing numbers of these events take place before the 
first event corresponding to $\Gamma$ even occurs.  This motivates the 
description of the genealogy given above in terms of an alternation of 
very short scattering phases driven by $\Psi$ and of longer collecting 
phases ending with the first event generated by $\Gamma$ at which 
$\mcP^D$ leaves $\Pi_n$.  Viewing all of the transitions occurring
during a given scattering phase as a single, more complex event, 
and exploiting the fact that these phases are vanishingly short, 
it is plausible that there is a genealogical process $\mcP$ with 
values in $\Pi_n$ such that for each fixed time $t>0$, $\mcP^D_t 
\Rightarrow \mcP_t$ as $D\rightarrow \infty$. Our main result makes 
these heuristic arguments rigorous:
\begin{theo}\label{convergence p-dim}Let $\zeta\in \Pn^s$.
Under the conditions described in Section \ref{section conditions},
the finite-dimensional distributions of the structured genealogical
process $\mcP^D$ starting at $\zeta$ converge to those of a
$\Pi_n$-valued Markov process $\mcP$ starting at $\und{\zeta}$,
except at time $0$.
\end{theo}

The proof that $\mcP^D$ converges in law to $\mcP$ in the Skorokhod space
$D_{\Pn^s}[0,\infty)$ of all c\`adl\`ag paths with values in $\Pn^s$
requires tightness of the corresponding sequence of distributions. 
We shall show in Proposition \ref{prop non-tight} that this property holds 
if and only if the rate at which the genealogical process leaves the 
set $\Pi_n$ tends to zero as $D$ grows to infinity.  Indeed, if this
condition is not satisfied, then two or more jumps can accumulate 
during a scattering phase: the jump out of $\Pi_n$ followed by the 
events needed to bring $\mcP^D$ back into $\Pi_n$.  Fortunately, the 
proof that the unstructured genealogical processes are tight is less 
demanding, since these processes do not change state when lineages 
move between demes.  In this case, an accumulation of jumps due to 
the fast within-deme dynamics will be ruled out if we can show that the 
probability that the process $\mcP^D$ re-enters $\Pi_n$ in a single 
jump converges to one as $D$ tends to infinity.

The limiting process $\mcP$ with values in $\Pi_n$ is introduced 
and investigated in Section \ref{section construction}, and we show in
Proposition \ref{xi-coal} that, under the assumptions of Theorem 1.1, 
the unstructured genealogical process induced by $\mcP$ is the restriction 
to $\Pn$ of a $\Xi$-coalescent.  We also identify the limiting process 
$\xi$ for the genealogy on the fast time scale in Section \ref{section
construction}, and state in Proposition \ref{convergence fast
time scale} the convergence of $\mcP^D_{r_D^{-1}\cdot}$ to $\xi$ 
as processes with values in $D_{\Pn^s}[0,\infty)$.
The proofs of Theorem \ref{convergence p-dim} and Proposition
\ref{convergence fast time scale} are given in Section \ref{section
convergence}, along with a discussion of the tightness of $\mcP^D$.
Although the conditions of Theorem \ref{convergence p-dim} are somewhat contrived, 
we show in Section \ref{section collapse} that these are necessary and sufficient for the 
unstructured genealogical process of a generalized island model to converge to a
$\Xi$-coalescent on the slow time scale.  In Section \ref{section example}, we apply
these results to a particular class of models incorporating mass extinction events.
Based on our analysis of this class, we suggest that families of $\Xi$-coalescents 
may often interpolate between $\Lambda$-coalescents and Kingman's coalescent in 
structured population models, and that it may be a generic property of such models
that they admit simultaneous mergers whenever they admit multiple mergers.

\section{Construction of the limiting genealogical processes}\label{section construction}
\subsection{A generalized Island-Cannings Model}\label{section
cannings}
To motivate the genealogical processes considered in this paper, we
begin by introducing a general model for the demography of a
subdivided population which combines features of the Cannings model
\citep{CAN1974} with those of the classical Island model
\citep{WRI1931}.

Suppose that the population is subdivided into $D$ demes, each of
which contains $N$ haploid individuals. Islands are labeled
$1,\dots,D$, while individuals within each island are labeled
$1,\ldots,N$. At rate $1$, an $ND^2$-dimensional random vector
$R\equiv \big(R^{i,j}_k,\ i,j\in\{1,\ldots,D\},
k\in\{1,\ldots,N\}\big)$ is chosen, such that for all $i,j,k$,
$R^{i,j}_k$ is the number of descendants of the $k$'th individual in
deme $j$ which settle into deme $i$ during the event. In keeping
with the spirit of the Cannings' model, we use the term `descendant'
both to refer to the offspring of reproducing individuals as well as
to individuals which were alive both before and after the event (as
in Cannings' formulation of the Moran model). We impose the
following conditions on the random variables $R^{i,j}_k$:
\begin{enumerate}
\item \textbf{Constant deme size:} With probability $1$, for all $i\in [D]$
we have $\sum_{j,k}R^{i,j}_k =N$.
\item \label{forwards_exch}\textbf{Exchangeable dynamics:}
The law of $R$ is invariant under any permutation $\sigma$ of
$[D]^2\times [N]$ such that for every $i\in [D]$,
$\sigma(i,i,k)_1=\sigma(i,i,k)_2$, i.e., $\sigma$ conserves the
relation \emph{source deme} $=$ \emph{destination deme}. (Here,
$\sigma(i,j,k)_l$ denotes the $l$'th component of the permuted
vector.)
\end{enumerate}
Then, in each deme the current population is replaced by the $N$
offspring coming into this deme during the event, which we label in
an exchangeable manner.

Let us comment on the above conditions. The first one simply
guarantees that the number of individuals in each deme is constant
and equal to $N$. For the second condition, let us first fix $i,\ j$
and a permutation $\tau$ of $[N]$, and look at the permutation
$\sigma$ given by $\sigma(i,j,k)=(i,j,\tau(k))$ and
$\sigma(i',j',k')=(i',j',k')$ whenever $i\neq i'$ or $j\neq j'$.
Then, condition \ref{forwards_exch} corresponds to the
exchangeability of the contribution of the inhabitants of deme $j$
in repopulating deme $i$. Second, fix $i$ and choose a permutation
$\tau$ of $[D]\setminus\{i\}$. Set $\sigma(i,j,k)=(i,\tau(j),k)$ if
$j\neq i$, $\sigma(i,i,k)=(i,i,k)$ and $\sigma(i',j',k')=(i',j',k')$
whenever $i'\neq i$. In this case, condition \ref{forwards_exch}
states that the demes different from deme $i$ contribute in an
exchangeable manner to the repopulation of deme $i$. Finally, let
$\tau$ be a permutation of $[D]$ and define
$\sigma(i,j,k)=(\tau(i),j,k)$ if $j\notin\{i,\tau(i)\}$,
$\sigma(i,j,k)=(\tau(i),\tau(i),l)$ if $j=i$, and
$\sigma(i,j,k)=(\tau(i),i,l)$ if $j=\tau(i)$. For such permutations,
condition \ref{forwards_exch} asserts that the dispersal mechanism
is exchangeable with respect to the destination of dispersing
individuals (provided that this differs from the source deme).
Overall, our assumptions aim at making the dynamics depend on the
labels as weakly as possible, but we allow the repopulation
mechanism of a deme to differ according to whether the new
individuals are produced within this deme or come from one of the
$D-1$ other demes.

\begin{exam}\label{ex0}If $R$ is invariant under all permutations
$\sigma$ of $[D]^2\times [N]$ (not just those
satisfying condition \ref{forwards_exch}), then the dynamics are
those of a Cannings' model for a panmictic population of size $DN$,
i.e., there is no population subdivision.
\end{exam}

\begin{exam}\label{ex1} If all demes evolve independently of each other,
then $R^{i,j}\equiv 0$ whenever $j\neq i$.
Condition~\ref{forwards_exch} imposes that $(R^{i,i},i\in [D])$
should be an exchangeable $D$-tuple of exchangeable $N$-tuples, a
situation corresponding to a continuous-time Cannings model acting
within each deme.
\end{exam}

\begin{exam}\label{ex2}
Let $m\in [0,1]$ and assume that, with probability $1-m$, $R$ is
chosen as in Example \ref{ex1}. With probability $m$, four numbers
$i,j,l,k$ are sampled uniformly at random in $[D]^2\times [N]^2$,
and the $k$'th individual in deme $j$ produces an offspring that
replaces the $l$'th individual in deme $i$. In this case,
$R^{i,i}=(1,\ldots,0,\ldots,1)$, where the unique zero is in the
$l$'th coordinate; $R^{i,j}=(0,\ldots,1,\ldots,0)$, where the unique
$1$ is in the $k$'th coordinate; $R^{i,j'}= (0,\ldots,0)$ if
$j'\notin \{i,j\}$ and for $i'\neq i$, $R^{i',j'}= (1,\ldots,1)$ if
$i'=j'$ and $(0,\ldots,0)$ otherwise. This model gives a simple
example including within-deme reproduction and individual migration.
Alternatively, individuals could be exchanged between demes during a
migration event, in which case a descendant of individual $l$ in
deme $i$ (in the above notation) also replaces individual $k$ in
deme $j$.
\end{exam}

\begin{exam}\label{ex3}An event during which one deme goes extinct
and is recolonized by the offspring
of individuals coming from other demes has the following formulation: 
$R^{i,i}= (0,\ldots,0)$ if deme $i$ goes extinct, $R^{l,l}=
(1,\ldots,1)$ if $l\neq i$ and the repopulation of deme $i$
satisfies the exchangeability condition \ref{forwards_exch}. For
instance, $N$ individuals are chosen uniformly at random among the
$N(D-1)$ inhabitants of the other demes and contribute one offspring
in the new population of deme $i$.
\end{exam}

Many other kinds of events can be imagined, but these three
mechanisms (reproduction, migration and extinction/recolonization)
will be the building blocks of the models we shall consider in this
paper. Viewed backwards in time, reproduction events as in Example
\ref{ex1} will correspond to the merger of several lineages if they
are produced (forwards in time) by the same individual during the
event considered. A migration event such as in Example \ref{ex2}
will correspond to the movement of one or a few lineages from their
demes to other subpopulations, if these lineages happen to have
their parents in the source demes. An extinction event will also
typically result in the movement of lineages among demes, and could
involve much larger numbers of individuals or demes than simple
migration events. Note that lineages can both move and merge during
the same event, if their common parent lies in a different deme.

\subsection{Genealogy on the fast time scale} \label{paragraphe xi}
Let us start by constructing a structured genealogical process
$(\xi_t,t\geq 0)$ such that its restriction to $\Pn^s$ describes the
genealogy of $n$ individuals on the fast time scale  of individual
demes. This process will incorporate mergers of lineages occupying
the same deme as well as dispersal of lineages into empty demes
(i.e., those not containing other ancestral lineages), but no events
where geographically separated lineages end up in identical demes
and possibly merge. In fact, if the rate at which such events occur
is very large, then it is not difficult to see that the structure of
the population effectively disappears on the fast time scale  and the
model collapses to that of a panmictic population. We thus rule out
this kind of situations to keep a structured population.

We construct the process $\xi$ by specifying its restriction to
$\Pn^s$. As $\Pn^s$ is a finite set, we can define a continuous-time
Markov process on this space by specifying its transition rates.
Because a block represents a single ancestor, whose descendance at
time $0$ is made of the individuals contained in the block, we shall
ask that the rates at which blocks move and merge do not depend on
the number or labels of these individuals. Hence, these rates will
only depend on the collection $\{k_1,\ldots,k_p,0,\ldots,0\}$ giving
the numbers of blocks contained in the different components of
$\xi$. In order to describe the possible transitions, we need the
following definition.
\begin{defi}\label{defi compatible}Let $\hat{k}=\{k_1,\ldots,k_p\}$
and $\hat{k}'=\{k_1',\ldots,k_q'\}$ be two collections of (non-zero)
integers. We shall write $\hat{k} \rhd\hat{k}'$ if $q\geq p$,
$\sum_{i=1}^q k_i'\leq \sum_{j=1}^p k_j$, and we can arrange the
elements of $\hat{k}'$ so that for each $i\in\{1,\ldots,p\}$, we
have $1\leq k_i'\leq k_i$ and at least one of such $k_i'$ is
strictly less than $k_i$.
\end{defi}
\noindent Note that no collection $\hat{k}$ of integers satisfies
$\{1,\ldots,1\}\rhd \hat{k}$.

For all pairs $(\hat{k}, \hat{k}')$ such that $\hat{k}\rhd\hat{k}'$,
let $\vt_{\hat{k},\hat{k}'}\in \mathbf{R}_+$. In addition, if
$\zeta\in \Pn^s$, let $\hat{k}(\zeta)$ be the collection of integers
which gives the number of blocks within each non-empty component of
$\zeta$. Define the infinitesimal rate $q_{(\xi)}(\eta|\zeta)$ of a
particular transition $\zeta\rightarrow \eta$ (when $\eta \neq
\zeta$ and both belong to $\Pn^s$) as:
\begin{itemize}
\item $q_{(\xi)}(\eta|\zeta)=\vt_{\hat{k}(\zeta),\hat{k}(\eta)}$,
if $\eta$ can be obtained from $\zeta$ by first merging some
number (possibly zero) of blocks contained in the same component 
of $\zeta$, and then moving some blocks to formerly empty demes 
with the restriction that only blocks originating from the same deme 
can be gathered into the same destination deme (again, we
allow the number of blocks moved to be zero).  In this case, we
easily see that we must have $\hat{k}(\zeta)\rhd\hat{k}(\eta)$.
\item $q_{(\xi)}(\eta|\zeta)=0$ otherwise.
\end{itemize}

In the following, we shall assume that for any $\zeta\in\Pn^s$
containing more than one block in at least one component, the rates
satisfy the condition
$$\sum_{\eta\in \Pn^s}q_{(\xi)}(\eta|\zeta)>0.$$
These conditions ensure that, whenever a deme contains more than one
lineage, a scattering or a coalescence event will happen in the
future with probability one. Recall the definition of $\Pi_n$ given
in (\ref{definition Pi}). From the form of the rates given above, we
see that any $\eta \in \Pi_n$ is an absorbing state for $\xi$.
Moreover, we have the following result, saying in essence that the
process $\xi$ with values in $\Pn^s$ reaches a final state in a
finite number of steps, and this final state is a random variable
with values in $\Pi_n$.

\begin{lemm} \label{description xi}Let $\tau_{\pi}$ be the stopping
time defined by $\tau_{\pi}\equiv \inf\{t\geq 0: \xi_t \in \Pi_n\}.$
Then, $\tau_{\pi}$ is a.s. finite and for all $t\geq \tau_{\pi}$,
$\xi_t = \xi_{\tau_{\pi}}$.
\end{lemm}

\begin{proof} From our assumptions on the rates, the only absorbing states
of the process $\xi$ are the structured partitions contained in
$\Pi_n$. Moreover, any transition results in a coarsening of the
corresponding unstructured partition or in the movement of some
lineages to different empty demes, so the number of transitions
for $\xi$, starting at any $\xi_0 \in \Pn^s$ is bounded by $n$.
Since the time between two events is exponentially distributed with
a non-zero parameter (the sum of the rates of the possible
transitions) as long as the process has not reached an absorbing
state, the finiteness of the number of transitions undergone by
$\xi$ imposes that $\tau_{\pi}$ is a.s. finite.
\end{proof}

Let us introduce the following notation, justified by the result of
Lemma \ref{description xi}.
\begin{nota}\label{notation und}If $\zeta \in
\Pn^s$, let $\und{\zeta}$ denote a random variable with values in
$\Pi_n$, whose distribution is that of the final state of the
structured genealogical process $\xi$ started at $\zeta$.
\end{nota}
We end this subsection with the following lemma, whose main purpose
is to introduce the notion of consistency for structured
genealogical processes. If $\zeta \in \Pk^s$ and $\tilde{\zeta}\in
\mathrm{P}_{k+1}^s$ for some $k\in \N$, let us write $\zeta \prec
\tilde{\zeta}$ if the projection of $\tilde{\zeta}$ onto $\Pk^s$
(the $k$-tuple describing the structured partition of $1, \ldots,k$)
equals $\zeta$.
\begin{lemm}\label{consistence}Suppose that $\xi$ is defined on $\Pk^s$ for
every $k\in \N$. The following conditions are equivalent:

(i) For each $k\geq 1$, $\zeta,\eta \in \Pk^s$ and
$\tilde{\zeta}\in \mathrm{P}_{k+1}^s$ such that $\zeta \prec
\tilde{\zeta}$,
$$q_{(\xi)}(\eta|\zeta)= \sum_{\tilde{\eta}}q_{(\xi)}(\tilde{\eta}|\tilde{\zeta}),$$
where the sum is over all $\tilde{\eta}\in \mathrm{P}_{k+1}^s$ such
that $\eta \prec \tilde{\eta}$.

(ii) The process $\xi$ is consistent in the sense that, for all
$k\geq 1$, if $\zeta\in \Pk^s$ and $\zeta'\in \mathrm{P}_{k+1}^s$
satisfy $\zeta\prec \zeta'$, then the law of the restriction to
$\Pk^s$ of the process $\xi$ started at $\zeta'$ is the same as the
law of $\xi$ started at $\zeta$.

In particular, if both conditions are fulfilled and if $\eta\in
\Pi_k$ has $r$ blocks, then
$$\mbP\big[\ \und{\zeta}=\eta\big]= \sum_{j=1}^{r+1}\mbP\big[\ \und{\zeta'}=\eta^{(j)}\big],$$
where for each $j\in \{1,\ldots,r\}$, $\eta^{(j)}\in \Pi_{k+1}$ is
obtained from $\eta$ by adding an empty $(k+1)$-st component to
$\eta$ (which turns it into an element $\eta'$ of $\Pi_{k+1}$), and
adding $k+1$ in the $j$'th block of $\eta'$. Likewise,
$\eta^{(r+1)}$ is obtained by adding the singleton $\{k+1\}$ in an
empty component of $\eta'$.
\end{lemm}
\begin{proof} Let $\xi_{(k)}$ (resp. $\xi_{(k+1)}$) denote the process
$\xi$ started at $\zeta\in \Pk^s$ (resp. $\zeta'\in
\mathrm{P}_{k+1}^s$), and call $\xi'_{(k)}$ the projection onto
$\Pk^s$ of $\xi_{(k+1)}$. Since we work with finite state spaces and
discrete jump processes, $(ii)$ is equivalent to the fact that for
all $\gamma,\eta \in \mathrm{P}_{k}^s $ the rate at which
$\xi'_{(k)}$ jumps from $\gamma$ to $\eta$ is equal to the
corresponding transition rate for $\xi_{(k)}$. By construction, the
former is the sum of the rates of all the transitions from the
current state of $\xi_{(k+1)}$ to a state $\eta'$ such that $\eta
\prec \eta'$, and so $(ii)$ holds if and only if $(i)$ does.

The second part of Lemma \ref{consistence} is a direct consequence
of the consistency of the process.
\end{proof}

\subsection{Limiting process on the slow time scale}\label{descript
P}

Let us now describe the form that we would expect the genealogical
process to take on the slow time scale  as the number of demes tends
to infinity. This process $\mcP$ will have values in $\Pi_n$, so
once again we construct it by specifying its transition rates.

Recall the two ingredients of the description of the structured
genealogical processes indexed by $D<\infty$, given in Section
\ref{section resultats}. Coalescence and movement of blocks to
formerly empty demes are the two kinds of events that constitute
the fast process $\xi$, and we saw in Lemma \ref{description xi} that the
final state reached by $\xi$ belongs to $\Pi_n$ a.s. Therefore, we
now need to describe how the resulting geographically separated
lineages are gathered into identical demes and, potentially, merge
during the same event. As in the definition of $\xi$, the rate at
which such an event occurs will only depend on the number $r$ of
demes containing at least one lineage just after the event, on the
numbers $k_1,\ldots,k_r$ of blocks brought together into these
components, and on the number and sizes of the groups of blocks
ending up in the same demes which subsequently merge into a bigger
block. Hence, we shall use the following terminology.
\begin{defi}\label{defi geogr collision} Let $k\geq 2$, and $k_1,\ldots,k_r \geq 1$
such that $\sum_{i=1}^rk_i =k$ and at least one of the $k_i$'s is
greater than $1$. Let also $L_1= \{l_{1,1}, \ldots,l_{1,i_1}\}$,
$\ldots\ $, $L_r=\{l_{r,1},\ldots,l_{r,i_r}\}$ be $r$ (unordered)
sets of integers such that for $j\in \{1,\ldots,r\}$, we have
$\sum_{u=1}^{i_j}l_{j,u}=k_j$. We call an event in which $k$
lineages spread in $k$ different demes become grouped into $k_1$
lineages in one deme, $k_2$ lineages in another deme, $\ldots\ $,
and for all $j\in\{1,\ldots,r\}$, $l_{j,1}$ lineages in deme $j$
merge into one, $l_{j,2}$ into another, and so on (all mergers occur
between lineages which landed in the same deme) a
\textbf{$(k;k_1,\ldots,k_r;L_1,\ldots,L_r)$-geographical collision}.
\end{defi}
\begin{rema} A geographical collision is to be understood as a particular
transition. Because the order of $k_1,\ldots,k_r$ does not matter, a
$(k;k_1,\ldots,k_r;L_1,\ldots,L_r)$-geographical collision is also a
$(k;k_{\sigma(1)},\ldots,k_{\sigma(r)};$
$L_{\sigma(1)},\ldots,L_{\sigma(r)})$-geographical collision for any
permutation $\sigma$ of $\{1,\ldots,r\}$. Thus, for a given
$(k;k_1,\ldots,k_r;$ $L_1,\ldots,L_r)$, the number of
$(k;k_1,\ldots,k_r;L_1,\ldots,L_r)$-geographical collisions is
$$\mathrm{A}(k;k_1,\ldots,k_r)\prod_{m=1}^r \mathrm{A}(k_m;l_{m,1},\ldots,l_{m,i_m}),$$
where if $k,k_1,\ldots,k_r$ are such that $\sum_{i=1}^r k_i=k$ and
$b_j$ is the number of $k_i$'s equal to $j$, then
\begin{equation}\label{number geo_coll}\mathrm{A}(k;k_1,\ldots,k_r)=
\binom{k}{k_1,\ldots,k_r}\frac{1}{\prod_{j=1}^k b_j!}.\end{equation}
Indeed, the binomial term gives the number of ways of choosing $k_1$
blocks to form a family numbered $1$, $k_2$ other blocks to form
family $n^o2$, and so on. But any permutation of the labels of families
having the same size gives the same unordered structured partition,
hence the normalization by the fraction in the right-hand side
of (\ref{number geo_coll}).
\end{rema}
Let us now define the structured genealogical process $\mcP$. The
relation between the coefficients of $\mcP$ and the sequence of
structured genealogical processes will be given in the next section,
and we simply give the form of the limiting process here. Recall that
$|\zeta|_a = k$ if the $a$'th component of $\zeta$ contains $k$ blocks, 
and write $|\zeta|$ for the total number of blocks of $\zeta\in\Pn^s$, 
that is $|\zeta|=\sum_{a=1}^n|\zeta|_a$.  Furthermore, let $\und{\zeta}$ 
be a $\Pi_{n}$-valued random variable with the distribution specified in 
Notation \ref{notation und}.
\begin{defi}\label{definition P}
For all integers and sets $k,\ k_i$ and $L_i$ satisfying the
conditions of Definition \ref{defi geogr collision}, let
$\lambda^g_{k;k_1,\ldots,k_r;L_1,\ldots,L_r}\geq 0$. Then,
$(\mcP_t,t\geq 0)$ is the Markov process with values in $\Pi_n$
which evolves as follows: when $\mcP_t = \chi \in \Pi_n$, any
$(|\chi|;k_1,\ldots,k_r; L_1,\ldots,L_r)$-geographical collision
occurs at rate $\lambda^g_{|\chi|;k_1,\ldots,k_r;L_1,\ldots,L_r}$.
Given that this collision has outcome $\zeta \in \Pn^s$, the new
value of $\mcP$ is drawn from the distribution of $\und{\zeta}$.
\end{defi}
We can recover the expression for the rate of any given transition in the form
$$q(\eta|\chi)= \sum_{\zeta\in \Pn^s} \lambda^g_{|\chi|;k_1,\ldots,k_r;
L_1,\ldots,L_r} \und{\zeta}[\eta],$$ where the rate
$\lambda^g_{|\chi|;k_1,\ldots,k_r; L_1,\ldots,L_r}$ in the term of
the sum labeled by a given $\zeta$ is the rate of occurrence of the
only possible geographical collision turning $\chi$ into $\zeta$, if
such a collision exists. If it does not, we set the rate to $0$.
Consequently, the previous description does specify a Markov process
on $\Pi_n$.

Observe that this description allows `ghost events' in which
lineages are gathered in identical demes by a geographical collision
and then scattered again in different demes without coalescing, so
that the actual transition is of the form $\chi\rightarrow \chi$.
However, we shall need to keep track of these ghost events in the
proof of convergence of the structured genealogical processes.
Therefore, we shall always consider them as events which do occur at
a certain rate but have no effect on $\mcP$.

To finish the description of our limiting process, we have the
following result, which in fact describes the unstructured genealogy
under some additional conditions.
\begin{prop}\label{xi-coal}For each $\zeta \in \Pi_n$,
let us define $\zeta^u \in \Pn$ as the unstructured partition of $n$
induced by $\zeta$. Then the unstructured genealogical process
$(\mcP^u_t,t\geq 0)$ induced by $\mcP$ is a Markov process with
values in $\Pn$. Suppose in addition that condition $(i)$ of Lemma
\ref{consistence} holds and that the $\lambda^g$'s satisfy the
following consistency equations: for all $k\in \N$ and compatible
$k_1,\ldots,k_r,L_1,\ldots,L_r$,
\begin{equation}\label{consist lambda g}\lambda^{g}_{k;k_1,\ldots,k_r;L_1,\ldots,L_r}=
\sum_{u=1}^r\sum_{j=1}^{i_u+1}\lambda^{g}_{k+1;k_1,\ldots,k_u+1,
\ldots,k_r;L_1,\ldots,L_u^{(j)},\ldots,L_r}+
\lambda^{g}_{k+1;k_1,\ldots,k_r,1;L_1,\ldots,L_r,\{1\}},
\end{equation}
 where for each $u\leq r$
$$L_u^{(j)}=\left\{ \begin{array}{ll}
\{l_{u,1},\ldots,(l_{u,j}+1),\ldots,l_{u,i_u}\} & \mathrm{if\ }j\leq
i_u \\ \{l_{u,1},\ldots,l_{u,i_u},1\} & \mathrm{if\ }j= i_u+1.
\end{array}\right. $$
(In particular, if instantaneous coalescence after the gathering of
lineages is forbidden, then the $\lambda^g$'s are associated to a
$\Xi$-coalescent.)

Then $(\mcP^u_t,t\geq 0)$ is the restriction to $\Pn$ of a
$\Xi$-coalescent on the partitions of $\N$.
\end{prop}
\begin{rema}\label{tilde lambda} By fixing $k,k_1,\ldots,k_r$ and summing
over all compatible integer sets $L_1,\ldots,L_r$, we see that the
rates $\tilde{\lambda}^g_{k;k_1,\ldots,k_r}$ at which $k$ lineages
lying in $k$ different demes end up in a configuration where $k_1$
lineages are in the same deme, $k_2$ in another one, and so on
(regardless of how many of them merge instantaneously thereafter), are
associated to a $\Xi$-coalescent whenever condition (\ref{consist
lambda g}) holds.
\end{rema}
\begin{proof}
Any component of $\gamma \in \Pi_n$ contains at most one block and
all $n$-tuples are defined up to a permutation of their components,
so the map $\Pi_n\rightarrow \Pn :\ \gamma \mapsto \gamma^u$ is a
measurable bijection between $\Pi_n$ and $\Pn$. Thus, $\mcP^u$
inherits the Markov property of $\mcP$, and its transition rates
$q^u(\eta^u|\ \gamma^u)$ are obviously given by $q^u(\eta^u|\
\gamma^u)= q(\eta|\ \gamma)$.

Let us turn to the second part of Proposition \ref{xi-coal}. By
assumption, $\mcP^u$ only coarsens as time goes on and it is easy to
check that all transition rates $\rho(k;k_1,\ldots,k_r)$ from a
partition with $k$ blocks to a partition obtained by merging $k_1$
of those blocks into one, $k_2$ into a second one, $\ldots\ $
($k_1,\ldots, k_r \in \N,\ \sum_{i=1}^r k_i=k$), are equal and
depend only on $k, k_1,\ldots,k_r$ (the order of $k_1,\ldots,k_r$
does not matter). Therefore, we need only check the consistency
condition given in \citet{SCH2000} to identify $\mcP^u$ as the
restriction to the partitions of $[n]$ of a $\Xi$-coalescent. As the
rates do not depend on $n$, let us rather work in $\Pi_{k}$ with
$\gamma= \big(\big\{\{1\}\big\},\ldots, \big\{\{k\}\big\}\big)$ and
$\eta=\big(\big\{\{1,\ldots,k_1\}\big\},\big\{\{k_1+1,\ldots,k_1+k_2\}\big\},\ldots,
\big\{\{k_1+\ldots+k_{r-1}+1,\ldots,k\}\big\},\emptyset,\ldots,\emptyset\big)$,
and check that
$$\rho(k;k_1,\ldots,k_r)= \sum_{i=1}^r \rho(k+1;k_1,\ldots,k_i+1,\ldots,k_r)
+ \rho(k+1;k_1,\ldots,k_r,1).$$
Since the $\lambda^g$'s satisfy (\ref{consist lambda g}), we have
\setlength\arraycolsep{1pt}
\begin{eqnarray}
\rho(k;k_1,\ldots,k_r) &=& q(\eta|\ \gamma) \nonumber\\
& = & \sum_{\zeta\in \Pk^s}
\lambda^g_{k;l_1,\ldots,l_s;L_1,\ldots,L_s}\
\und{\zeta}[\eta] \label{expression rho}\\
& = & \sum_{\zeta\in \Pk^s} \sum_{v=1}^s \sum_{j=1}^{i_v+1}
\lambda^g_{k+1;l_1,\ldots,l_v+1,\ldots,l_s;L_1,\ldots,L_v^{(j)},\ldots,L_s}\
\und{\zeta}[\eta]+\ \sum_{\zeta\in
\Pk^s}\lambda^g_{k+1;l_1,\ldots,l_s,1;L_1,\ldots,L_s,\{1\}}
\und{\zeta}[\eta]. \nonumber
\end{eqnarray}
We wish to compare this rate to the rates corresponding to $k+1$
blocks. To this end, let us define $\zeta_{(v,j)}\in
\mathrm{P}_{k+1}^s$ for all $\zeta\in \Pk^s$ with $l$ non-empty
components and $v\in \{1,\ldots,l+1\},\ j\in \{1,\ldots,i_v+1\}$
($i_v$ being the number of blocks in the $v$'th non-empty component
of $\zeta$) by turning $\zeta$ into a $(k+1)$-tuple and adding
individual $k+1$ in the $j$'th block of the $v$'th component of the
new vector ($v=l+1$ means that we add the block $\{k+1\}$ in the
extra component, and likewise $j=i_v+1$ means that we add the block
$\{k+1\}$ in the $v$'th component of the new vector). For example,
with the previous notation $\gamma$,
$$\gamma_{(1,2)}=\big(\big\{\{1\},\{k+1\}\big\},\ldots,\big\{\{k\}\big\},\emptyset
\big)\quad \mathrm{and}\quad
\gamma_{(k+1,1)}=\big(\big\{\{1\}\big\},\ldots,\big\{\{k\}\big\},
\big\{\{k+1\}\big\}\big).$$ Define also $\gamma^{(j)}\in \Pi_{k+1}$,
for all $\gamma\in \Pi_k$ with $r$ blocks and $j\in \{1,\ldots,r+1\}$,
by turning $\gamma$ into a $(k+1)$-tuple and adding individual $k+1$
in the block of the $j$'th component of the new vector. Once again,
$j=k+1$ means that we add a block $\{k+1\}$ in the extra component.
For instance,
$$\eta^{(1)}=\big(\big\{\{1,\ldots,k_1,k+1\}\big\},\big\{\{k_1+1,\ldots,k_1+k_2\}\big\},\ldots,
\big\{\{k_1+\ldots+k_{r-1}+1,\ldots,k\}\big\},\emptyset,\ldots,\emptyset\big)$$
and
$$\eta^{(k+1)}=\big(\big\{\{1,\ldots,k_1\}\big\},\big\{\{k_1+1,\ldots,k_1+k_2\}\big\},\ldots,
\big\{\{k_1+\ldots+k_{r-1}+1,\ldots,k\}\big\},\big\{\{k+1\}\big\},\emptyset,\ldots,\emptyset\big).$$
With this notation, we see that
\begin{equation}\label{sum rates}
\sum_{i=1}^r \rho(k+1;k_1,\ldots,k_i+1,\ldots,k_r)
+ \rho(k+1;k_1,\ldots,k_r,1)=
\sum_{i=1}^{r+1} q(\eta^{(i)}|\ \gamma_{(k+1,1)}).
\end{equation}
For all $\zeta'\in \mathrm{P}_{k+1}^s$, there exists a unique
triplet $(\zeta,v,j)$ where $\zeta \in \Pk^s$ has $l$ non-empty
components, $v \in \{1,\ldots,l+1\}$ and $j\in \{1,\ldots,i_v+1\}$
such that $\zeta'= \zeta_{(v,j)}$. Indeed, $\zeta$ is given by the
partition of $\{1,\ldots,k\}$ induced by $\zeta'$, $v$ is the
component containing $k+1$ and $j$ is the block of that component in
which $k+1$ lies. Therefore, the right-hand side of (\ref{sum
rates}) is equal to
\begin{equation}\label{sum rates 2}
\sum_{i=1}^{r+1}\sum_{\zeta\in
\Pk^s}\sum_{v=1}^s\sum_{j=1}^{i_v+1}\lambda^g_{k+1;l_1,\ldots,l_v+1,\ldots,l_s;
L_1,\ldots,L_v^{(j)},\ldots,L_s}\ \und{\zeta_{(v,j)}}[\eta^{(i)}]
+\sum_{i=1}^{r+1}\sum_{\zeta\in
\Pk^s}\lambda^g_{k+1;l_1,\ldots,l_s,1;L_1,\ldots,L_s,\{1\}}\
\und{\zeta_{(s+1,1)}}[\eta^{(i)}],
\end{equation}
where $s$ and the coefficients
$\lambda^g_{k+1;l_1,\ldots,l_v+1,\ldots,l_s;
L_1,\ldots,L_v^{(j)},\ldots,L_s}$ correspond to the particular
$\zeta$ indexing the term of the sum. Let us look at a particular
$\zeta$ in the second sum. The block $\{k+1\}$ remains a singleton
just after the geographical collision, so it is not affected by a
following genealogical event and
$\und{\zeta_{(s+1,1)}}[\eta^{(i)}]=0$ for all $i\in \{1,\ldots,r\}$.
Lemma \ref{consistence} hence implies that
$\und{\zeta_{(s+1,1)}}[\eta^{(r+1)}]=\und{\zeta}[\eta]$, and the
second term of (\ref{sum rates 2}) is equal to
$$\sum_{\zeta\in \Pk^s}\lambda^g_{k+1;l_1,\ldots,l_s,1;L_1,\ldots,L_s,\{1\}}\ \und{\zeta}[\eta].$$
Let us look at a particular $\zeta$ in the first sum, now. When
$v\leq s$, the corresponding geographical collision brings $k+1$ in
a block of the $v$'th component of $\zeta$. By the second part of
Lemma \ref{consistence}, the probability that the final state of all
the blocks different from $k+1$ is given by $\eta$ is equal to the
sum over all corresponding final states of these blocks and ${k+1}$.
But taking the sum over $i$ in $\sum_{i=1}^{s+1}
\und{\zeta_{(v,j)}}[\eta^{(i)}]$ boils down to considering all such
final states, since $\und{\zeta_{(v,j)}}[\eta^{(i)}]=0$ if the
individuals in the $i$'th block of $\eta$ were not in the $v$'th
component of $\zeta$ before their rearrangement by the genealogical
process (recall that, under the action of $\xi$, lineages can merge
only if they start in the same deme). Therefore, we obtain that, for
all $\zeta \in \Pk^s$ and compatible $v,j$,
$$\sum_{i=1}^{r+1}\und{\zeta_{(v,j)}}[\eta^{(i)}]= \und{\zeta}[\eta].$$
Coming back to expressions (\ref{sum rates}) and (\ref{sum rates 2}), we obtain
\setlength\arraycolsep{1pt}
\begin{eqnarray*}\sum_{i=1}^r &\rho&(k+1;k_1,\ldots,k_i+1,\ldots,k_r) + \rho(k+1;k_1,\ldots,k_r,1)\\
&=&\sum_{\zeta\in
\Pk^s}\sum_{v=1}^s\sum_{j=1}^{i_v+1}\lambda^g_{k+1;l_1,\ldots,l_v+1,\ldots,l_s;
L_1,\ldots,L_v^{(j)},\ldots,L_s}\ \und{\zeta}[\eta] +\sum_{\zeta\in \Pk^s}
\lambda^g_{k+1;l_1,\ldots,l_s,1;L_1,\ldots,L_s,\{1\}}\ \und{\zeta}[\eta]\\
&=& \rho(k;k_1,\ldots,k_r),
\end{eqnarray*}
where the last equality follows from (\ref{expression rho}). This
completes the proof of Proposition \ref{xi-coal}.
\end{proof}

\section{Convergence of the structured genealogical processes}\label{section convergence}
Now that we have constructed the potential limits for our sequence
of structured genealogical processes on the fast and slow
time scales, let us state precisely what conditions we impose and in
which sense $\mcP^D_{r_D^{-1}\cdot}$ and $\mcP^D$ converge.

\subsection{Description of the conditions}\label{section conditions}
Let $n\geq 1$ be the sample size and define two types of events:
\begin{itemize}
\item Type 1: some lineages contained in the same demes merge 
and some move (potentially in groups) to empty islands. The number 
of lineages involved in either step can be zero (meaning that only
coalescence or only scattering has occurred), and lineages starting from
different demes are not gathered into the same deme by the event.

\item Type 2: $k$ lineages move, but at least one of them lands in a
non-empty deme or at least two dispersing lineages not coming from
the same deme are gathered. During that event, $k_1$ lineages end up
in the same deme, $k_2$ lineages in another, and so on. This is
immediately followed by the coalescence of some lineages lying in
identical demes (the number of such mergers can be zero, meaning
that the lineages have only moved).
\end{itemize}
By our assumptions on the genealogical processes, these two types
describe all kinds of events which can happen to the structured
genealogical process $\mcP^D$, for each $D$. For conciseness, we
shall call an event of type $i$ an $i$-event. Assume now that, when
$\mcP^D$ has value $\zeta \in \Pn^s$ and $\eta \in \Pn^s$ is a
possible new value compatible with the type of the event ($o$'s hold
as $D$ goes to infinity):
\begin{enumerate}
\item \label{ass1} The rate of occurrence of a particular 1-event
$\zeta\rightarrow \eta$ can be written
$$r_D\vt_{\hat{k}(\zeta),\hat{k}(\eta)} + \nu^{(n)}(\zeta,\eta)+o(1)
\qquad \mathrm{as\ }D\rightarrow \infty,$$
where for each $n$, $\nu^{(n)}(\cdot,\cdot)$ is
a bounded function on $(\Pn^s)^2$ and $r_D\rightarrow \infty$ as
$D\rightarrow \infty$.
\item \label{ass2} Consider a 2-event involving $k$ lineages, for which
there exist $k_1,\ldots,k_r \geq 1$ such that $\sum_{i=1}^rk_i
=|\zeta|$ and there exist $r$ sets of integers
$L_1=\{l_{1,1},\ldots,l_{1,i_1}\},\ \ldots,\
L_r=\{l_{r,1},\ldots,l_{r,i_r}\}$ such that for all $j\in
\{1,\ldots,r\}$ we have $\sum_{u=1}^{i_j}l_{j,u}=k_j$, satisfying:
in the new structured partition, $k_1$ lineages end up in one deme,
$k_2$ in another deme, $\ldots\ $, and for all $j\in
\{1,\ldots,r\}$, $l_{j,1}$ lineages in deme $j$ merge into one,
$l_{j,i_2}$ into another one, and so on (once again, all mergers
occur between lineages lying in the same deme). Then the rate of
occurrence of any such event is of the form
$$\mathit{l}^{(n)}_k(\zeta,\eta) + o(1),$$
where for each $n$ and all $k\leq n$, $\mathit{l}^{(n)}_k$ is a
bounded function on $\Pn^s \times \Pn^s$, and in particular if $\zeta \in \Pi_n$,
$$\mathit{l}^{(n)}_k(\zeta,\eta)= \lambda^g_{|\zeta|;k_1,\ldots,k_r;L_1,\ldots,L_r}.$$
 Here again, the order of $k_1,\ldots,k_r$ does
not matter.
\item \label{ass4}The $\vt$'s correspond to a structured genealogical
process $\xi$ as described in the last section, and the
$\lambda^g$'s satisfy the consistency equations (\ref{consist lambda
g}).
\end{enumerate}
Morally, the coalescence of lineages occupying common demes and the
scattering of such lineages into empty demes occur more and more
rapidly as $D$ tends to infinity, whereas events collecting lineages
into common demes occur at bounded rates. Other events are less and
less frequent, so that in the limit we obtain a separation of
time scales between the instantaneous structured genealogical process
and the slow collecting phase of the limiting unstructured
genealogical process. Notice that 1-events do not affect a
structured partition contained in $\Pi_n$.

Let $G^{n,D}$ denote the generator of the genealogical process of a
sample of $n$ individuals when the number of demes is $D$. For each
$D$, the domain $\mathcal{D}(G^{n,D})$ of $G^{n,D}$ contains the
measurable symmetric functions of $n$ variables (by symmetric, we
mean invariant under all permutations of the variables). From the
last remark, we see that for all $f\in \mathcal{D}(G^{n,D})$, the
parts of $G^{n,D}f$ corresponding to 1-events vanish on $\Pi_n$.
Furthermore, we can define linear operators $\Psi^n$, $\Gamma^n$ and
$R_D^n$ such that $G^{n,D}$ has the following form:
$$G^{n,D}= r_D\Psi^n + \Gamma^n + R_D^n.$$
More precisely, for every function $f$ as above and each $\zeta\in
\Pn^s$, we have
$$
\Psi^{n}f(\zeta)=\sum_{\eta\in
\Pn^s}\vt_{\hat{k}(\zeta),\hat{k}(\eta)}\big(f(\eta)-f(\zeta)\big)\
\quad \mathrm{and}\quad\
\Gamma^{n}f(\zeta)=\sum_{\eta\in\Pn^s}\big\{\nu^{(n)}(\zeta,\eta)+
\mathit{l}_k^{(n)}(\zeta,\eta)\big\}\big(f(\eta)-f(\zeta)\big),
$$
and by the nonnegativity of their coefficients, these two operators
can each be viewed as generating a jump process independent of $D$.
In particular, we can define the structured genealogical process
$\xi$ on $\Pn^s$ as the process generated by $\Psi^{n}$. The
remaining terms $o(1)$ in Assumptions \ref{ass1} and \ref{ass2}
constitute the coefficients of the (not necessarily positive)
operator $R_D^n$, and so if we again use the operator norm
introduced in (\ref{norm}),
the finiteness of the number of possible transitions guarantees that
$\langle R_D^n\rangle=o(1)$ as $D\rightarrow \infty$.

\subsection{Convergence of the structured genealogical processes}
The main result of this section is the convergence of the
finite-dimensional distributions of the $\Pn^s$-valued structured
genealogical processes $\mcP^D$ to the corresponding ones of $\mcP$,
except at time $t=0$. The difficulty stems from the fact that 
the sequence of generators $G^{n,D}$ is unbounded because
of the fast genealogical events driven by $\Psi^n$. The proof
consists in essence in showing that the dynamics of the genealogical
processes become very close to the description of the dynamics of
$\mcP$, in that for $D$ large enough, once a $\Gamma^n$-event (i.e.,
a geographical collision) occurs, enough $\Psi^n$-events happen in a
very short period of time to bring the structured partition back
into $\Pi_n$. During that short period, the probability that a
$\Gamma^n$- or an $R_D^n$-event occurs is vanishingly small so that
at the time when $\mcP^D$ re-enters $\Pi_n$, with a high probability
it has the distribution of the final state of $\xi$ started at the
structured partition created by the geographical collision. Overall,
$R_D^n$-events are more and more infrequent and do not occur in the
limit.

Before stating the results of this section, let us define the
probability measures of interest. We take for granted the fact that
the processes $\xi$ and $\mcP^D$ for each $D\in \N$ and all $n\in\N$
can be constructed on the same probability space
$(\Omega,\mbP,\mathcal{F})$. For all $\zeta \in \Pn^s$, we thus
denote the probability measure under which these processes start at
$\zeta$ by $\mbP_{\zeta}$. Likewise, let
$(\Omega',\mbbP,\mathcal{F}')$ be the probability space on which the
processes $\mcP$ and $\chi$ (see Definition \ref{defi chi}) are
defined for all $n\in \N$. $\mbbP_{\eta}$ denotes the probability
measure under which these processes start at $\eta\in \Pi_n$.

With this notation, Theorem \ref{convergence p-dim} can be restated
as:

\bigskip
\noindent \textbf{Theorem \ref{convergence p-dim}'.} \emph{Suppose
that the conditions stated in Section \ref{section conditions} hold,
and let $\zeta \in \Pn^s$. Then, the structured genealogical
processes $\mcP^D$ started at $\zeta$ converge to the process $\mcP$
started at $\und{\zeta}$ as $D$ tends to infinity, in the sense that
for all $0<t_1<\ldots<t_p$,
$$\mbP_{\zeta}(\mcP^D_{t_1},\ldots, \mcP^D_{t_p}) \Rightarrow
\mbbP_{\und{\zeta}}(\mcP_{t_1},\ldots, \mcP_{t_p}) \qquad \qquad
\mathrm{as\ }D \rightarrow \infty,
$$ where $\mbP_{\zeta}(X)$ stands for the law of the random
variable $X$ under $\mbP_{\zeta}$ and $\mbbP_{\und{\zeta}}(X)$ is defined similarly.}
\bigskip

We also have the following result.
\begin{prop}\label{convergence fast time scale} Assume
again that the conditions of Section \ref{section conditions} hold.
Then the sequence of $D_{\Pn^s}[0,\infty)$-valued processes
$\{\mcP^D_{r_D^{-1}t},t\geq 0\}$ converges in distribution to the
structured genealogical process $\xi$ introduced in Section
\ref{paragraphe xi}.
\end{prop}

The proof of Proposition \ref{convergence fast time scale} is a
direct consequence of the uniform convergence of the generator of
$\mcP^D_{r_D^{-1}\cdot}$ (namely $r_D^{-1}G^{n,D}$) to the generator of
$\xi$ and the finiteness of the state space. A coupling with $\xi$
shows that the first time at which both processes differ when started
from the same value tends to infinity in probability, which is the
main argument to obtain the desired convergence. The proof being
immediate, we turn instead to the proof of Theorem \ref{convergence
p-dim}.

Let us first introduce the following notation, for each $D\in \N$:
$$\sigma_1^D\equiv \inf \{t\geq 0: \mcP^D_t \in \Pi_n\}, \qquad
\tau_1^D\equiv \inf \{t> \sigma^D_1: \mathrm{\ a\ }2\mathrm{-event\
occurs\ at\ }t\},$$ and for all $i\geq 2$,
\begin{equation}\label{defi sigma tau} \sigma_i^D\equiv \inf \{t\geq
\tau_{i-1}^D: \mcP^D_t \in \Pi_n\}, \qquad \tau_i^D\equiv \inf \{t>
\sigma^D_i: \mathrm{\ a\ }2\mathrm{-event\ occurs\ at\ }t\},
\end{equation}
with the convention that $\inf\emptyset = +\infty$ and if
$\sigma_i^D$ or $\tau_i^D=+\infty$, then the following random times
are all equal to $+\infty$. Note that if a 2-event occurs, its
outcome may still be in $\Pi_n$ (if all lineages gathered in
identical demes merge into one lineage in each of these demes). In
that case, $\sigma_{i+1}^D=\tau_i^D$. Let us also denote the ranked
epochs of events occurring to the process $\mcP$ by $\sigma_i,\
i\geq1$, including what we previously called the `ghost events',
with the conventions that $\sigma_1=0$ and $\sigma_{k}=+\infty$ for
$k\geq j+1$ if there are no more events after the $j$'th transition.

\medskip
\noindent \textit{Proof of Theorem \ref{convergence p-dim}.} We
start by proving the convergence of the one-dimensional
distributions, then establish the convergence of the finite
dimensional distributions by induction on their dimension. Since the
sample size is fixed, we drop the superscript $n$ in our notation.

As a first step, let us state the following definition and two
lemmas, which will be useful in the course of the proof. For the
sake of clarity, the proofs of the lemmas are postponed until after
the proof of Theorem \ref{convergence p-dim}.
\begin{defi}\label{defi chi} Let $(\chi_t, t\in [0,T))$ denote a
$\Pi_n$-valued Markov process generated by $\Gamma^n$, where $T$ is
defined as
$$T\equiv \inf \big\{t\geq 0:\ \chi_t \notin \Pi_n \big\}.$$
Then, for all $\eta\in \Pi_n$, $\chi(\eta)$ is defined as a
$\Pn^s$-random variable distributed like the outcome of the first
geographical collision when $\chi$ starts at $\eta$ (this event is
always defined if the $\lambda^g$'s satisfy (\ref{consist lambda g})
and $\eta$ has at least two blocks, since the coefficients
$\tilde{\lambda}^g$ are the rates of a $\Xi$-coalescent as mentioned
in Remark \ref{tilde lambda}).
\end{defi}
\begin{lemm}\label{conv tau sigma}
Let $i\geq 1$. Then for all bounded measurable functions $f$ on
$\mathbf{R}_+ \times \Pn^s$, we have
\begin{eqnarray*}
\lim_{D\rightarrow \infty}
\mbE_{\zeta}\big[f(\sigma_i^D,\mcP^D_{\sigma_i^D})\
\ind_{\{\sigma_i^D<\infty\}}\big]&=&
\mbbE_{\und{\zeta}}\big[f(\sigma_i,\mcP_{\sigma_i})\
\ind_{\{\sigma_i<\infty\}}\big],\\
\lim_{D\rightarrow \infty}
\mbE_{\zeta}\big[f(\tau_i^D,\mcP^D_{\tau_i^D})\
\ind_{\{\tau_i^D<\infty\}} \big]&=&
\mbbE_{\und{\zeta}}\big[f\big(\sigma_{i+1},\chi(\mcP_{\sigma_i})\big)\
\ind_{\{\sigma_{i+1}<\infty\}}\big].
\end{eqnarray*}
In particular, by taking $f(t,\eta)=\ind_{\{t\leq s\}}$ for all
$s>0$, we obtain that the law under $\mbP_{\zeta}$ of the
$[0,+\infty]$-valued random variable $\sigma_i^D$ (resp. $\tau_i^D$)
converges to the law under $\mbbP_{\und{\zeta}}$ of $\sigma_i$ (resp.
$\sigma_{i+1}$).
\end{lemm}

\begin{lemm}\label{value P_t}
Let $t\in (0,\infty)$ and let $i\in \N$ be such that
$\mbbP_{\und{\zeta}}[\sigma_i<\infty]>0$. By Lemma \ref{conv tau
sigma}, we also have for $D$ large enough
$\mbP_{\zeta}[\sigma^D_i<\infty]>0$. Let $f$ be a real-valued function
on $\Pn^s$. Then
$$\lim_{D\rightarrow \infty}\mbE_{\zeta}\big[f(\mcP^D_t)\ \ind_{[\sigma_i^D,\tau_i^D)}(t)\ \big|\
\sigma_i^D<\infty \big] = \mbbE_{\und{\zeta}}\big[f(\mcP_t)\
\ind_{[\sigma_i,\sigma_{i+1})}(t)\ \big|\ \sigma_i<\infty \big].$$
\end{lemm}

Fix $t>0$, let $f$ be a real-valued function on $\Pn^s$ and denote
the supremum norm of $f$ by $\|f\|$. We have for each $D$ and all $N
\in \N$: \setlength\arraycolsep{1pt}
\begin{eqnarray}
\Big|\mbE_{\zeta}\big[f(\mcP^D_t)\big]&-&\mbbE_{\und{\zeta}}\big[f(\mcP_t)\big]\Big|
\nonumber \\
&=& \Bigg|\mbE_{\zeta}\Big[\sum_{i=1}^{\infty}f(\mcP^D_t)\
\ind_{[\tau_{i-1}^D,\sigma_i^D)}(t) +
\sum_{i=1}^{\infty}f(\mcP^D_t)\ \ind_{[\sigma_i^D,\tau_i^D)}(t)\Big]
-\mbbE_{\und{\zeta}}
\Big[\sum_{i=1}^{\infty}f(\mcP_t)\ \ind_{[\sigma_i,\sigma_{i+1})}(t)\Big]\Bigg|\nonumber \\
&\leq&\ \sum_{i=1}^{N}\left|\mbE_{\zeta}\big[f(\mcP^D_t)\
\ind_{[\sigma_i^D,\tau_i^D)}(t)\big]-
\mbbE_{\und{\zeta}}\big[f(\mcP_t)\
\ind_{[\sigma_i,\sigma_{i+1})}(t)\big]\right|
\label{main ineq}\\
& &+\ \sum_{i=1}^{N}\mbE_{\zeta}\big[\big|f(\mcP^D_t)\big|\
\ind_{[\tau_{i-1}^D,\sigma_i^D)}(t)\big]+
\Big|\mbE_{\zeta}\big[f(\mcP^D_t)\ind_{\{t\geq
\tau^D_{N}\}}\big]\Big|+
\Big|\mbbE_{\und{\zeta}}\big[f(\mcP_t)\ind_{\{t\geq
\sigma_{N+1}\}}\big]\Big|, \nonumber
\end{eqnarray}
where $\tau_0 \equiv 0$. Let $\epsilon >0$. The random variables
$\sigma_i$ are the jump times of $\mcP$, the rates of which are
bounded above by a constant $b\geq 0$. Thus, for each $N\geq 1$,
$\sigma_N$ is bounded below by the sum of $N$ independent
exponentials with parameter $b$, and so there exists $N\geq 1$ such
that
$$\mbbP_{\und{\zeta}}[\sigma_{N+1}< t] <
\frac{\epsilon}{4\|f\|}\ .$$
In addition, $\tau^D_N \Rightarrow
\sigma_{N+1}$ by Lemma \ref{conv tau sigma}, so there exists a $D_0$
such that for $D\geq D_0$, $$\mbP_{\zeta}[\tau^D_{N+1}<\infty] <
\frac{\epsilon}{4\|f\|}.$$ Consequently, for $D\geq D_0$ we have
\begin{equation}\label{part 1}\Big|\mbE_{\zeta}\big[f(\mcP^D_t)\ind_{\{t\geq \tau^D_{N}\}}\big]\Big|
+\Big|\mbbE_{\und{\zeta}}\big[f(\mcP_t)\ind_{\{t\geq
\sigma_{N+1}\}}\big]\Big| \leq \|f\|
\left(\mbP_{\zeta}[\tau^D_{N}\leq
t]+\mbbP_{\und{\zeta}}[\sigma_{N+1}\leq t]\right) <
\frac{\epsilon}{2}.\end{equation}
Let $i\in \{1,\ldots,N\}$. We have
$$\mbE_{\zeta}\big[\big|f(\mcP^D_t)\big|\ \ind_{[\tau_{i-1}^D,\sigma_i^D)}(t)\big]\leq \|f\|\
\mbP_{\zeta}\big[\tau_{i-1}^D\leq t < \sigma_i^D\big]= \|f\|\
\left(\mbP_{\zeta}\big[\tau_{i-1}^D\leq t \big]
-\mbP_{\zeta}\big[\sigma_i^D \leq t \big]\right).$$ By Lemma
\ref{conv tau sigma}, both $\tau_{i-1}^D$ and $\sigma_i^D$ converge
in law towards $\sigma_i$ (whose distribution function is continuous
on $\mathbf{R}_+$), so the right-hand side of the last inequality
tends to 0 when $D\rightarrow \infty$. Hence, there exists a $D_1$
such that for all $D\geq D_1$,
\begin{equation}\label{part 2}\sum_{i=1}^{N}\mbE_{\zeta}\big[\big|f(\mcP^D_t)\big|\
\ind_{[\tau_{i-1}^D,\sigma_i^D)}(t)\big] \leq
\frac{\epsilon}{4}.\end{equation}
Once again, let $i\in
\{1,\ldots,N\}$. If $\mbbP_{\und{\zeta}}\big[\sigma_i \leq t\big]=0$,
then $\mbbE_{\und{\zeta}}\big[f(\mcP_t)\
\ind_{[\sigma_i,\sigma_{i+1})}(t)\big]=0$ and
$$\left|\mbE_{\zeta}\big[f(\mcP^D_t)\ \ind_{[\sigma_i^D,\tau_i^D)}(t)\big]\right|\leq \|f\|\
\mbP_{\zeta}\big[\sigma^D_i \leq t\big] \rightarrow 0 $$ as $D$
tends to infinity, by Lemma \ref{conv tau sigma} and the continuity
of the distribution function of $\sigma_i$ in $t$. If
$\mbbP_{\und{\zeta}}\big[\sigma_i\leq t\big]>0$, we also have
$\mbP_{\zeta}\big[\sigma^D_i \leq t\big]>0$ for $D$ large enough, so
we can write
\begin{eqnarray*} \mbE_{\zeta}\big[f(\mcP^D_t)\ \ind_{[\sigma_i^D,\tau_i^D)}(t)\big]&=&
\mbE_{\zeta}\big[f(\mcP^D_t)\ \ind_{[\sigma_i^D,\tau_i^D)}(t)\big|\
\sigma_i^D <\infty\big]
\mbP_{\zeta}[\sigma_i^D <\infty] \\
&\rightarrow & \mbbE_{\und{\zeta}}\big[f(\mcP_t)\
\ind_{[\sigma_i,\sigma_{i+1})}(t)\big|\
\sigma_i <\infty\big]\mbbP_{\und{\zeta}}[\sigma_i <\infty] \\
& &= \ \mbbE_{\und{\zeta}}\big[f(\mcP_t)\
\ind_{[\sigma_i,\sigma_{i+1})}(t)\big],
\end{eqnarray*}
where the convergence on the second line stems from Lemma \ref{value
P_t} and the convergence in distribution of $\sigma^D_i$ towards
$\sigma_i$. Consequently, there exists $D_2$ such that for all
$D\geq D_2$,
\begin{equation}\label{part 3}\sum_{i=1}^{N}\left|\mbE_{\zeta}\big[f(\mcP^D_t)\
\ind_{[\sigma_i^D,\tau_i^D)}(t)\big]-\mbbE_{\und{\zeta}}\big[f(\mcP_t)\
\ind_{[\sigma_i,\sigma_{i+1})}(t)\big]\right| \leq
\frac{\epsilon}{4}.\end{equation}
Combining to (\ref{main ineq}),
(\ref{part 1}), (\ref{part 2}) and (\ref{part 3}), we obtain for
all $D\geq \max \{D_0,D_1,D_2\}$
$$ \left|\mbE_{\zeta}\big[f(\mcP^D_t)\big]-\mbbE_{\und{\zeta}}\big[f(\mcP_t)\big]\right|
\leq \epsilon. $$ We can hence conclude that
$$\lim_{D\rightarrow \infty} \mbE_{\zeta}\big[f(\mcP^D_t)\big]=\mbbE_{\und{\zeta}}
\big[f(\mcP_t)\big],$$ which completes the proof of the convergence
of the one-dimensional distributions of $\mcP^D$ under
$\mbP_{\zeta}$ to the corresponding ones of $\mcP$ under
$\mbbP_{\und{\zeta}}$.

Let us now turn to the convergence of the finite-dimensional
distributions. We prove by induction on $p$ that, for all
$0<t_1<\ldots<t_p$, $\mbP_{\zeta}(\mcP^D_{t_1},\ldots, \mcP^D_{t_p})
\Rightarrow \mbbP_{\und{\zeta}}(\mcP_{t_1},\ldots, \mcP_{t_p})$ as $D
\rightarrow \infty$. By the preceding step, the case $p=1$ is
already established. Let $p\geq 2$, and suppose that the convergence
holds for the $(p-1)$-dimensional distributions. Let $0< t_1<\ldots<
t_p$, and let $f_1,\ldots,f_p$ be real-valued functions on
$\Pn^s$. We denote the $\sigma$-field generated
by $\big\{\mcP^D_s,\ s\in [0,t]\big\}$ by $\mathcal{F}^D_t$. Then,
\begin{eqnarray*}\mbE_{\zeta}\Big[\prod_{i=1}^p f_i(\mcP^D_{t_i})\Big] &=& \mbE_{\zeta}
\Big[\mbE\Big[\prod_{i=1}^p f_i(\mcP^D_{t_i})\Big|\mathcal{F}^D_{t_{p-1}}\Big]\Big] \\
 &=&  \mbE_{\zeta}\Big[\prod_{i=1}^{p-1} f_i(\mcP^D_{t_i})\mbE_{\mcP^D_{t_{p-1}}}
 \big[f_p(\tilde{\mcP}^D_{t_p-t_{p-1}})\big]\Big] \qquad \mathrm{by\ the\ Markov\ property} \\
&=& \mbE_{\zeta}\Big[\prod_{i=1}^{p-1} f_i(\mcP^D_{t_i})
\Big(\sum_{\eta \in \Pn^s} p^D(\mcP^D_{t_{p-1}},\eta, t_p-t_{p-1})
f_p(\eta)\Big)\Big],
\end{eqnarray*}
where here and in the following $\tilde{X}$ denotes an independent
version of the random variable $X$, the second expectation is taken
with regards to $\tilde{X}$, and $p^D(\cdot,\cdot,s)$ is the
transition kernel of $\mcP^D$ corresponding to time $s$. Continuing
the preceding equalities, we obtain
\begin{eqnarray}\mbE_{\zeta}\Big[\prod_{i=1}^p f_i(\mcP^D_{t_i})\Big] &=& \sum_{\eta \in \Pn^s}
f_p(\eta)\ \mbE_{\zeta}\Big[\prod_{i=1}^{p-1} f_i(\mcP^D_{t_i})\
p^D(\mcP^D_{t_{p-1}},\eta, t_p- t_{p-1})\Big] \nonumber\\ &=&
\sum_{\eta \in \Pn^s}f_p(\eta)\ \mbE_{\zeta}\Big[\prod_{i=1}^{p-1}
f_i(\mcP^D_{t_i})\ p^D(\mcP^D_{t_{p-1}},\eta, t_p-t_{p-1})\
\ind_{\{\mcP^D_{t_{p-1}}\notin\
\Pi_n\}}\Big] \label{ind non-pi}\\
& & +\ \sum_{\eta \in \Pn^s}f_p(\eta)\
\mbE_{\zeta}\Big[\prod_{i=1}^{p-1} f_i(\mcP^D_{t_i})\
p^D(\mcP^D_{t_{p-1}},\eta, t_p-t_{p-1})\
\ind_{\{\mcP^D_{t_{p-1}}\in\ \Pi_n\}}\Big]. \label{ind pi}
\end{eqnarray}
For all $\eta \in \Pn^s$,
\begin{equation} \label{conv non-pi}\left|\mbE_{\zeta}\Big[\prod_{i=1}^{p-1} f_i(\mcP^D_{t_i})\
p^D(\mcP^D_{t_{p-1}},\eta, t_p-t_{p-1})\
\ind_{\{\mcP^D_{t_{p-1}}\notin\ \Pi_n\}}\Big]\right| \leq
\Big(\prod_{i=1}^{p-1}\|f_i\|\Big)\
\mbP_{\zeta}\big[\mcP^D_{t_{p-1}}\notin\ \Pi_n\big] \rightarrow 0
\end{equation}
by the convergence of $\mcP^D_{t_{p-1}}$ to $\mcP_{t_{p-1}}$ in
distribution and the finiteness of $\Pn^s$. As the sum in (\ref{ind
non-pi}) is finite, (\ref{conv non-pi}) implies that this sum tends
to $0$ when $D$ grows to infinity. Moreover, the convergence in law
of $\mcP^D_{t_p-t_{p-1}}$ to $\mcP_{t_p-t_{p-1}}$, the finiteness of
$\Pn^s$ and the fact that $\und{\gamma}=\gamma$ a.s. if $\gamma \in
\Pi_n$ enable us to write
\begin{equation}\label{pd-p to 0}\max_{\gamma \in \Pi_n} \max_{\eta \in \Pn^s} \big|p^D(\gamma,\eta,t_p-t_{p-1})-
p(\gamma,\eta,t_p-t_{p-1})\big| \rightarrow 0 \qquad \mathrm{as\ } D\rightarrow \infty,
\end{equation}
where $p(\gamma,\eta,t_p-t_{p-1})$ is the transition kernel of
$\mcP$ corresponding to time $t_p-t_{p-1}$, extended to $\eta \notin
\Pi_n$ by $p(\gamma,\eta,t_p-t_{p-1})=0$. Now, we have for all $\eta
\in \Pn^s$ \setlength\arraycolsep{1pt}
\begin{eqnarray}
\mbE_{\zeta}&\Big[& \prod_{i=1}^{p-1} f_i(\mcP^D_{t_i})\
p^D(\mcP^D_{t_{p-1}},\eta, t_p-t_{p-1})\
\ind_{\{\mcP^D_{t_{p-1}}\in\ \Pi_n\}}\Big] \nonumber \\
&=& \mbE_{\zeta}\Big[\prod_{i=1}^{p-1} f_i(\mcP^D_{t_i})\
\big(p^D(\mcP^D_{t_{p-1}},\eta, t_p-t_{p-1})
-p(\mcP^D_{t_{p-1}},\eta, t_p-
t_{p-1})\big)\ \ind_{\{\mcP^D_{t_{p-1}}\in\ \Pi_n\}}\Big] \label{conv pi1}\\
& & +\ \mbE_{\zeta}\Big[\prod_{i=1}^{p-1} f_i(\mcP^D_{t_i})\
p(\mcP^D_{t_{p-1}},\eta, t_p-t_{p-1})\ \ind_{\{\mcP^D_{t_{p-1}}\in\
\Pi_n\}}\Big] \label{conv pi2}
\end{eqnarray}
The expression in (\ref{conv pi1}) tends to $0$ by (\ref{pd-p to 0})
and dominated convergence. As for the quantity in (\ref{conv pi2}),
for each $\eta \in \Pn^s$ the function $\gamma \mapsto
p(\gamma,\eta, t_p-t_{p-1})\ \ind_{\{\gamma \in\ \Pi_n\}}$
(vanishing on $\Pn^s \setminus \Pi_n$) is necessarily continuous and
bounded on the finite set $\Pn^s$, so by the induction hypothesis
for $p-1$, we have
$$\lim_{D\rightarrow \infty} \mbE_{\zeta}\Big[\prod_{i=1}^{p-1} f_i(\mcP^D_{t_i})\
p(\mcP^D_{t_{p-1}},\eta, t_p-t_{p-1})\ \ind_{\{\mcP^D_{t_{p-1}}\in\
\Pi_n\}}\Big] = \mbbE_{\und{\zeta}}\Big[\prod_{i=1}^{p-1}
f_i(\mcP_{t_i})\ p(\mcP_{t_{p-1}},\eta, t_p- t_{p-1})\Big].$$ The
two latter results, together with (\ref{ind pi}), (\ref{conv
non-pi}), the finiteness of the sums and the Markov property applied
to $\mcP$ lead to
\begin{eqnarray*} \lim_{D\rightarrow \infty} \mbE_{\zeta}\Big[\prod_{i=1}^p f_i(\mcP^D_{t_i})\Big] &=&
\sum_{\eta \in \Pn^s}f_p(\eta)\
\mbbE_{\und{\zeta}}\Big[\prod_{i=1}^{p-1} f_i(\mcP_{t_i})\
p(\mcP_{t_{p-1}},\eta, t_p-t_{p-1})\Big] \\
&=& \mbbE_{\und{\zeta}}\Big[\prod_{i=1}^{p} f_i(\mcP_{t_i})\Big].
\end{eqnarray*}
As any real-valued function on $\big(\Pn^s\big)^p$ can be
obtained as a uniform limit of product functions, the convergence of
the $p$-dimensional distributions is proven. The proof of Theorem
\ref{convergence p-dim} is complete by the induction principle.
\begin{flushright} $\Box$ \end{flushright}

\noindent \textit{Proof of Lemma \ref{conv tau sigma}.}
Let us start by proving that $\sigma_1^D$ converges in probability
to $0$. If $\zeta \in \Pi_n$, then $\sigma_1^D=0$ a.s. for all $D$
so the convergence trivially holds. If $\zeta \notin \Pi_n$, then
$\sigma^D_1>0$ a.s. and with the notation introduced previously, we
have for each function $f$ on $\Pn^s$
$$G^Df(\zeta)=r_D\Psi f(\zeta)+ \Gamma f(\zeta) + R_Df(\zeta),$$
where $\Psi f(\zeta)\neq 0$ (in fact, this holds for any $\eta
\notin \Pi_n$, and consequently for all values of $\mcP^D_t,\ t \in
[0,\sigma^D_1)$). Let us write $r_D c_{\Psi}(\zeta)$ (resp.
$c_{\Gamma}(\zeta)$, $c_{R_D}(\zeta)$) the total rate of the
non-trivial events generated by $r_D \Psi$ (resp. $\Gamma$, $R_D$)
when $G^Df$ is applied to $\zeta$. As events are discrete for each
$D$, we can write for $s>0$ \setlength\arraycolsep{1pt}
\begin{eqnarray}\mbP_{\zeta}&\big[&\sigma_1^D >s \big] \nonumber \\
& \leq & \mbP_{\zeta}\big[\mathrm{at\ most\ }n\ \Psi\mathrm{-events\
and\ then\ a\ }\Gamma \mathrm{-\ or\ }
R_D\mathrm{-event\ occur\ in\ }[0,s]\mathrm{\ and\ }\mcP_u^D \notin
\Pi_n\ \forall u\in [0,s]\big]\nonumber \\
& & +\mbP_{\zeta}\big[\mathrm{at\ most\ }n\ \Psi\mathrm{-events\
and\ no\ }\Gamma \mathrm{-\ or\ } R_D
\mathrm{-events\ occur\ in\ }[0,s]\mathrm{,\ and\ }\mcP_u^D \notin \Pi_n\
\forall u\in [0,s]\big]\nonumber \\
& & + \mbP_{\zeta}\big[\mathrm{more\ than\ }n\ \Psi\mathrm{-events\
occur\ before\ the\ first\ }\Gamma \mathrm{-\ or\ }
R_D\mathrm{-event}\big]. \label{conv sigma1}
\end{eqnarray}
Since the events generated by $\Psi$ correspond to the structured
genealogical process $(\xi_t,t\geq 0)$ started at $\zeta$  as long
as no $\Gamma$- or $R_D$-events occurred, by the bound on the number
of transitions of $\xi$ ($n$, see the previous section), the third
term on the right-hand side of (\ref{conv sigma1}) vanishes.
Moreover, the probability that the next event generated by $G^D$ is
a $\Gamma$ or an $R_D$-event when the current value of $\mcP^D$ is
$\eta \notin \Pi_n$ is given by
$$ \frac{c_{\Gamma}(\eta) + c_{R_D}(\eta)}{c_{\Gamma}(\eta) + c_{R_D}(\eta)
+r_Dc_{\Psi}(\eta)} \rightarrow 0,
\qquad D\rightarrow \infty,$$ since $c_{\Psi}(\eta)>0$ for such an
$\eta$, and this is precisely the kind of situation required to be
in the configuration given by the first term of (\ref{conv sigma1}).
So by bounding this term by the maximum over $\eta\notin \Pi_n$ of
the probabilities calculated just before, we obtain that the first
term of (\ref{conv sigma1}) tends to $0$ as $D$ grows to infinity.
To finish, for each $D$ and all $k\in \{1,\ldots,n\}$ let us call
$U^D_k$ the random time of the $k$'th event occurring to $\mcP^D$,
with the convention that $U^D_k=+\infty$ if there are less than $k$
such events. If $k$ events occur (i.e., $U^D_k<\infty$) and $\mcP^D$
stays out of $\Pi_n$, then $U_{k+1}$ is stochastically bounded by
the sum of $k+1$ i.i.d. exponential variables with parameter $r_D
\min_{\eta\notin \Pi_n} c_{\Psi}(\eta)$, whose distribution becomes
concentrated close to $0$ as $D$ grows since $\min_{\eta\notin
\Pi_n} c_{\Psi}(\eta)>0$. Consequently, \setlength\arraycolsep{1pt}
\begin{eqnarray*}
\mbP_{\zeta}&\big[&\mathrm{exactly\ }k\ \Psi\mathrm{-events\ and\
no\ }\Gamma \mathrm{-\ or\ } R_D\mathrm{-events\
occur\ in\ }[0,s]\mathrm{,\ and\ }\mcP_u^D \notin \Pi_n\ \forall u\in [0,s]\big] \\
&\leq & \mbP_{\zeta}\big[U_{k+1}^D>s,\ U_k^D<\infty \mathrm{\ and\ }
\mcP_u^D \notin \Pi_n\ \forall u\in [0,s] \big] \rightarrow 0.
\end{eqnarray*}
As the second term in (\ref{conv sigma1}) is bounded by the sum over
$k\in \{0,\ldots,n\}$ of the preceding quantities, it converges to
zero. Hence, $\mbP_{\zeta}[\sigma_1^D>s]\rightarrow 0$ for all $s>0$
and $\sigma_1\rightarrow 0$ in probability.

Now, let $f$ be a function on $\Pn^s$. For each
$s>0$, we have
\begin{equation}\label{conv P1 intro} \mbE_{\zeta}\big[f(\mcP^D_{\sigma_1^D})\
\ind_{\{\sigma_1^D<s\}} \big]
=\mbE_{\zeta}\big[f(\mcP^D_{\sigma_1^D})\big] -
\mbE_{\zeta}\big[f(\mcP^D_{\sigma_1^D})\ \ind_{\{\sigma_1^D\geq
s\}}\big].
\end{equation}
By the convergence in probability of $\sigma^D_1$ to $0$ and the
fact that $f$ is bounded, the second term in the right-hand side of
(\ref{conv P1 intro}) vanishes as $D$ grows to infinity.
Furthermore, we have
\begin{eqnarray}\mbE_{\zeta}\big[f(\mcP^D_{\sigma_1^D})\big] &=& \mbE_{\zeta}
\big[f(\mcP^D_{\sigma_1^D});\ \mathrm{only\ }\Psi\mathrm{-events\ before\ }\sigma_1^D\big] \nonumber\\
& & +\ \mbE_{\zeta}\big[f(\mcP^D_{\sigma_1^D});\ \mathrm{at\ least\
one\ }\Gamma\mathrm{- or\ } R_D\mathrm{-events\ before\
}\sigma_1^D\big]. \label{conv P1}
\end{eqnarray}
The second term in (\ref{conv P1}) is bounded by
$\|f\|\mbP_{\zeta}\big[\mathrm{at\ least\ one\ }\Gamma \mathrm{- or\
}R_D\mathrm{-events\ before\ }\sigma_1^D\big]$ which tends to $0$ by
the preceding calculations, giving as a by-product that
$\mbP_{\zeta}\big[\mathrm{only\ }\Psi\mathrm{-events\ before\
}\sigma_1^D\big]>0$ for $D$ large enough. Moreover, when only
$\Psi-$events occurred between $0$ and $\sigma_1^D$, then the
evolution of $\mcP^D$ between these two times is driven by the
structured genealogical process $\xi$ started at $\zeta$, so
$\mcP^D_{\sigma_1^D}$ has the same distribution as $\und{\zeta}$.
Thus, \setlength\arraycolsep{1pt}
\begin{eqnarray*}
\mbE_{\zeta}\big[f(\mcP^D_{\sigma_1^D})&;&\ \mathrm{only\ }\Psi\mathrm{-events\ before\ }\sigma_1^D\big]\\
&=& \mbE_{\zeta}\big[f(\mcP^D_{\sigma_1^D}) \big|\ \mathrm{only\
}\Psi\mathrm{-events\ before\ }
\sigma_1^D\big]\mbP_{\zeta}\big[\mathrm{only\ }\Psi\mathrm{-events\ before\ }\sigma_1^D\big] \\
&=& \mbE\big[f(\und{\zeta}) \big]\mbP_{\zeta}\big[\mathrm{only\
}\Psi\mathrm{-events\ before\ }
\sigma_1^D\big]\\
& \rightarrow & \mbE\big[f(\und{\zeta}) \big].
\end{eqnarray*}
Together with (\ref{conv P1 intro}) and (\ref{conv P1}), we obtain that
$$\lim_{D\rightarrow \infty}\mbE_{\zeta}\big[f(\mcP^D_{\sigma_1^D})\ \ind_{\{\sigma_1^D<s\}}\big]
=\mbE[f(\und{\zeta})]= \mbbE_{\und{\zeta}}[f(\mcP_0)\ \ind_{\{0<s\}}].$$
A monotone class
argument of enables us to conclude the same result for any bounded
measurable function $f$ on $\mathbf{R}_+\times \Pn^s$.

Let us now investigate the convergence of $\tau_1^D$. Recall that if
$\eta \in \Pi_n$ and $f\in \mathcal{D}(G^D)$, then
$$G^Df(\eta)= \Gamma f(\eta)+ R_Df(\eta).$$
If $s>0$, by the strong Markov property applied to $\mcP^D$ at time
$\sigma^D_1$ we have
\begin{equation}\label{tau1}\mbP_{\zeta}[\tau_1^D>s] = \mbE_{\zeta}
\big[\mbP_{\mcP^D_{\sigma_1^D}}[\tilde{\tau}_1^D>s-\sigma_1^D]\
\ind_{\{s>\sigma_1^D\}}\big]+
\mbE_{\zeta}\big[\mbP_{\mcP^D_{\sigma_1^D}}[\tilde{\tau}_1^D>s-\sigma_1^D]\
\ind_{\{s\leq \sigma_1^D<\infty\}}\big].
\end{equation}
The second term in (\ref{tau1}) is equal to $\mbP_{\zeta}[s\leq
\sigma_1^D <\infty]$ which tends to $0$ when $D$ grows to infinity.
If a $\Gamma$-event occurs when the current value of $\mcP^D$ lies
in $\Pi_n$, it is necessarily a $2$-event, hence the first term is
equal to \setlength\arraycolsep{1pt}
\begin{eqnarray*}
\mbE_{\zeta}&\big[&\ind_{\{s>\sigma_1^D\}}\
\mbP_{\mcP^D_{\sigma_1^D}}[\mathrm{no\ }\Gamma
\mathrm{-\ or\ }R_D\mathrm{-events\ between\ }0\ \mathrm{and\ }s-\sigma_1^D] \big]\\
&+ & \mbE_{\zeta}\big[\ind_{\{s>\sigma_1^D\}}\
\mbP_{\mcP^D_{\sigma_1^D}}[\mathrm{no\ }\Gamma \mathrm{-events\ and\
at\ least\ one\ }R_D\mathrm{-event\ between\ }0\ \mathrm{and\
}s-\sigma_1^D;\ \tilde{\tau}_1^D>s- \sigma_1^D] \big].
\end{eqnarray*}
But for all $\eta \in \Pi_n$,
\setlength\arraycolsep{1pt}
\begin{eqnarray*}\mbP_{\eta}&[&\mathrm{no\ }\Gamma\mathrm{-events\ and\
at\ least\ one\ }R_D\mathrm{-event\
between\ }0\ \mathrm{and\ }s-\sigma_1^D;\ \tilde{\tau}_1^D>s-\sigma_1^D]\\
&\leq & \mbP_{\eta}[\mathrm{no\ }\Gamma\mathrm{-events\ and\ at\
least\ one\ }R_D\mathrm{-event\ between\ }0\
\mathrm{and\ }s] \\
&\leq & 1-\exp\big(-s\ \max_{\gamma \in \Pn^s}c_{R_D}(\gamma)\big)
\rightarrow 0,
\end{eqnarray*}
so by dominated convergence,
$$\mbE_{\zeta}\big[\ind_{\{s>\sigma_1^D\}}\ \mbP_{\mcP^D_{\sigma_1^D}}[\mathrm{no\ }\Gamma
\mathrm{-event\ and\ at\ least\ one\ }R_D\mathrm{-event\ between\
}0\ \mathrm{and\ }s-\sigma_1^D; \ \tilde{\tau}_1^D>s- \sigma_1^D]
\big] \rightarrow 0.$$ Consequently,
\begin{eqnarray*} \mbP_{\zeta}[\tau_1^D >s] &=& \mbE_{\zeta}\Big[\ind_{\{s>\sigma_1^D\}}\ \exp-
\Big\{\big(c_{\Gamma}(\mcP^D_{\sigma_1^D})+c_{R_D}(\mcP^D_{\sigma_1^D})\big)
(s-\sigma_1^D)\Big\}\Big] \\
&\rightarrow& \mbbE[e^{-s c_{\Gamma}(\und{\zeta})}]=
\mbbP_{\und{\zeta}}[\sigma_2>s]
\end{eqnarray*}
by the preceding convergence result for $(\sigma^D_1,
\mcP^D_{\sigma^D_1})$ and the uniform convergence of $c_{R_D}$
towards $0$. We can thus conclude that the law of $\tau_1^D$ under
$\mbP_{\zeta}$ converges to the law of $\sigma_2$ under
$\mbbP_{\und{\zeta}}$.

Now, by the strong Markov property
applied to $\mcP^D$ at time $\sigma^D_1$, we have
\setlength\arraycolsep{1pt}
\begin{eqnarray}\mbE_{\zeta}&[&\ind_{\{\tau_1^D<s\}}\ f(\mcP^D_{\tau_1^D})\
\ind_{\{\tau_1^D<\infty\}}] \label{conv P_tau1}\\
&=& \mbE_{\zeta}\big[\
\ind_{\{\sigma_1^D<\infty\}}\mbE_{\mcP^D_{\sigma^D_1}}[\
\ind_{\{\tilde{\tau}_1^D<s-\sigma^D_1\}}\
f(\tilde{\mcP}^D_{\tilde{\tau}_1^D});\
\mathrm{the\ first\ event\ is\ an\ }R_D\mathrm{-event}]\big] \nonumber \\
& & +\ \mbE_{\zeta}\big[\
\ind_{\{\sigma_1^D<\infty\}}\mbE_{\mcP^D_{\sigma^D_1}}[\
\ind_{\{\tilde{\tau}_1^D<s-\sigma^D_1\}}\
f(\tilde{\mcP}^D_{\tilde{\tau}_1^D});\ \mathrm{the\ first\ event\
is\ a\ }\Gamma \mathrm{-event}]\big]. \nonumber
\end{eqnarray}
The absolute value of the first term in the right-hand side of
(\ref{conv P_tau1}) is bounded by
$$\|f\|\ \max_{\eta \in \Pi_n}\mbP_{\eta}[\mathrm{a\ first\ event\ occurs\ and\ is\ an\ }
R_D\mathrm{-event}] \rightarrow 0.$$
Moreover, if $\tilde{\mcP}^D_{0}=\eta \in \Pi_n$ and the first event
is a $\Gamma$-event, then $\tilde{\tau}_1^D$ is the time of that
first event and $\tilde{\mcP}^D_{\tilde{\tau}_1^D}$ its outcome.
Therefore, both are independent and
$\tilde{\mcP}^D_{\tilde{\tau}_1^D}$ is distributed like
$\chi(\eta)$, so the second term in (\ref{conv P_tau1}) is equal to
\setlength\arraycolsep{1pt}
\begin{eqnarray*}
\mbE_{\zeta} &\big[ &
\ind_{\{\sigma_1^D<\infty\}}\mbE_{\mcP^D_{\sigma^D_1}}[\
\ind_{\{\tilde{\tau}_1^D<s-\sigma^D_1\}} f(\chi(\tilde{\mcP}^D_0))\ ]\big] +o(1) \\
& = & \mbE_{\zeta}\big[\
\ind_{\{\sigma_1^D<\infty\}}\mbP_{\mcP^D_{\sigma^D_1}}
[\tilde{\tau}_1^D<s-\sigma^D_1]\
\mbE_{\mcP^D_{\sigma^D_1}}[f(\chi(\tilde{\mcP}^D_0))]\ \big] +o(1).
\end{eqnarray*}
Let us write
$$\mbP_{\mcP^D_{\sigma^D_1}}[\tilde{\tau}_1^D<s-\sigma^D_1] =
\mbP_{\mcP^D_{\sigma^D_1}}[\tilde{\tau}_1^D<s]-\mbP_{\mcP^D_{\sigma^D_1}}
\big[\tilde{\tau}_1^D \in [s-\sigma^D_1,s]\big]$$ and fix $\epsilon
>0$. For any $\delta>0$, we have \setlength\arraycolsep{1pt}
\begin{eqnarray*}
\mbE_{\zeta}\big[\ind_{\{\sigma_1^D<\infty\}}&\mbP_{\mcP^D_{\sigma^D_1}}&
\big[\tilde{\tau}_1^D\in [s-\sigma^D_1,s]\big]\
\mbE_{\mcP^D_{\sigma^D_1}}
[f(\chi(\tilde{\mcP}^D_0))]\ \big] \\
&\leq & \mbE_{\zeta}\big[\
\ind_{\{\sigma_1^D<\delta\}}\mbP_{\mcP^D_{\sigma^D_1}}\big[\tilde{\tau}_1^D
\in [s-\delta,s]\big]\
\mbE_{\mcP^D_{\sigma^D_1}}[|f(\chi(\tilde{\mcP}^D_0))|]\ \big] +
\|f\|\ \mbP_{\zeta}[\sigma_1^D \in [\delta,\infty)].
\end{eqnarray*}
By the convergence in probability of $\sigma_1^D$ to $0$, there
exists $D_0 \geq 1$ such that for all $D\geq D_0$,
$$\mbP_{\zeta}[\sigma_1^D \in [\delta,\infty)] < \frac{\epsilon}{3 \|f\|}.$$
Let $\eta \in \Pi_n$. By the continuity of the distribution function
of $\sigma_2$, there exists $\delta_0>0$ such that
$$\mbP_{\eta}\big[\tilde{\sigma}_2 \in [s-\delta_0,s]\big] < \frac{\epsilon}{3\|f\|}.$$
In addition, $\tilde{\tau}_1^D$ converges in distribution to $\sigma_2$,
hence there exists $D_1 \geq 1$ such that for all $D\geq D_1$,
$$\Big|\mbP_{\eta}\big[\tilde{\tau}_1^D \in [s-\delta_0,s]\big]-\mbP_{\eta}
\big[\tilde{\sigma}_2 \in [s-\delta_0,s]\big]\Big| <
\frac{\epsilon}{3\|f\|}.$$ Since $\Pi_n$ is a finite set, we can
conclude that for $\delta >0$ small enough, and $D$ large enough, we
have \setlength\arraycolsep{1pt}
\begin{eqnarray*}
\mbE_{\zeta}&\big[&
\ind_{\{\sigma_1^D<\delta\}}\mbP_{\mcP^D_{\sigma^D_1}}
\big[\tilde{\tau}_1^D \in [s-\delta,s]\big]\
\mbE_{\mcP^D_{\sigma^D_1}}
[|f(\chi(\tilde{\mcP}^D_0))|]\ \big] +\|f\|\ \mbP_{\zeta}[\sigma_1^D \in [\delta,\infty)]\\
 &\leq & \mbE_{\zeta}\big[\ \ind_{\{\sigma_1^D<\delta\}}\max_{\eta \in \Pi_n}
 \mbP_{\eta}\big[\tilde{\tau}_1^D \in [s-\delta,s]\big]\ \mbE_{\mcP^D_{\sigma^D_1}}
 [|f(\chi(\tilde{\mcP}^D_0))|]\ \big] + \frac{\epsilon}{3} \\
& \leq & \epsilon.
\end{eqnarray*}
Now, $\eta \mapsto \mbP_{\eta}[\tilde{\tau}_1^D <s]$ converges
uniformly in $\eta \in \Pi_n$ to $\eta \mapsto
\mbbP_{\eta}[\tilde{\sigma}_2 <s]$ and
$$(s,\eta)\mapsto \ind_{\{s<\infty\}}\ind_{\{\eta \in \Pi_n\}}\mbbP_{\eta}[\tilde{\sigma}_2<s]
\mbbE_{\eta}\big[f(\chi(\tilde{\mcP}))\big]$$ is a bounded measurable
function, so by the convergence in distribution of
$(\sigma_1^D,\mcP^D_{\sigma_1^D})$ proven above, for $D$ large
enough we have \setlength\arraycolsep{1pt}
\begin{eqnarray*}
\Big| \mbE_{\zeta}&[&\ind_{\{\tau_1^D<s\}}\ f(\mcP^D_{\tau_1^D})\
\ind_{\{\tau_1^D<\infty\}}]-\mbbE_{\und{\zeta}}[\ind_{\{\sigma_2<s\}}\
f(\chi(\mcP_{\sigma_1}))\ \ind_{\{\sigma_2<\infty\}}]\Big| \\
&=& \Big| \mbE_{\zeta}[\ind_{\{\tau_1^D<s\}}\ f(\mcP^D_{\tau_1^D})\
\ind_{\{\tau_1^D<\infty\}}]-\mbbE_{\und{\zeta}}\big[\ind_{\{\sigma_1<\infty\}}
\mbbP_{\mcP_{\sigma_1}}[\tilde{\sigma}_2<s]\ \mbbE_{\mcP_{\sigma_1}}
\big[f(\chi(\tilde{\mcP}_0))]\big]\Big| \\
&<& 3 \epsilon.
\end{eqnarray*}
Letting $\epsilon$ tend to zero yields the desired result (once
again by invoking monotone classes) and completes the step $i=1$ of
the proof of Lemma \ref{conv tau sigma}.

Suppose that the desired properties hold for $i-1$. Let $f$ be a
bounded continuous function on $\mathbf{R}_+\times \Pn^s$. Since
$\ind_{\{\sigma_i^D<\infty\}}=\ind_{\{\sigma_i^D<\infty\}}
\ind_{\{\tau_{i-1}^D<\infty\}}$, the strong Markov property applied
to $\mcP^D$ at time $\tau_{i-1}^D$ gives
$$
\mbE_{\zeta}\big[f(\sigma_i^D,\mcP^D_{\sigma_i^D})\ind_{\{\sigma^D_i<\infty\}}\big]=
\mbE_{\zeta}\big[\ind_{\{\tau_{i-1}^D<\infty\}}\mbE_{\mcP^D_{\tau^D_{i-1}}}[f(\tau^D_{i-1}+
\tilde{\sigma}_1^D,\tilde{\mcP}^D_{\tilde{\sigma}_1^D})\
\ind_{\{\tilde{\sigma}^D_1<\infty\}}]\big]
$$
But, for all $\eta \notin \Pi_n$, if $X$ denotes a random
variable whose distribution under $\mbP_{\eta}$ is that of
$\und{\eta}$ (e.g. $\xi_{\tau_{\pi}}$ in the notation of Proposition
\ref{description xi}), then
\begin{eqnarray*}\big|\mbE_{\eta}[f(t+\tilde{\sigma}_1^D,\tilde{\mcP}^D_{\tilde{\sigma}_1^D})
\ind_{\{\tilde{\sigma}^D_1<\infty\}}]-\mbE_{\eta}[f(t,X)]\big|&\leq&
\big|\mbE_{\eta}
[f(t+\tilde{\sigma}_1^D,\tilde{\mcP}^D_{\tilde{\sigma}_1^D})\
\ind_{\{\tilde{\sigma}^D_1<\infty\}}]-
\mbE_{\eta}[f(t+\tilde{\sigma}_1^D,X)\ \ind_{\{\tilde{\sigma}^D_1<\infty\}}]\big| \\
& & +\big|\mbE_{\eta}[f(t+\tilde{\sigma}_1^D,X)\
\ind_{\{\tilde{\sigma}^D_1<\infty\}}]-\mbE_{\eta}
[f(t,X)]\big|.\end{eqnarray*} Since
$\tilde{\mcP}^D_{\tilde{\sigma}_1^D}$ has the same distribution as
$\und{\eta}$ under $\mbP_{\eta}$ if only $\Psi$-events occurred
between $0$ and $\tilde{\sigma}_1^D$, the first term is equal to
\setlength\arraycolsep{1pt}
\begin{eqnarray*}\big|\mbE_{\eta}\big[\ind_{\{\tilde{\sigma}^D_1<\infty\}}
\big(f(t+\tilde{\sigma}_1^D,
\tilde{\mcP}^D_{\tilde{\sigma}_1^D})&-&f(t+\tilde{\sigma}_1^D,X)\big);\
\mathrm{not\ only\ } \Psi\mathrm{-events\ between\ }0\mathrm{\ and\
}\tilde{\sigma}_1^D
\big]\big|\\
&\leq & 2\|f\| \max_{\eta\notin
\Pi_n}\mbP_{\eta}[\tilde{\sigma}^D_1<\infty, \mathrm{not\ only\
}\Psi \mathrm{-events\ between\ }0\mathrm{\ and\
}\tilde{\sigma}_1^D]\rightarrow 0
\end{eqnarray*}
by the calculations done in the proof of the convergence of
$\sigma_1^D$. Moreover, $$\mbE_{\eta}[f(t+ \tilde{\sigma}_1^D,X)\
\ind_{\{\tilde{\sigma}^D_1<\infty\}}]-\mbE_{\eta}[f(t,X)]\rightarrow
0$$ uniformly in $\eta$ by the convergence in probability of
$\tilde{\sigma}_1^D$ towards $0$ and the finiteness of the number of
states that $X$ can take. Therefore,
$$\mbE_{\eta}[f(t+\tilde{\sigma}_1^D,\tilde{\mcP}_{\tilde{\sigma}_1^D}^D)\
\ind_{\{\tilde{\sigma}^D_1 <\infty\}}] \rightarrow
\mbE_{\eta}[f(t,X)]= \mbbE_{\und{\eta}}[f(t,\mcP_0)]$$ uniformly in
$(t,\eta)$. This uniform convergence (which trivially holds also for
$\eta\in \Pi_n$ since $\und{\eta}=\eta$ and $\sigma_1^D=0$ a.s.),
together with the induction hypothesis for $i-1$ yields
$$\lim_{D\rightarrow \infty}\mbE_{\zeta}\big[\ind_{\{\tau_{i-1}^D<\infty\}}
\mbE_{\mcP^D_{\tau^D_{i-1}}}[f(\tau^D_{i-1}+\tilde{\sigma}_1^D,\tilde{\mcP}^D_{\tilde{\sigma}_1^D})\
\ind_{\{\tilde{\sigma}^D_1<\infty\}}]\big]=\mbbE_{\und{\zeta}}\big[\ind_{\{\sigma_i<\infty\}}
\mbbE_{\chi(\mcP_{\sigma_{i-1}})}[f(\sigma_i,X)]\big].$$ But, from
the description of the evolution of $\mcP$ in terms of the
geographical events followed by the instantaneous action of the
structured genealogical process $\xi$, we see that the random
variable $\und{\chi(\mcP_{\sigma_{i-1}})}$ is distributed precisely
like $\mcP_{\sigma_i}$ (if $\sigma_{i}< \infty$). Consequently,
$$ \lim_{D\rightarrow \infty}\mbE_{\zeta}[f(\sigma_i^D,\mcP^D_{\sigma_i^D})
\ind_{\{\sigma^D_i<\infty\}}]=
\mbbE_{\und{\zeta}}\big[\ind_{\{\sigma_i<\infty\}}
f(\sigma_i,\mcP_{\sigma_i})\big].$$

The same technique applies to $(\tau_i^D,\mcP^D_{\tau_i^D})$, where
this time we use the strong Markov property at time $\sigma_i^D$ and
the following convergence result:
$$\mbE_{\eta}[f(t+\tau_1^D,\mcP^D_{\tau_1^D})\ \ind_{\{\tau_1^D<\infty\}}]\rightarrow
\mbbE_{\eta}[f(t+\sigma_2,\chi(Y))\ \ind_{\{\sigma_2<\infty\}}]$$
uniformly in $(t,\eta)\in \mathbf{R}_+\times \Pi_n$, where under
$\mbP_{\eta}$, $Y$ is a.s. equal to $\eta$.
\begin{flushright} $\Box$ \end{flushright}

\noindent \textit{Proof of Lemma \ref{value P_t}.}
If an event occurs in the (random) interval $[\sigma_i^D,\tau_i^D)$,
the first such event can be neither a $\Psi$-event since
$\mcP_{\sigma_i^D}^D \in \Pi_n$, nor a $\Gamma$-event since $\mcP^D$
would then undergo a $2$-event before time
$\tau_i^D$, contradicting the definition of $\tau_i^D$, so it
must be an $R_D$-event. Consequently, if we write
\begin{eqnarray}
\mbE_{\zeta}\big[f(\mcP^D_t)\ \ind_{[\sigma_i^D,\tau_i^D)}(t)\big|\
\sigma_i^D<\infty \big] &=& \mbE_{\zeta}\big[f(\mcP^D_t)\
\ind_{[\sigma_i^D,\tau_i^D)}(t);\ \mathrm{nothing\ happens\ in\ }
[\sigma_i^D,t]\big|\ \sigma_i^D<\infty \big] \label{eq1 lemme P_t}\\
& &+\ \mbE_{\zeta}\big[f(\mcP^D_t)\
\ind_{[\sigma_i^D,\tau_i^D)}(t);\ \mathrm{something\ occurs\ in\ }
[\sigma_i^D,t]\big|\ \sigma_i^D<\infty \big], \nonumber
\end{eqnarray}
then the absolute value of the second term of the right-hand side of
(\ref{eq1 lemme P_t}) is bounded by
$$\|f\|\ \mbP_{\zeta}\big[\sigma_i^D\leq t \mathrm{\ and\ an\ }
R_D\mathrm{-event\ occurs\ in\ }[\sigma_i^D,t]\big|
\ \sigma_i^D<\infty \big]\leq \|f\|\big(1-\exp(-t\
\max_{\eta}c_{R_D}(\eta))\big) \rightarrow 0, $$ where in the
exponential the maximum is over $\eta \in \Pn^s$ and recall that
$c_{R_D}(\eta)$ is the total rate at which $R_D$-events occur when
the current value of $\mcP^D$ is $\eta$.

The first term of (\ref{eq1 lemme P_t}) is equal to
\setlength\arraycolsep{1pt}
\begin{eqnarray*}
\mbE_{\zeta}&\big[&f(\mcP^D_{\sigma_i^D})\
\ind_{[\sigma_i^D,\tau_i^D)}(t);\ \mathrm{nothing\ happens\
in\ }[\sigma_i^D,t]\big|\ \sigma_i^D<\infty \big] \\
&=& \mbE_{\zeta}\big[f(\mcP^D_{\sigma_i^D})\
\ind_{[\sigma_i^D,\tau_i^D)}(t)\ \big|\ \sigma_i^D<\infty
\big]-\mbE_{\zeta}\big[f(\mcP^D_{\sigma_i^D})\
\ind_{[\sigma_i^D,\tau_i^D)}(t);\ \mathrm{something\ happens\ in\
}[\sigma_i^D,t]\big|\ \sigma_i^D<\infty \big]
\end{eqnarray*}
As before,
$$\mbE_{\zeta}\big[f(\mcP^D_{\sigma_i^D})\ \ind_{[\sigma_i^D,\tau_i^D)}(t);
\ \mathrm{something\ happens\ in\ }[\sigma_i^D,t]\big|\
\sigma_i^D<\infty \big]\leq \|f\|\big(1-\exp(-t\
\max_{\eta}c_{R_D}(\eta))\big) \rightarrow 0$$ and furthermore
\begin{equation}\label{eq2 lemme P_t}
\mbE_{\zeta}\big[f(\mcP^D_{\sigma_i^D})\
\ind_{[\sigma_i^D,\tau_i^D)}(t)\ \big|\ \sigma_i^D<\infty
\big]=\mbE_{\zeta}\big[f(\mcP^D_{\sigma_i^D})\ \ind_{\{\sigma_i^D
\leq t\}}\ \big|\ \sigma_i^D<\infty
\big]-\mbE_{\zeta}\big[f(\mcP^D_{\sigma_i^D})\ \ind_{\{\tau_i^D\leq
t\}}\ \big|\ \sigma_i^D<\infty \big].
\end{equation}
On the one hand, by Lemma \ref{conv tau sigma} and the fact that
$\mbP_{\zeta}[\sigma_i^D<\infty]\rightarrow
\mbbP_{\und{\zeta}}[\sigma_i<\infty]>0$, the first term in (\ref{eq2
lemme P_t}) converges as $D$ tends to infinity to
$$\mbbE_{\und{\zeta}}\big[f(\mcP_{\sigma_i})\ \ind_{\{\sigma_i \leq t\}}\ \big|\
\sigma_i<\infty \big].$$

On the other hand, by the strong Markov property applied to $\mcP^D$
at time $\sigma_i^D$, the second term in (\ref{eq2 lemme P_t}) is
equal to
$$\mbE_{\zeta}\big[f(\mcP^D_{\sigma_i^D})\ \mbP_{\mcP^D_{\sigma_i^D}}\big[\tilde{\tau}_1^D
\leq t\big]\ \big|\ \sigma_i^D<\infty \big].$$ The function
$\eta\mapsto \mbP_{\eta}\big[\tilde{\tau}_1^D\leq t\big]$ converges
uniformly in $\eta \in \Pi_n$ to $\eta\mapsto
\mbbP_{\eta}\big[\tilde{\sigma}_2\leq t\big]$, so by Lemma \ref{conv
tau sigma} we obtain
\begin{eqnarray*} \lim_{D\rightarrow \infty}\mbE_{\zeta}\big[f(\mcP^D_{\sigma_i^D})\
\mbP_{\mcP^D_{\sigma_i^D}}\big[\tilde{\tau}_1^D\leq t\big]\ \big|\
\sigma_i^D<\infty \big] &=&
\mbbE_{\und{\zeta}}\big[f(\mcP_{\sigma_i})\
\mbbP_{\mcP_{\sigma_i}}\big[\tilde{\sigma}_2
\leq t\big]\ \big|\ \sigma_i<\infty \big]\\
& = &\mbbE_{\und{\zeta}}\big[f(\mcP_{\sigma_i})\
\ind_{\{\sigma_{i+1}\leq t\}}\ \big|\ \sigma_i<\infty
\big]\end{eqnarray*} by the strong Markov property applied this time
to $\mcP$ at time $\sigma_i$. Combining the above, we obtain
the desired result.
\begin{flushright} $\Box$ \end{flushright}

The results obtained in this section are similar in spirit to
perturbation theorems such as Theorem 1.7.6 in \citet{EK1986}.
Indeed, in our case the existence of a projector $p$ corresponding
to $\Psi$ and the convergence of the semigroup required (see
condition (7.12) and Remark 1.7.5 in \citealp{EK1986}, p.39) easily
follows from Lemma \ref{description xi} and the finiteness of
$\Pn^s$. Furthermore, the existence of a limit for $r_D^{-1}G^{D}$
is obvious from the form of $G^D$. However, condition (7.17) of
Theorem 1.7.6 requires the existence of a subspace $E$ of functions
on $\Pn^s$ such that for every $f\in E$, there exist functions $g,
f_1,f_2,\ldots$ satisfying
$$\|f-f_D\|\rightarrow 0 \qquad \mathrm{and} \qquad\|g-G^Df_D\|\rightarrow 0 \qquad
\mathrm{as\ }D\rightarrow \infty.$$
The condition on $(f_D)_{D\geq 1}$ and the finiteness of $\Pn^s$ yield
$$G^Df_D = r_D\Psi f +o(r_D),$$
implying that a corresponding function $g$ can exist only if $\Psi
f\equiv 0$. Although $\Psi f(\zeta)=0$ if $\zeta \in \Pi_n$, this
condition would also require that $f(\zeta)=0$ whenever $\zeta
\notin \Pi_n$. Hence, to fit into Ethier and Kurtz' framework, an
obvious candidate for $E$ would be
$$E= \{f:f(\zeta)=0 \ \mathrm{for\ all\ }\zeta \in \Pn^s\setminus \Pi_n\},$$
where we then define a bounded linear transformation
$\wp_n:\Pn^s\rightarrow \Pi_n$ such that $\wp_n (\eta)=\eta$ for
every $\eta\in \Pi_n$. We may then apply Theorem 1.7.6 of
\citet{EK1986} and obtain convergence of the semigroup corresponding
to (or equivalently here of the finite dimensional distributions of)
$\wp_n(\mcP^D)$ to that of $\mcP$. However, it is unclear how to
define $\wp_n$ on the set $\Pn^s\setminus \Pi_n$, that is to specify
how to project $\Pn^s$ down onto its subset $\Pi_n$, in such a way
that the operator $\{\big(f\circ\wp_n, (\Gamma^n\circ p)
(f\circ\wp_n)\big),\ f:\Pn^s\rightarrow \Pn^s\}$ generates a Markov
process. Unfortunately, unless this condition is satisfied, Theorem
1.7.6 cannot be used to deduce the convergence result given in our
Theorem \ref{convergence p-dim}.

\subsection{Tightness}\label{section tightness}
The convergence of the finite-dimensional distributions relies on
the fact that the time required for the process to re-enter $\Pi_n$
following a geographical collision is vanishingly small as $D$ tends
to infinity. On the other hand, multiple changes to the
configuration of the genealogy can occur during this short period
with high probability, so that the conditions for $\mcP^D$ to
converge as processes in $D_{\Pn^s}([0,\infty))$ are much more
delicate.

Recall the definition of the stopping times $\sigma_i^D$ and
$\tau_i^D$ given in (\ref{defi sigma tau}). Suppose that the
probability that a $2$-event results in a configuration not in
$\Pi_n$ vanishes as $D\rightarrow 0$, or equivalently that
\begin{equation}\label{cond tightness P}\lim_{D\rightarrow
\infty}\mbP_{\zeta}[\tau_1^D<\infty,\ \mcP^D_{\tau_1^D}\notin
\Pi_2]=0,
\end{equation}
where $\zeta=\big(\big\{\{1\}\big\},\big\{\{2\}\big\}\big)$ and the
equivalence stems from the consistency equations (\ref{consist
lambda g}). Then, we easily see that the first time $\tau$ after
$\sigma_1^D$ such that $\mcP^D\notin \Pi_n$ converges to $+\infty$
in probability. Since
$$G^Df(\eta)=\Gamma f(\eta)+o(1) \qquad \mathrm{as\ }D\rightarrow
\infty$$ if $\eta \in \Pi_n$, we readily obtain that for any $a>0$,
the sequence of processes $(\{\mcP^D_t, t\geq a\},\ D\geq 1)$ is
tight (recall that $\sigma_1^D$ converges in probability towards
$0$). Let us now show that if condition (\ref{cond tightness P})
does not hold, the sequence $\mcP^D$ is not tight. It will be easier
to work with a metric on $\Pn^s$, the associated topology still
being the discrete topology.
\begin{prop}\label{prop non-tight}Assume that (\ref{cond tightness P}) does not hold.
Let $d$ be a discrete metric on $\Pn^s$, and suppose that $\zeta \in
\Pn^s$ is such that $\mbbP_{\und{\zeta}}[\sigma_2 <\infty]>0$. Then
the sequence of processes $\mcP^D$ under $\mbP_{\zeta}$ is not tight
in $D_{\Pn^s}([0,\infty))$ endowed with the Skorokhod topology
corresponding to $d$.
\end{prop}

\begin{proof}
First, recall the definition of the modulus of continuity $w'$ given
in \citet{EK1986}, p.122: for $X\in D_{\Pn^s}([0,\infty))$,
$\delta>0$ and $T>0$,
\begin{equation}
\label{defi w'}w'(X,\delta,T)\equiv \inf_{\{t_i\}}
\max_i\sup_{s,t\in [t_{i-1},t_i)}d(X_s,X_t),
\end{equation}
where the infimum is over all finite sets of times of the form
$0=t_0<t_1<\cdots<t_{n-1}<T\leq t_n$ such that
$\min_{1\leq i \leq n}(t_i-t_{i-1})>\delta$ and $n\geq 1$.

Suppose that the sequence $\mcP^D$ is tight. $\Pn^s$ is a finite
set, so the discrete topology on $(\Pn^s,d)$ turns it into a
complete and separable metric space, therefore $\mcP^D$ is also
relatively compact. By Corollary 3.7.4 of \citet{EK1986}, this
implies that for every $\gamma\in \Pn^s$, all $\eta>0$ and $T>0$,
there exists $\delta >0$ such that
\begin{equation} \label{non-tight}
\limsup_{D\rightarrow \infty} \mbP_{\gamma}\big[w'(\mcP^D,\delta,T)\geq
\eta\big]\leq \eta.
\end{equation}
Besides, the finiteness of $\Pn^s$ guarantees the existence of
$\epsilon >0$ such that, if $\gamma\neq \gamma' \in \Pn^s$, then
$d(\gamma, \gamma')>\epsilon$.

Let $T=1$, $\eta \in (0,\epsilon)$ and $\delta \in (0,1)$. We have
\begin{eqnarray}
\mbP_{\zeta}\big[w'(\mcP^D,\delta,1)\geq \eta \big] &=&
\mbP_{\zeta}\Big[w'(\mcP^D,\delta,1) \geq \eta \ \Big|\ \tau_1^D
\geq\frac{1}{2}\Big]\mbP_{\zeta}\Big[\tau_1^D \geq \frac{1}{2}\Big]
\nonumber\\ & & +\ \mbP_{\zeta}\Big[w'(\mcP^D,\delta,1)\geq \eta \
\Big|\ \tau_1^D <\frac{1}{2}\Big] \mbP_{\zeta}\Big[\tau_1^D <
\frac{1}{2}\Big]. \label{borne w}
\end{eqnarray}
On the one hand, $\sigma_1^D$ converges to $0$ in probability and
$\tau_1^D\Rightarrow \sigma_2$, so by Slutsky's lemma (see Lemma 2.8 in \citealp{vdV1998})
$\tau_1^D-\sigma_1^D\Rightarrow \sigma_2$, which is an exponential
random variable with positive parameter, so we have
\begin{eqnarray*}
\mbP_{\zeta}\Big[\tau_1^D < \frac{1}{2}\Big] &\geq&
\mbP_{\zeta}\Big[\sigma_1^D< \frac{1}{3}\mathrm{\
and\ }\tau_1^D-\sigma_1^D \leq \frac{1}{6}\Big] \\
&=& \mbP_{\zeta}\Big[\tau_1^D-\sigma_1^D \leq
\frac{1}{6}\Big]-\mbP_{\zeta}\Big[\sigma_1^D\geq \frac{1}{3}
\mathrm{\ and\ }\tau_1^D-\sigma_1^D \leq \frac{1}{6}\Big] \\
&\rightarrow & \mbbP_{\und{\zeta}}\Big[\sigma_2\leq \frac{1}{6}\Big] \equiv
C>0
\end{eqnarray*}
since the last term in the second line vanishes by the convergence
in probability of $\sigma_1^D$ to $0$. On the other hand,
\begin{eqnarray}
\mbP_{\zeta}\Big[w'(\mcP^D,\delta,1)\geq \eta
 \ \Big|\ \tau_1^D <\frac{1}{2}\Big] &\geq&  \mbP_{\zeta}\Big[w'(\mcP^D,\delta,1)\geq \eta
 \ \Big|\ \tau_1^D <\frac{1}{2},\ \sigma_2^D-\tau_1^D <\frac{\delta}{2},\ \mcP^D_{\tau_1^D}\notin \Pi_n \Big]
 \nonumber
 \\ & & \qquad \times \mbP_{\zeta}\Big[\sigma_2^D-
 \tau_1^D <\frac{\delta}{2},\ \mcP^D_{\tau_1^D}\notin \Pi_n\Big|\ \tau^D_1 <\frac{1}{2}\Big]. \label{module cont}
 \end{eqnarray}
By the convergence in probability of $\sigma_{1}^D$ to $0$,
uniformly in $\eta \in \Pn^s$, and the strong Markov property
applied to $\mcP^D$ at time $\tau_1^D$, we obtain that
$\sigma_2^D-\tau_1^D$ converges in probability to $0$. Furthermore,
on the event that no $R_D-$events occurred between the times
$\sigma_1^D$ and $\tau_1^D$ (the probability of which is growing to
one), $\tau_1^D$ is the epoch of the first event after $\sigma_1^D$
and $\mcP_{\tau_1^D}^D$ its outcome so, by the strong Markov
property, $\tau_1^D$ and $\mcP_{\tau_1^D}^D$ are independent
conditionally on $\mathcal{F}^D_{\sigma_1^D}$. Since (\ref{cond
tightness P}) does not hold, we can write
$$\liminf_{D\rightarrow \infty} \mbP_{\zeta}\Big[\sigma_2^D-
 \tau_1^D <\frac{\delta}{2},\ \mcP^D_{\tau_1^D}\notin \Pi_n\Big|\ \tau^D_1 <\frac{1}{2}\Big] \geq
 \liminf_{D\rightarrow \infty}\mbP_{\zeta}\big[\mcP^D_{\tau_1^D}\notin \Pi_n\big] >0.$$
 Now, if $\tau_1^D <\frac{1}{2}$, $\mcP^D_{\tau_1^D}\notin \Pi_n$ and
$\sigma_2^D-\tau_1^D <\frac{\delta}{2}$, then by definition of
$\epsilon$,
$$d\big(\mcP^D_{\tau^D_1-},\mcP^D_{\tau^D_1}\big)>\epsilon
\quad\mathrm{\ and\ } \quad d\big(\mcP^D_{\tau^D_1},
\mcP^D_{\sigma^D_2}\big)>\epsilon,$$ and by assumption
$\sigma_2^D-\tau_1^D <\frac{\delta}{2}$, so $w'(\mcP^D,\delta,1)\geq
\epsilon \geq \eta.$ Consequently,
$$\liminf_{D\rightarrow \infty}\mbP_{\zeta}\Big[w'(\mcP^D,\delta,1)\geq \eta
 \ \Big|\ \tau_1^D <\frac{1}{2}\Big] >0. $$
Therefore, we see from (\ref{borne w}) that, for all $\delta \in (0,1)$,
$$\liminf_{D\rightarrow \infty}\mbP_{\zeta}\big[w'(\mcP^D,\delta,1)\geq \eta \big] > C' > \eta $$
for any $\eta \in (0,\epsilon \wedge C')$. This yields a contradiction with (\ref{non-tight}).
\end{proof}
>From the last proof, we see that what prevents the sequence of
structured genealogical processes from being tight is that at each
geographical collision, at least two jumps accumulate: the
geographical collision itself and one or more transitions generated
by $\xi$ to bring $\mcP^D$ back into $\Pi_n$.

Yet the unstructured genealogical process, which is not a Markov
process for $D< \infty$, is not modified by movements of blocks.
Thus, if the number of jumps needed by $\mcP^D$ to re-enter $\Pi_n$
after a geographical collision were at most one with a probability
growing to $1$, we would expect tightness to hold for $\mcP^{D,u}$
(recall that $\zeta^u$ denotes the unstructured partition generated
by $\zeta$). The next proposition in fact gives an equivalence
between the behaviour of the latter probability and tightness of
$\big\{\mcP^{D,u}, D\geq 1\big\}$.
\begin{prop}\label{prop tightness}
For each $D\in \N$, let $U_1^D$ be the random time defined by
$$U^D_1\equiv \inf\big\{t>0:\ \mcP^D_{t-}\neq \mcP^D_t\big\},$$
with the convention that $\inf \emptyset =+\infty$. Note that, if
$\mcP^D_0 \notin \Pi_n$, then $U^D_1\leq \sigma_1^D$. Let also $\chi(\Pi_n)$
denote the image of $\Pi_n$ by the first
geographical collision (when it exists), that is
$$\chi(\Pi_n)\equiv \big\{\gamma\in \Pn^s\ :\ \exists\ \zeta\in \Pi_n,\ \mbbP[\chi(\zeta)=\gamma] >0\big\}.$$
Suppose that for all $\gamma\notin \Pi_n$, $\mbP[\und{\gamma}^u\neq
\gamma^u]>0$ (meaning that the process $\xi$ started at
$\gamma$ has at least one coalescence with positive
probability).

Then the following are equivalent:

(i) For all $\gamma \in \chi(\Pi_n)\setminus \Pi_n$,
$\lim_{D\rightarrow \infty}\mbP_{\gamma}[U^D_1=\sigma_1^D] = 1$.

(ii) For all $\zeta \in \Pn^s$ and $a>0$, the sequence of
$D_{\Pn}([a,\infty))$-valued random variables $\mcP^{D,u}$, started
at $\zeta^u$ at time $0$, is tight.

Furthermore, if $\zeta \in \Pi_n\cup \chi(\Pi_n)$, then condition
(i) is equivalent to the tightness in $D_{\Pn}([0,\infty))$ of
$\mcP^{D,u}$ started at $\zeta^u$.

As a consequence of Theorem \ref{convergence p-dim}, if conditions
(i) and (ii) hold, then for all $\zeta \in \Pn^s$ and $a>0$, the law
of $\big(\mcP^{D,u}_t,\ t\geq a\big)$ under $\mbP_{\zeta}$ converges
to the law under $\mbbP_{\und{\zeta}^u}$ of $\big(\mcP^{u}_t,\ t\geq
a\big)$. Furthermore, if $\zeta \in \Pi_n$, then the convergence
holds for $a =0$.
\end{prop}

\begin{rema} Assuming that for all $\gamma\notin \Pi_n$, $\mbP[\und{\gamma}^u\neq
\gamma^u]>0$ is actually not required, but not
supposing it makes the proof unnecessarily more involved.
\end{rema}

\begin{proof} Once again we work with a metric $d$ on $\Pn$, so that
$D_{\Pn}([0,\infty))$ is a complete and separable metric space and
the sequence $(\mcP^{D,u})_{D\geq 1}$ is tight if and only if it is
relatively compact. We call $\epsilon>0$ the minimum distance
between two different partitions. Let us first show that if
condition $(i)$ is not fulfilled, then neither is condition $(ii)$.
The following proof is highly reminescent to the proof of
Proposition \ref{prop non-tight}, so let us adopt directly the same
notation. In particular, we work with $T=1$ and $\zeta$ such that
$\mbbP_{\und{\zeta}}[\sigma_2 < \infty]
>0$.

For each $a>0$, let us write $w'_a$ the modulus of continuity of a
process corresponding to times $t\geq a$, defined as in (\ref{defi
w'}) with the condition on the finite sets $\{t_i\}$ replaced by
$a=t_0 <\cdots < T\leq t_n$. Fix $a\in (0,1/3)$, and let $\eta \in
(0,\epsilon)$ and $\delta\in (0,1)$. The same calculation as in the
proof of Proposition \ref{prop non-tight} holds by replacing the
event $\{w'(\mcP^D,\delta,1)\geq \eta\}$ by
$\{w'_a(\mcP^{D,u},\delta,1)\geq \eta\}$ and $\mbP_{\zeta}[\tau_1^D<
1/2]$ by $\mbP_{\zeta}[1/3\leq \tau_1^D < 1/2]$. Hence, by
(\ref{module cont}) and the argument directly following it, we just
need to prove that
\begin{equation}\label{key non-tightness}\mbP_{\zeta}\Big[w'_a(\mcP^{D,u},\delta,1)\geq \eta
 \ \Big|\ \frac{1}{3}\leq \tau_1^D <\frac{1}{2},\ \sigma_2^D-\tau_1^D <\frac{\delta}{2},\
 \mcP^D_{\tau_1^D}\notin \Pi_n \Big] \end{equation}
is bounded below by a positive constant for $D$ large enough.
If we define $V^D_1$ by
$$V^D_1\equiv \inf \big\{t>\tau_1^D\ :\ \mcP^D_{t-}\neq \mcP^D_t\big\},$$
then the expression in (\ref{key non-tightness}) is equal to
\setlength\arraycolsep{1pt}
\begin{eqnarray}\mbP_{\zeta}&\Big[&w'_a(\mcP^{D,u},\delta,1)\geq \eta\ ;\ \mcP^{D,u}_{\tau^D_1}
=\mcP^{D,u}_{V^D_1}\mathrm{\ or\
}\mcP^{D,u}_{V^D_1}=\mcP^{D,u}_{\sigma^D_2}
 \ \Big|\ \frac{1}{3}\leq \tau_1^D <\frac{1}{2},\ \sigma_2^D-\tau_1^D <\frac{\delta}{2},\
 \mcP^D_{\tau_1^D}\notin \Pi_n \Big] \label{eq non-tight}\\
& & + \mbP_{\zeta}\Big[w'_a(\mcP^{D,u},\delta,1)\geq \eta\ ;\
\mcP^{D,u}_{\tau^D_1}\neq\mcP^{D,u}_{V^D_1}\ ;\
\mcP^{D,u}_{V^D_1}\neq \mcP^{D,u}_{\sigma^D_2}\Big|\ \frac{1}{3}\leq
\tau_1^D <\frac{1}{2},\ \sigma_2^D-\tau_1^D <\frac{\delta}{2},\
\mcP^D_{\tau_1^D}\notin \Pi_n \Big] \nonumber
\end{eqnarray}
The first term in (\ref{eq non-tight}) is nonnegative, and if we are
in the conditions given by the second term, then
$$d(\mcP^{D,u}_{\tau_1^D},\mcP^{D,u}_{V^D_1})\geq \epsilon,\quad d(\mcP^{D,u}_{V_1^D},
\mcP^{D,u}_{\sigma_2^D}) \geq \epsilon\quad \mathrm{and}\quad
\sigma_2^D-\tau_1^D <\frac{\delta}{2},$$ implying that
$w'_a(\mcP^D,\delta,1)\geq \epsilon > \eta$. Therefore, the second
term in (\ref{eq non-tight}) is equal to \setlength\arraycolsep{1pt}
\begin{eqnarray*}
\mbP_{\zeta}\Big[\mcP^{D,u}_{\tau^D_1}\neq\mcP^{D,u}_{V^D_1}&\ ;\
&\mcP^{D,u}_{V^D_1}\neq \mcP^{D,u}_{\sigma^D_2}\Big|\
\frac{1}{3}\leq \tau_1^D <\frac{1}{2},\ \sigma_2^D-\tau_1^D
<\frac{\delta}{2},\
\mcP^D_{\tau_1^D}\notin \Pi_n \Big]\\
&=&\mbP_{\zeta}\Big[\mcP^{D,u}_{\tau^D_1}\neq\mcP^{D,u}_{V^D_1}\ ;\
\mcP^{D,u}_{V^D_1}\neq \mcP^{D,u}_{\sigma^D_2}\Big|\ \frac{1}{3}\leq
\tau_1^D <\frac{1}{2},\ \mcP^D_{\tau_1^D}\notin
\Pi_n\Big]\big(1+o(1)\big).
\end{eqnarray*}
Now, by the strong Markov property applied to $\mcP^D$ at time
$\tau_1^D$, we have \setlength\arraycolsep{1pt}
\begin{eqnarray}\mbP_{\zeta}\Big[\mcP^{D,u}_{\tau^D_1}\neq\mcP^{D,u}_{V^D_1}\ &;&
\mcP^{D,u}_{V^D_1}\neq \mcP^{D,u}_{\sigma^D_2}\ ;\
\frac{1}{3}\leq \tau_1^D <\frac{1}{2};\ \mcP^D_{\tau_1^D}\notin \Pi_n\Big] \label{U/V} \\
& & =
\mbE_{\zeta}\Big[\mbP_{\mcP^D_{\tau_1^D}}\big[\mcP^{D,u}_{\tau^D_1}\neq
\tilde{\mcP}^{D,u}_{\tilde{U}^D_1}\ ;\
\tilde{\mcP}^{D,u}_{\tilde{U}^D_1}\neq
\tilde{\mcP}^{D,u}_{\tilde{\sigma}^D_1},\ \big] \ \ind_{\{1/3\leq
\tau_1^D < 1/2\}}\ind_{\{\mcP^D_{\tau_1^D}\notin \Pi_n\}}\Big].
\nonumber
\end{eqnarray}
Since we assumed that condition $(i)$ did not hold, there exists
$\eta \in \chi(\Pi_n)\setminus \Pi_n$ such that
$\mbP_{\eta}[\tilde{U}^D_1<\tilde{\sigma}_1^D] \geq C_1$ for a
constant $C_1>0$ and $D$ large enough. As $\eta \in \chi(\Pi_n)$, we
can choose $\zeta$ such that $\mbbP[\chi(\und{\zeta})=\eta]>0$ (and
$\mbbP_{\und{\zeta}}[\sigma_2 < \infty]
>0$). Now, since we assumed
that $\mbP[\und{\gamma}^u\neq \gamma^u]>0$ for all $\gamma \notin
\Pi_n$, the probability that a coalescence event occurs before a
scattering event in the structured genealogical process $\xi$ started
at any value not in $\Pi_n$ is greater than a constant $C_2$.
Therefore, we can write
$$\mbP_{\eta}\big[\eta^u \neq \tilde{\mcP}^{D,u}_{\tilde{U}^D_1}\ ;\
\tilde{\mcP}^{D,u}_{\tilde{U}^D_1}\neq
\tilde{\mcP}^{D,u}_{\tilde{\sigma}^D_1}\big]> C_1'$$ for a constant
$C_1'>0$. By the distribution of the epochs of the geographical
collisions, the convergence in law of $(\tau_1^D,\mcP^D_{\tau_1^D})$
to $\big(\sigma_2,\chi(\und{\zeta})\big)$ (c.f. Lemma \ref{conv tau
sigma}) and the fact that $\eta \notin \Pi_n$, we have for $\zeta$
chosen as above
$$\mbP_{\zeta}\Big[\mcP^D_{\tau_1^D}=\eta\ ;\ \frac{1}{3}\leq \tau_1^D <\frac{1}{2};\
\mcP^D_{\tau_1^D}\notin \Pi_n\Big]>C_3$$ for a constant $C_3>0$ and
$D$ large enough, so the expression in the right-hand side of
(\ref{U/V}) is bounded below by $C_1'C_3$, and so is (\ref{key
non-tightness}). Hence, $(ii)\Rightarrow (i)$.

Suppose now that condition $(i)$ is fulfilled. Condition (a) of
Corollary 3.7.4 in \citet{EK1986} trivially holds, so we
only need to check condition (b) on the modulus of continuity. Fix
$\zeta \in \Pn^s$ and $a>0$, and let $T>a$ and $\eta>0$. Firstly, by
the convergence in probability of $\sigma_1^D$ to $0$, there exists
$D_1\in \N$ such that for all $D\geq D_1$, $\mbP_{\zeta}[\sigma_1^D
\geq a]< \frac{\eta}{5}$. Secondly, we have
$$\mbP_{\zeta}\big[\mathrm{at\ least\ one\ }R_D\mathrm{-event\ in\ }[0,T]\big]\leq
1-\exp\big(-T\max_{\xi\in \Pn^s}c_{R_D}(\xi)\big)\rightarrow
0,\qquad D\rightarrow \infty,$$ so there exists $D_2\geq 1$ such
that for all $D\geq D_2$, the previous quantity is less than
$\frac{\eta}{5}$. Thirdly, by the same argument as in the beginning
of the proof of Theorem \ref{convergence p-dim}, there exists $N\in
\N$ such that $\mbbP_{\und{\zeta}}[\sigma_N \leq T]<\frac{\eta}{5}$.
Hence, by Lemma \ref{conv tau sigma}, there exists $D_3 \geq 1$ such
that for all $D\geq D_3$, $\mbP_{\zeta}[\sigma^D_N \leq
T]<\frac{\eta}{5}$.

Consequently, we can write for each $D\geq \max \{D_1,D_2, D_3\}$ and all $\delta>0$
\begin{eqnarray*}
\mbP_{\zeta}[w'_a(\mcP^{D,u},\delta,T)\geq \eta]&\leq &
\mbP_{\zeta}[\sigma_1^D \geq a] + \mbP_{\zeta}[\sigma^D_N \leq T]
+\mbP_{\zeta}\big[\mathrm{at\ least\ one\ }R_D\mathrm{-event\ in\ }
[0,T]\big]\\
& & \qquad + \mbP_{\zeta}\big[w'_a(\mcP^{D,u},\delta,T)\geq \eta;\
\sigma_1^D <a;\ \sigma^D_N > T;\
\mathrm{no\ }R_D\mathrm{-events\ in\ }[0,T] \big] \\
&< & \frac{3\eta}{5} +
\mbP_{\zeta}\big[w'_a(\mcP^{D,u},\delta,T)\geq \eta;\ \sigma_1^D
<a;\ \sigma^D_N > T;\ \mathrm{no\ }R_D\mathrm{-events\ in\ }[0,T]
\big].
\end{eqnarray*}
Furthermore, there exists $\delta>0$ such that
$\mbbP_{\gamma}[\sigma_2<3\delta]<\frac{\eta}{5N}$ for all $\gamma
\in \Pi_n$. Now, for all $i\in \{1,\ldots,N\}$, by the strong Markov
property applied to $\mcP^D$ at time $\tau_{i-1}^D$ and the
convergence of $\mbP_{\gamma}[\tau^D_1<3\delta]$ to
$\mbbP_{\und{\gamma}} [\sigma_2<3\delta]$, uniformly in $\gamma$, we
have
$$\mbP_{\zeta}[\tau_{i-1}^D< \infty;\ \tau_{i}^D-\tau_{i-1}^D<3\delta] = \mbE_{\zeta}\big[
\ind_{\{\tau_{i-1}^D< \infty\}}\
\mbP_{\mcP^D_{\tau^D_{i-1}}}[\tilde{\tau}_1^D <3\delta]\big] \leq
\frac{\eta}{5N}$$ for $D$ large enough. Therefore,
\setlength\arraycolsep{1pt}
\begin{eqnarray}
\mbP_{\zeta}&\big[&w'_a(\mcP^{D,u},\delta,T)\geq \eta;\ \sigma_1^D
<a;\ \sigma^D_N \leq T;\ \mathrm{no\ }
R_D\mathrm{-event\ in\ }[0,T] \big] \nonumber\\
&\leq& \sum_{i=1}^N \mbP_{\zeta}[\tau_{i-1}^D< \infty;\ \tau_{i}^D-\tau_{i-1}^D<3\delta] \nonumber\\
& & + \mbP_{\zeta}\big[w'_a(\mcP^{D,u},\delta,T)\geq \eta;\
\sigma_1^D <a;\ \sigma^D_N \leq T;\ \mathrm{no\ }R_D\mathrm{-event\
in\ }[0,T];\ \tau_{i}^D-\tau_{i-1}^D\geq 3\delta \mathrm{\ for\ all\
}i\leq N \nonumber\\ & & \qquad \qquad \qquad\qquad \mathrm{\ s.t.\
}\tau_{i-1}^D< \infty \big] \label{tau's loin}
\end{eqnarray}
and the first sum is less than $\frac{\eta}{5}$. To finish, let $V_i^D$
denote the epoch of the next event after $\tau_i^D$ if
$\mcP^D_{\tau_i^D}\notin \Pi_n$ (if it exists, $V_i^D=+\infty$
otherwise), and set $V_i^D=\tau_i^D=\sigma_{i+1}^D$ if
$\mcP^D_{\tau_i^D}\in \Pi_n$. Since we assume that condition $(i)$
holds, for all $i\in \{1,\ldots,N\}$ we have by the strong Markov
property applied at time $\tau_i^D$ and the fact that the
distribution of $\mcP^D_{\tau_{i}^D}$ concentrates on $\chi(\Pi_n)$
as $D$ grows to infinity by Lemma \ref{conv tau sigma},
$$\mbP_{\zeta}[\tau_i^D<\infty; V_i^D<\sigma_{i+1}^D ] \rightarrow 0,\qquad D\rightarrow \infty,$$
so the last term in (\ref{tau's loin}) is less than
\setlength\arraycolsep{1pt}
\begin{eqnarray*}
\sum_{i=1}^{N} &\mbP_{\zeta}&[\tau_i^D<\infty; V_i^D<\sigma_{i+1}^D] \\
& & + \mbP_{\zeta}\big[w'_a(\mcP^{D,u},\delta,T)\geq \eta;\
\sigma_1^D <a;\ \sigma^D_N > T;\ \mathrm{no\ } R_D\mathrm{-event\
in\ }[0,T];\ \tau_{i}^D-\tau_{i-1}^D\geq 3\delta \mathrm{\ and}\\ &
& \qquad \qquad \qquad\qquad V_i^D=\sigma_{i+1}^D\mathrm{\ for\ all\
}i\leq N\mathrm{\ s.t.\ }\tau_{i-1}^D< \infty  \big],
\end{eqnarray*}
where the first sum is less than $\frac{\eta}{5}$ for $D$ large
enough. But on that last event, $\sigma_1^D$ is less than $a$ and no
$R_D$-events occur so $\tau_1^D$ is the epoch of the event directly
after $\sigma_1^D$, then all geographical collisions are at least
$3\delta$ away from each other and the $\sigma_i^D$'s are the only
times in between at which an event occurs, so necessarily
$w'_a(\mcP^D,\delta,T)=0$. Assembling all the pieces, we obtain that
$$\mbP_{\zeta}[w'_a(\mcP^{D,u},\delta,T)\geq \eta]
<\eta, $$ completing the proof of $(i)\Rightarrow (ii)$.

If $\zeta \in \chi(\Pi_n)\cup \Pi_n$, then we only need to show that
$(i)$ implies the tightness of $(\mcP^{D,u})_{D\geq 1}$ on
$[0,\infty)$. Let us directly use the same notation as in the last
proof. In the last paragraph, we proved that with a high
probability, there is no accumulations of jumps between the random
time $\tau_1^D$ and $T$. Also, we can make
$\mbP_{\zeta}[\tau_1^D\leq 2a]$ as small as we want by adjusting $a$
and taking $D$ large enough, and the probability that at least one
$R_D$-event occurs is vanishingly small, so we are left with proving
that, if  $\delta$ is such that
$\mbP_{\zeta}[w'_a(\mcP^{D,u},\delta,T)\geq \eta] <\eta$, $\tau_1^D>
2a$ and no $R_D$-events occur between $0$ and $T$, then
$\mbP_{\zeta}[w'(\mcP^{D,u},\delta',a)\geq \eta] <\eta, $ for some
$\delta' \in (0,\delta)$. If $\zeta\in \Pi_n$, $\tau_1^D> 2a$ and no
$R_D$-events occur, then $\tau_1^D$ is the epoch of the first event
occurring to $\mcP^D$ so $w'(\mcP^{D,u},\delta',a)=0$ for all
$\delta' \in (0,\delta \wedge a)$. If $\zeta\in \chi(\Pi_n)$ and the
other conditions hold, then by condition (i) we have
$$\mbP_{\zeta}[U^D_1=\sigma_1^D] \rightarrow 1, \qquad D\rightarrow \infty, $$
and furthermore $\mbP_{\zeta}[\sigma_1^D<a]\rightarrow 1$, so with a
probability tending to one as $D$ grows to infinity, one event
occurs between $0$ and $a$, then nothing happens between $a$ and
$2a$ (there is no $R_D$-events, so the next event after $\sigma_1^D$
must occur at time $\tau_1^D>2a$) and the condition on the modulus
of continuity is fulfilled after time $a$ so, for any $\delta' \in
(0,\delta \wedge a)$, we do have
$$\mbP_{\zeta}[w'(\mcP^{D,u},\delta',a)\geq \eta] <\eta.$$
This completes the proof of the case  $\zeta \in \chi(\Pi_n)\cup \Pi_n$.

Now, by Theorem 3.7.8 in \citet{EK1986}, the two ingredients to
obtain the convergence of the processes $(\mcP^D)_{D\geq 1}$ are
tightness, given by the first part of Proposition \ref{prop
tightness} for any $a>0$, and convergence of the finite-dimensional
distributions, given by Theorem \ref{convergence p-dim} and the
bijective correspondence between $\Pi_n$ and $\Pn$. For $\zeta \in
\Pi_n$ and $a=0$, tightness still holds by virtue of the last
paragraph, and an easy modification (namely allowing $t=0$ in the
proof of the convergence of the one-dimensional distributions) of
the proof of Theorem \ref{convergence p-dim} in that case, where
$\und{\zeta} =\zeta$ and $\sigma_1^D=0$ a.s., gives the convergence
of the finite-dimensional distributions of $\mcP^D$, including at
time $t=0$.
\end{proof}

Let us briefly comment on the condition
$\mbP_{\xi}[U^D_1=\sigma_1^D] \rightarrow 1$. If the fast
within-deme coalescence is given by a $\Xi$-coalescent (including
Kingman's coalescent) occurring in one deme at a time, the condition
is fulfilled if and only if at most two lineages can be collected
into the same deme during a single event. Indeed, in that case the
next step of the genealogical process is either to scatter these two
lineages into two different demes or to merge them into one lineage,
the outcome of which is always in $\Pi_n$. If more than $2$ lineages
are gathered in the same deme and do not merge during the
geographical collision, then with a positive probability only two of
them are involved in the next genealogical event and at least two
rapid steps are needed for $\mcP^D$ to re-enter $\Pi_n$. The same
conclusion holds if two pairs of lineages are gathered in two demes
(meaning $2$ lineages per deme), since the genealogical process acts
in one deme at a time by assumption.

\section{Collapse of structured genealogical processes}\label{section collapse}
The next proposition states that the only reasonable structured
genealogies which collapse to an unstructured genealogy (given by a
$\Xi$-coalescent) when the number of demes tends to infinity are the
genealogies that we have described before, subject to certain
conditions.

Note that if we want the lineages to be exchangeable in the limit,
the limiting process needs to take its values in $\bigcup_{n\geq
1}\Pi_n$. Indeed, since the rates of intra- and inter-demes mergers
greatly differ, we should observe only inter-demes events on the
slow time scale. This requires that each deme contains at most one
lineage at any given time in the limit.
\begin{prop}\label{collapse}Let $(\mcP^D_t,\ t\geq 0)_{D\geq 1}$ be a sequence of
structured genealogical processes with values in $\bigcup_{n\geq 1}
\Pn^s$. Then the following are equivalent
\begin{enumerate}
\item There exists a sequence $r_D$ such that $r_D\rightarrow +\infty$
as $D\rightarrow \infty$ and two structured genealogical processes,
$(\xi_t,t\geq 0)$ (resp. $(\mcP_t,t\geq 0)$) with values in
$\bigcup_{n\geq 1}\Pn^s$ (resp. $\bigcup_{n\geq 1}\Pi_n $)
satisfying
\begin{enumerate}
\item for each $n\in \N$, the sequence of structured genealogical processes
$(\mcP^D_{r_D^{-1}t}, t\geq 0)_{D\geq 1}$ on the fast time scale,
with initial value in $\Pn^s$, converges to $\xi$ as a process in
$D_{\Pn^s}[0,\infty)$. In addition, $\xi$ is consistent in the sense
of Lemma \ref{consistence};
\item the sequence $(\mcP^D_t,t\geq 1)_{D\geq 1}$ on the slow time scale
converges towards $\mcP$ in that the finite-dimensional
distributions (except possibly at time $0$) converge as in Theorem
\ref{convergence p-dim} for every sample size $n$;
\item there exists a $\Xi$-coalescent $(R_t,t\geq 1)$ such that for all $n\geq 1$, the
unstructured genealogical process $\mcP^u$ induced by
$\mcP|_{\Pi_n}$ has the law of the restriction of $R$ to $\Pn$.
\end{enumerate}
\item The rates associated to $\mcP^D$ satisfy conditions (\ref{ass1}),
(\ref{ass2}) and (\ref{ass4}) of Section 3.1, and condition $(i)$ of
Lemma \ref{consistence} holds.
\end{enumerate}
\end{prop}

We shall see in the proof that the consistency of $\xi$ is a key
ingredient to obtain the desired equivalence. In fact, if we did not
impose it, it would certainly be possible to construct particular
examples in which the unstructured genealogy on the slow time scale
is also a $\Xi$-coalescent, but the genealogies within a deme are
not consistent. We would need to impose `good' values for the corresponding
rates. We rather chose here to emphasize more biologically relevant
models, for which the within-deme genealogical process is also
consistent and which can be described as part of an entire class of
models rather than special cases.
\begin{proof}
The implication $2\Rightarrow 1$ in a consequence of Theorem
\ref{convergence p-dim}, Proposition \ref{convergence fast
time scale} and Proposition \ref{xi-coal}.

Let us prove that $1\Rightarrow 2$. From the definition of a
structured genealogical process, blocks can only move and coalesce.
Furthermore $\mcP^D$ stays in $\Pn^s$ whenever its initial value
lies in this set, so we just need to fix $n\geq 0$ and look at the
corresponding rates of scattering, gathering and coalescence. From
the description of the limiting processes $\xi$ and $\mcP$, the
rates of $\mcP^D$ must be of the form
$$r_D\rho_D^{(1)}(\eta|\zeta) + \rho_D^{(2)}(\eta|\zeta)+o(1),$$
where for $i\in \{1,2\}$, $\rho_D^{(i)}(\eta|\zeta) \rightarrow
\rho^{(i)}(\eta|\zeta)$ as $D$ tends to infinity. (To simplify
notation, we shall write $\rho_D^{(i)}(\eta|\zeta)\equiv
\rho^{(i)}(\eta|\zeta)$.) Thus, $\rho^{(1)}(\eta|\zeta)$ are the
rates associated to the generator $\Psi$ of the process $\xi$. Let
us check that all cited conditions necessarily hold:
\begin{itemize}
\item If $\zeta \rightarrow \eta$ is a $1$-event, then by adding an $(n+1)$-st
individual in one of the existing blocks (therefore changing the
sizes of the blocks but not their number), we see that the
consistency of $\xi$ imposes that the part of the rate corresponding
to the fast time scale  depends neither on $n$, nor on the sizes of
the blocks. By exchangeability of the demes, this rate is thus
characterized by the number of lineages present in each deme before
and after the transition, the order of these numbers being
irrelevant. Therefore, condition (1) holds. By Lemma
\ref{consistence}, the consistency of $\xi$ implies that condition
(i) of Lemma \ref{consistence} is also satisfied.

\item Once again by consistency of $\xi$, the rate of a $2$-event
must be of order $1$. Indeed, it may otherwise lead to an additional
$1$-event for the restriction of the process with the $(n+1)$-st
lineage (if this additional lineage lands in a non-empty deme or in
the same deme as another moving lineage coming from a different
subpopulation , and the other dispersing lineages land in different
demes), or involve at least two lineages alone in their demes on the
fast time scale. If such an event was allowed, then by
exchangeability of the islands the fast dynamic could act on a
structured partition in $\Pi_n$ and merge two lineages starting from
different demes. Again by exchangeability, any pair of lineages
could merge on the fast time scale  and so the outcome of $\xi$ would
be a single lineage with probability one, a trivial situation which
is of no interest here. Now, since we want to keep exchangeability
of the lineages in the unstructured genealogy (on the slow
time scale), the rates of $2$-events should depend only on the number
of lineages and their geographical distribution (and possibly on
$n$). But if $\zeta \in \Pi_n$, all lineages are in different demes,
so the corresponding rates are necessarily of the form given in
condition (2). If the rates were to depend on $n$, then as the rates
of the fast genealogical process which follows directly (for $D$
large enough, as in the proof of Theorem \ref{convergence p-dim})
are independent of $n$, the overall transition from $\eta \in \Pi_n$
to the value of $\mcP^D$ when it reenters $\Pi_n$ would eventually
give different rates for $\mcP$ acting on $\Pi_n$ and for the
restriction to $\Pi_n$ of $\mcP$ acting on $\Pi_{n+1}$ (recall the
convergence of $\tau_{i-1}^D$ and $\sigma_i^D$ towards $\sigma_i$ to
see that the transitions of $\mcP$ actually can be described as in
Section 3.2). This would contradict the fact that the process
$\mcP^{u}$ corresponds to a $\Xi$-coalescent. Finally, we obtain
that condition (2) must hold.

\item The last argument imposes also that geographical collisions
involving $k$ lineages occur at a rate which is the sum of all
corresponding geographical events involving $k+1$ lineages, which is
exactly writing the consistency equations (\ref{consist lambda g})
of condition (3).
\end{itemize}
Finally, we obtain that $2\Rightarrow 1$.
\end{proof}

\section{Example}\label{section example}
We now turn our attention to a particular class of metapopulation models which 
combine a (finite) $\Lambda$-coalescent within demes with migration between 
demes and sporadic mass extinction events.  We will use the results derived in the 
preceding sections to characterize the form that the genealogy takes in the infinitely
many demes limit.  This, in turn, will allow us to illustrate how the statistics of the 
population-wide $\Xi$-coalescent depend on the interplay between extinction/recolonization 
events and the local demographic processes occurring within demes.  While these 
models are quite contrived - in particular, we have simply imposed the condition that 
a small number of demes is responsible for repopulating vacant demes following a 
mass extinction - they will allow us to explicitly calculate some quantities of interest.

We first describe how the population evolves forwards-in-time.  Suppose that for each 
$D$, each deme contains exactly $N$ individuals.  Fix $K \in \N$, and let $\Lambda^d(dx)$ 
and $\Lambda^g(dy)$ be two probability measures on $[0,1]$ with no atom at $0$.
Then reproduction, migration, and extinction/recolonization events occur according
to the following rules.

\begin{itemize}
\item Each individual in each deme reproduces at rate $D$ according to
the following scheme.  If an individual in deme $i$ reproduces, then a 
number $x$ is sampled from $[0,1]$ according to the probability distribution
$\Lambda^{d}(dx)$, and then each occupant of that deme dies with 
probability $x$ and is replaced by an offspring of the reproducing individual.  
In terms of the notation of Section \ref{section cannings}, such an event
has the following representation when $k$ is the label of the reproducing
individual.  First, $R^{j,j'}=(0,\ldots,0)$ for all pairs of integers $j \neq j'\in [D]$
and $R^{j,j}=(1,\ldots,1)$ if $j\in [D]\setminus \{i\}$.  $R^{i,i}$ is a random 
vector obtained by choosing a number $x$ according to $\Lambda^d(dx)$, 
a number $m$ according to a binomial distribution with parameters $(N,x)$, 
and finally a set $\mathcal{O}\subset [N]$ of offspring of the reproducing 
individual by sampling $m$ labels in $[N]$ uniformly without replacement. 
Then, $R^{i,i}_k=m$, $R^{i,i}_{k'}=0$ for all $k'\in \mathcal{O}\setminus \{k\}$, 
and $R^{i,i}_l=1$ for all $l\notin \mathcal{O}\cup \{k\}$.

\item At rate $Dm_1$, each individual gives birth to a single migrant offspring 
which then moves to any one of the $D$ demes, chosen uniformly at random, 
and replaces one of the $N$ individuals within that deme, also uniformly at 
random. In this case, if $j$ is the label of the deme containing the parent and 
$k$ is its label, then a pair $(i,l)$ is sampled uniformly in $[D]\times [N]$ and
the vectors $R$ are as described in Example \ref{ex2} of Section
\ref{section cannings}.

\item Mass extinction events occur at rate $e$.  When such an event occurs,
a number $y$ is sampled from $[0,1]$ according to the probability distribution 
$\Lambda^g(dy)$.  Then, each deme goes extinct with probability $y$, independently 
of all the others, and is unaffected by the extinction otherwise.  Simultaneously, 
$K$ of the $D$ demes are chosen uniformly at random to be source demes, and 
the deceased occupants of the extinct island are replaced by offspring produced 
by individuals living in the source demes according to the following scheme.  The
parent of each individual recolonizing a deme left vacant by the mass extinction is
chosen uniformly at random and with replacement from among the $NK$ inhabitants 
of the source demes.  If a source deme is chosen from among the extinct ones, 
then the parents of the offspring emerging from that deme are the individuals that
occupied the deme immediately prior to the extinction.  To describe such an event 
using the notation of Section \ref{section cannings}, suppose that a number $y$ is
chosen according to $\Lambda^g(dy)$, a number $m$ is sampled according to a 
$\mathrm{Binom}(D,y)$-distribution and a (random) set $\mathcal{O}_{\mathrm{ext}}\subset [D]$ 
is constructed by sampling uniformly without replacement $m$ deme labels.  Independently,
another set $\mathcal{O}_{\mathrm{rec}}$ of $K$ recolonizing demes is also chosen 
by uniform sampling.  Then, for all $i\notin \mathcal{O}_{\mathrm{ext}}$ we have
$R^{i,i}=(1,\ldots,1)$ and each deme $j\in \mathcal{O}_{\mathrm{ext}}\setminus 
\mathcal{O}_{\mathrm{rec}}$ satisfies $R^{i,j}=(0,\ldots,0)$ for all $i\in [D]$. The vectors
$R^{i,j}$ with $j\in \mathcal{O}_{\mathrm{rec}}$ and $i\in \mathcal{O}_{\mathrm{ext}} 
\cup \mathcal{O}_{\mathrm{rec}}$ are not easily formulated explicitly (in particular, 
their description depends on whether the recolonizing demes also go extinct during 
the event), but it is clear that the evolution of the population satisfies the two conditions 
required in Section \ref{section cannings}.
\end{itemize}

Suppose that $n$ individuals are sampled from the population at time $0$, and let us 
consider the evolution (backwards-in-time) of the structured coalescent process $\mcP^D$ 
in $\Pn^s$.  From the description of the model forwards-in-time, the events affecting the
genealogy occur at the following rates:

\begin{enumerate}
\item If a deme contains $b$ lineages, then each $k$-tuple of lineages in this deme (for 
$k\leq b$) merges into one lineage in the same deme at rate $$D\lambda^d_{b;k,1,\ldots,1}
= DN\int_0^1 \Lambda^d(dx) x^{k}(1-x)^{b-k}.$$ Furthermore, any merger event occurs 
in one deme at a time.

\item Each lineage migrates (alone) at rate $Dm_1$.  Indeed, the total rate at which migrant 
offspring are produced forwards-in-time is $ND \times Dm_1$, but the probability that such
a migrant belongs to the lineage under consideration is $(ND)^{-1}$ (recall that the deme 
and the label of the individual replaced by the migrant are chosen uniformly at random).  
Consequently, the probability that a migrating lineage lands in a non-empty deme is $D^{-1}$ 
times the number of demes occupied by the other lineages of $\mcP^D_{t-}$.  When such an 
event occurs, the probability that the migrating lineage also merges with an ancestral lineage 
present in the source deme is $N^{-1}$ times the number of distinct ancestral lineages present 
in that deme.

\item Extinction events generate geographical collisions at rate $O(1)$.  Because the $K$ 
recolonizing demes are chosen uniformly from among the $D$ islands, recolonization by a deme 
containing at least one lineage of the genealogical process occurs with a probability of order 
$O(D^{-1})$, and so these events are negligible in the limit.  Suppose that $\mcP^D_{t-}\in 
\Pi_n$.  Let $k\leq |\mcP^D_{t-}|$, $r\leq K$, and let $k_1,\ldots,k_r$ be integers greater than 
$1$ and summing to $k$.  For each $i\in \{1,\ldots,r\}$, let $L_i=\{l_{i1},\ldots,l_{ij_i}\}$ be a 
collection of $j_i$ integers summing to $k_i$.  Then each $(|\mcP^D_{t-}|;k_1,\ldots,k_r,
1,\ldots,1;$ $L_1,\ldots,L_r,\{1\},\ldots,\{1\})$-geographical collision occurs at rate 

\setlength\arraycolsep{1pt}
\begin{eqnarray} &e&\int_0^1
\Lambda^g(dy)\sum_{s=0}^{|\mcP^D_{t-}|-k}\ind_{\{s \leq K-r\}}
\binom{|\mcP^D_{t-}|-k}{s}y^{k+s}(1-y)^{|\mcP^D_{t-}|-k-s}
\frac{K!}{(K-r-s)!}\frac{1}{K^{k+s}}\nonumber
\\
& & \qquad \qquad \times
\prod_{i=1}^r\left\{\frac{N!}{(N-j_i)!}\frac{1}{N^{k_i}}\right\}
+O\Big(\frac{1}{D}\Big)\label{rates example} \\ &\equiv &
\tilde{\lambda}^g_{|\mcP^D_{t-}|;k_1,\ldots,k_r,1,\ldots,1;L_1,\ldots,L_r,\{1\},\ldots,\{1\}}
+O\Big(\frac{1}{D}\Big).\nonumber\end{eqnarray}  
\end{enumerate}

The rate expression that appears in Eq.\ (35) can be interpreted in the following way.  As well 
as the $k$ ancestral lineages that are known to be affected by the disturbance (this is specified
by the type of event), an additional $s$ lineages may be caught up in the extinction event and 
moved to demes where they remain isolated (hence producing no changes in the structured 
genealogy).  In (\ref{rates example}), the first part in each term of the sum corresponds to the 
number of choices for these $s$ additional lineages, followed by the probability that only these 
$k+s$ lineages are affected.  The condition $r+s \leq K$ is imposed by the fact that the $r$ groups 
of lineages geographically gathered and the $s$ lineages affected but remaining alone in their 
demes must then belong to $r+s$ distinct recolonizing demes. The middle part of the term specifies 
the probability that the affected lineages are grouped in the desired way: regardless of the labels 
of the recolonizing demes, if the latter contain no lineages of the sample just before the extinction 
then $\frac{K!}{(K-r-s)!}$ is the number of (unordered) ways of choosing $r+s$ of them to receive
the affected lineages, while $K^{-k-s}$ is the probability that each of the $k+s$ lineages moves
to the prescribed recolonizing deme. Finally, the last product is obtained in a similar manner by
allocating as many distinct ancestors as required to the groups of lineages gathered into the same 
demes.  As explained above, the $O(D^{-1})$ remainder term accounts for the probability that at 
least one of the finitely-many recolonizing demes contains a lineage of $\mcP^D_{t-}$.

Let us say that a \textbf{simple collision} occurs when a single lineage moves into a non-empty 
deme, and possibly merges with one of the lineages present in this deme.  To verify that the 
convergence results from the previous sections apply to the example, it will be convenient
to introduce the following quantities, defined for all $\zeta,\eta \in \Pn^s$:

$$\phi_c(\zeta,\eta)= \left\{\begin{array}{ll}
1& \mathrm{if\ }\zeta \rightarrow \eta \mathrm{\ is\ a\ simple\
collision\ with\ coalescence,}\\
0& \mathrm{otherwise}, \end{array}\right.$$ and likewise
$$\phi_{nc}(\zeta,\eta)= \left\{\begin{array}{ll}
1& \mathrm{if\ }\zeta \rightarrow \eta \mathrm{\ is\ a\ simple\
collision\ without\ coalescence,}\\
0& \mathrm{otherwise}. \end{array}\right.$$ 
By `with coalescence' (resp. `without coalescence'), we mean that the migrating lineage merges
(resp. does not merge) with a lineage in the source deme during the same event.

Let us consider a particular $1$-event.  If this event involves a single lineage moving to an empty 
deme, it may be caused either by a migration event of the kind described in item 2 above (which 
occurs at rate $Dm_1(1-k/D)$ if $k$ is the number of demes occupied by the other lineages at the 
time of the event), or by a mass extinction event (whose rate is of order $O(1)$ according to item 3). Consequently, the overall rate of any $1$-event is of the form $Dm_1+O(1)$.  Groups of more than 
one lineage can also move simultaneously, but only through an extinction event and so at a rate of 
order $O(1)$.  If the event involves an intra-deme merger, then its rate is easily written in the form 
given in Section \ref{section conditions} with $r_D=D$; see item 1.  A $2$-event $\zeta \rightarrow \eta$ 
occurs at a rate of order $O(1)$, and in particular if $\zeta \in \Pi_n$, then this rate is given by
$$\tilde{\lambda}^g(\zeta,\eta)+ 2m_1\Big\{\phi_c(\zeta,\eta)\frac{1}{N} +
\phi_{nc}(\zeta,\eta)\frac{N-1}{N}\Big\} + O\Big(\frac{1}{D}\Big)\equiv
\lambda^g(\zeta,\eta)+O\Big(\frac{1}{D}\Big),$$ 
where $\tilde{\lambda}^g(\zeta,\eta)$ is the rate of the unique extinction event which turns $\zeta$ 
into $\eta$.  In this expression, the term in brackets is nonzero only if the event is a simple
collision involving two lineages that have been collected in the same deme through migration.
Such collisions occur at rate $2 m_{1}$, and then the two lineages either coalesce, with probability 
$N^{-1}$, or remain distinct, with probability $1-N^{-1}$.  Finally, we must check that the 
$\lambda^g$'s satisfy (\ref{consist lambda g}), and that the rates on the fast time scale satisfy 
condition $(i)$ of Lemma \ref{consistence}.  The latter condition follows from the description of 
the rates and the consistency of $\Lambda$-coalescents, and the validity of the former condition 
can be deduced from the fact that lineages choose independently of each other whether they are 
involved in the event, and which of the $NK$ individual they take as a parent.  We leave the 
straightforward but tedious calculations to the interested reader.

All conditions of Theorem \ref{convergence p-dim} and Proposition \ref{xi-coal} hold.  Thus, we 
can conclude that the finite dimensional distributions of $\mcP^D$ converge to those of a structured genealogical process $\mcP$ with values in $\Pi_n$, and that the unstructured process $\mcP^u$ 
is a $\Xi$-coalescent with values in $\Pn$.  Let us describe $\mcP^u$ as precisely as we can. 
To apply the results of Section \ref{descript P}, we need to know the distribution of the final state 
of the `fast' process $\xi$ that was introduced in Section \ref{paragraphe xi}.  Starting from a structured 
partition where all blocks are contained in the same component (i.e., all lineages lie initially in the 
same deme), this distribution coincides with the sampling distribution of the infinitely many alleles 
model of the generalized Fleming-Viot process dual to the $\Lambda$-coalescent with finite 
measure $x^2\Lambda^d(dx)$ acting within this deme.  Indeed, on the fast time scale, ancestral 
lineages belonging to a common deme migrate out to distinct, empty islands, a process analogous 
to mutation to unique types with a `mutation' rate equal to $m_1$.  Recursion formulae are given
in \citet{MOH2006} which can be used to compute the probability $p(\mathbf{n})$ of unordered 
allele configurations $\mathbf{n}=\{n_1,\ldots,n_k\}$ in the infinitely many alleles model when the 
genealogy is given by a $\Lambda$- or a $\Xi$-coalescent.  In our case, the formula of interest is 
(with $p(1)=1$):
$$p(\mathbf{n})=
\frac{nm_1}{g_n+nm_1}\sum_{j=1}^k\ind_{\{n_j=1\}}\frac{1}{k}\ p(\tilde{\mathbf{n}}_j)
+ \sum_{i=1}^{n-1}\frac{g_{n,n-i}}{g_n+nm_1}\sum_{j=1}^k
\ind_{\{n_j>i\}}\frac{n_j-i}{n-i}\ p(\mathbf{n}-i\mathbf{e}_j),$$
where $n\equiv \sum_{j=1}^kn_j \geq 2$, $\tilde{\mathbf{n}}_j \equiv 
(n_1,\ldots,n_{j-1},n_{j+1},\ldots,n_k)$, $\mathbf{e}_j$ denotes the $j$'th unit vector in 
$\mathbf{R}^k$ and $g_{nk}$ (resp. $g_n$) is the rate at which the number of lineages decreases 
from $n$ to $k$ (resp. the total rate at which the number of lineages changes when $n$ lineages 
are alive), given by $$g_{nk}=\binom{n}{k-1}\int_0^1\Lambda^d(dx)x^{n-k+1}(1-x)^{k-1}$$
and $$g_n=\sum_{k=1}^{n-1}g_{nk}=\int_0^1\Lambda^d(dx)\big(1-(1-x)^{n-1}(1-x+nx)\big).$$
These expressions are related to the distribution of $\und{\zeta}$ by the following formula:
$$\mbP\big[\und{\zeta}=\big(\{B_1\},\ldots,\{B_k\},\emptyset,\ldots,\emptyset \big)\big]
=p(|B_1|,\ldots,|B_k|),$$ where $\zeta= \big(\big\{\{1\},\ldots,\{n\}\big\},\emptyset,\ldots,\emptyset\big)$
and $|B_i|$ denotes the number of elements in the block $B_i$.  Indeed, because the dynamics 
on the fast time scale of lineages occupying different demes are independent, the final state of 
the fast genealogical process is the concatenation of all the final states of the groups of lineages 
starting in the same deme.  Hence, the preceding results are sufficient to describe $\und{\zeta}$ 
for any $\zeta \in \Pn^s$.  Unfortunately, with this level of generality, there does not appear to be 
a simple description of the measure $\Xi$ associated to $\mcP^u$, but the rate associated to its 
Kingman part (that is its mass at $\mathbf{0}$) is given by:
\begin{equation}\label{kingman part}
2m_1 \frac{1}{N}+ 2m_1
\frac{N-1}{N}\ p(2) = 2\frac{m_1}{N}
\left\{1+ (N-1)\ \frac{\int_0^1\Lambda^d(dx)x^2}{\int_0^1\Lambda^d(dx)x^2+2m_1}
\right\}.
\end{equation}
The first term in (\ref{kingman part}) corresponds to a simple collision with coalescence, and the 
second term to a simple collision without coalescence; the probability that the lineages then coalesce 
before one of them migrates is given by $p(2)$.

One case which can be characterized more thoroughly is when dispersal between demes only 
occurs during extinction-recolonization events ($m_{1} = 0$).  For example, this might be a 
reasonable approximation to make when modeling a population in which migrants are at a 
substantial competitive disadvantage relative to residents, so that dispersal is only successful 
into demes in which the resident population has gone extinct.  In this case, the Kingman component 
of the genealogy disappears (see (\ref{kingman part})).  Furthermore, viewed backwards in time, 
lineages gathered into common demes by mass extinction events cannot migrate away before 
the rapid within-deme coalescent reaches a common ancestor, and so any such group of lineages 
merges instantaneously into a single lineage.  The shape of the resulting global coalescent therefore 
is determined only by the way in which mass extinction events gather lineages together.  Recall 
the expression for the rates of geographical collisions given in (\ref{rates example}), and let us 
examine how $K$, the number of demes contributing colonists in the wake of a mass extinction, 
affects the shape of the genealogy.

If $K=1$, all lineages affected by a mass extinction event have parents within the same deme.  
The resulting genealogy is a $\Lambda$-coalescent, and the rate at which $k$ ancestral lineages 
merge when $m$ are present is equal to the rate at which exactly $k$ lineages are caught up in 
an extinction event when $m$ demes contain one lineage, that is
$$e\int_0^1 \Lambda^g(dy)y^k(1-y)^{m-k}.$$
On the other hand, if we let $K$ tend to infinity, then each term in the sum in (\ref{rates example}) 
is asymptotically equivalent to $\frac{K!}{(K-r-s)!}\ K^{-k-s}\sim K^{r-k}$, up to a constant (recall that 
the sample size $n$ is finite and bounds the number of lineages at any times).  Consequently, 
binary geographical collisions ($k=k_1=2$, $r=1$, $j_1\in \{1,2\}$) occur at a rate of order
$O(K^{-1})$, whereas the rate of a collision involving at least $3$ lineages is of order at most 
$O(K^{-2})$. Hence, for fixed sample size $n$, the probability that only binary mergers occur in 
the sample genealogy approaches $1$ as $K$ tends to infinity, and the rate of each binary merger 
(multiplied by $K$) converges to
\begin{equation}
	\label{K large} e\int_0^1 \Lambda^g(dy)y^2,
\end{equation}
where the term $y^2$ is obtained by observing that the condition $s\leq K-1$ in (\ref{rates example}) 
is always fulfilled for $n$ fixed and $K$ large enough, and that $\sum_{s=0}^{|\mcP_{t-}^D|-2}
\binom{|\mcP_{t-}^D|-2}{s}y^{2+s}(1-y)^{|\mcP_{t-}^D|-2-s}=y^2$.  Once the lineages are gathered 
into the same deme, they can only coalesce and they do so instantaneously on the slow time 
scale as $D\rightarrow \infty$.  It follows that if time is rescaled by a factor of $DK$, then the rate 
of a binary merger converges to that of Kingman's coalescent run at the rate shown in (\ref{K large}).   
Moreover, under this time rescaling, the rates of the finitely many possible multiple merger events 
converge to $0$ as $K$ grows to infinity, and so the limiting (as $D\rightarrow\infty$) unstructured 
genealogical process $\mcP^u$ corresponding to an evolution with $K$ recolonizing demes
converges to Kingman's coalescent as a process in $D_{\Pn}[0,\infty)$ as $K$ tends to infinity. 
(Note, however, that this does not imply that one can interchange the limits $D\rightarrow \infty$ 
and $K\rightarrow \infty$.)  Finally, if $K$ is finite but greater than $1$, then geographical collisions 
involving more than two lineages occur at a non-negligible rate, and so the resulting unstructured 
genealogy is a more general $\Xi$-coalescent.  

This example shows that a large class of coalescent processes can arise in the infinitely many 
demes limit of a subdivided population with sporadic mass extinctions.  Depending on both the 
migration and the extinction rates, as well as on the number of demes contributing to population 
recovery following a mass extinction, the limiting genealogical process can range from Kingman's 
coalescent ($K=\infty$), as derived by \citet{WAK2004}, to a $\Lambda$-coalescent  ($K=1,\ m_1=0$), 
with a family of $\Xi$-coalescents interpolating between these two extremes.  In this particular class 
of models, multiple mergers of ancestral lineages are more likely to occur when all three parameters, 
$K$, $N$ and $m_1$, are small, so that mass extinctions have a non-negligible probability of 
gathering lineages into a common deme where they undergo a series of rapid mergers before being 
scattered again by migration.  This observation suggests that it is a generic property of structured 
population models that if the limiting coalescent admits any multiple mergers, then it also admits 
simultaneous mergers.

\bigskip
\noindent \textbf{Acknowledgements.} We are grateful to Alison
Etheridge for helpful discussions and comments on this work, 
and to the referees for suggestions that have improved the 
presentation of the paper.  A. V\'eber would like to thank the 
Department of Statistics of the University of Oxford for hospitality.

\bibliographystyle{plainnat}

\begin{thebibliography}{1}
\bibitem[Cannings(1974)]{CAN1974}C. Cannings. The Latent Roots of Certain Markov
Chains Arising in Genetics: A New Approach, I. Haploid Models.
\emph{Adv. Appl. Prob.}, 6:260--290, 1974.
\bibitem[Cox(1989)]{COX1989}J.T. Cox. Coalescing random walks and
voter model consensus times on the torus in $\mathbb{Z}^d$.
\emph{Ann. Probab.}, 17:1333--1366, 1989.
\bibitem[Cox and Durrett(2002)]{CD2002}J.T. Cox and R. Durrett.
The stepping stone model: New formulas expose old myths. \emph{Ann. Appl. Probab.}, 12:1348--1377, 2002.
\bibitem[Eldon and Wakeley(2006)]{EW2006}B. Eldon and J. Wakeley. Coalescent
Processes When the Distribution of Offspring Number Among
Individuals is Highly Skewed. \emph{Genetics}, 172:2621--2633, 2006.
\bibitem[Ethier and Kurtz(1986)]{EK1986}S.N. Ethier and T.G. Kurtz. \emph{Markov processes:
characterization and convergence}. Wiley, 1986.
\bibitem[Ethier and Nagylaki(1980)]{EN1980}S.N. Ethier and T. Nagylaki. Diffusion
approximation of Markov chains with two time scales and applications
to population genetics. \emph{Adv. Appl. Probab.}, 12:14--49, 1980.
\bibitem[Fisher(1930)]{FIS1930}R. Fisher. \emph{The Genetical Theory of Natural Selection}.
Clarenson, Oxford, 1930.
\bibitem[Greven et al.(2007)]{GLW2007}A. Greven, V. Limic and A.
Winter. Coalescent processes arising in a study of diffusive
clustering. \emph{Preprint}, 2007.
\bibitem[Kingman(1982)]{KIN1982}J.F.C. Kingman. The coalescent. \emph{Stoch. Proc. Appl.}, 13:235--248, 1982.
\bibitem[M\"ohle(2000)]{MOH2000}M. M\"ohle. Ancestral Processes in Population Genetics - the Coalescent.
\emph{J. Theoret. Biol.}, 204:629--638, 2000.
\bibitem[M\"ohle(2006)]{MOH2006}M. M\"ohle. On sampling distributions for coalescent
processes with simultaneous multiple collisions. \emph{Bernoulli},
12:35--53, 2006.
\bibitem[M\"ohle and Sagitov(2001)]{MS2001}M. M\"ohle and S. Sagitov. A Classification
of Coalescent Processes for Haploid Exchangeable Population Models.
\emph{Ann. Probab.}, 29:1547--1562, 2001.
\bibitem[Nordborg(2001)]{NOR2001}M. Nordborg. \emph{Coalescent theory}, Chapter 7 in \emph{Handbook of
Statistical Genetics}. Wiley, Chichester, UK, 2001.
\bibitem[Nordborg and Krone(2002)]{NK2002}M. Nordborg and S.M. Krone. \emph{Separation of time scales and
convergence to the coalescent in structured populations}. In
\emph{Modern Developments in Theoretical Population Genetics: The
Legacy of Gustave Mal\'ecot}. Oxford University Press, Oxford, 2002.
\bibitem[Notohara(1990)]{NOT1990}M. Notohara. The coalescent and the genealogical process
in geographically structured populations. \emph{J. Math. Biol.},
29:59--75, 1990.
\bibitem[Pitman(1999)]{PIT1999}J. Pitman. Coalescents with Multiple Collisions.
\emph{Ann. Probab.}, 27:1870--1902, 1999.
\bibitem[Sagitov(1999)]{SAG1999}S. Sagitov. The general coalescent with asynchronous mergers
of ancestral lines. \emph{J. Appl. Probab.}, 36:1116--1125, 1999.
\bibitem[Sargsyan and Wakeley(2008)]{SW2008}O. Sargsyan and J. Wakeley. A coalescent process with simultaneous multiple mergers for approximating the gene genealogies of many marine organisms. \emph{Theor. Popul. Biol.}, 74:104--114, 2008.
\bibitem[Schweinsberg(2000)]{SCH2000}J. Schweinsberg. Coalescents with simultaneous multiple collisions.
\emph{Electr. J. Probab.}, 5:1--50, 2000.
\bibitem[Schweinsberg(2003)]{SCH2003}J. Schweinsberg. Coalescent processes obtained from
supercritical Galton-Watson processes. \emph{Stochastic Process.
Appl.}, 106:107--139, 2003.
\bibitem[Sousa(1984)]{SOU1984}W.P. Sousa. The role of disturbance in natural communities.
\emph{Annual Review of Ecology and Systematics}, 15:353--391, 1984.
\bibitem[van der Vaart(1998)]{vdV1998}A.W. van der Vaart. \emph{Asymptotic Statistics}.
Cambridge University Press, 1998.
\bibitem[Wakeley(1998)]{WAK1998}J. Wakeley. Segregating sites in Wright's island model. \emph{Theor. Popul. Biol.},
53:166--175, 1998.
\bibitem[Wakeley(1999)]{WAK1999}J. Wakeley. Nonequilibrium Migration in Human History. \emph{Genetics},
153:1863--1871, 1999.
\bibitem[Wakeley(2004)]{WAK2004}J. Wakeley. Metapopulation models for historical inference. \emph{Mol. Ecol.},
13:865--875, 2004.
\bibitem[Wakeley and Aliacar(2001)]{WA2001}J. Wakeley and N. Aliacar. Gene Genealogies in a Metapopulation.
\emph{Genetics}, 159:893--905, 2001.
\bibitem[Wilkinson-Herbots(1998)]{WIL1998}H.M. Wilkinson-Herbots. Genealogy and subpopulation differentiation
under various models of population structure. \emph{J. Math. Biol.}, 37:535--585, 1998.
\bibitem[Wright(1931)]{WRI1931}S. Wright. Evolution in Mendelian populations. \emph{Genetics}, 16:97--159, 1931.
\bibitem[Z\"ahle et al.(2005)]{ZCD2005}I. Z\"ahle, J.T. Cox and R. Durrett.
The stepping stone model II:Genealogies and the infinite sites model. \emph{Ann. Appl. Probab.}, 15:671--699, 2005.
\end{thebibliography}

\end{document}